\pdfoutput=1

\documentclass[preprint,10pt]{elsarticle}

\usepackage{fullpage}
\usepackage{amsmath,amssymb,amsfonts}
\usepackage[titletoc,toc,title]{appendix}

\usepackage{array} 
\usepackage{listings}
\usepackage{mathtools}
\usepackage{pdfpages}
\usepackage[textsize=footnotesize,color=green]{todonotes}
\usepackage{bm}
\usepackage[normalem]{ulem}
\usepackage{hhline}

\usepackage{algorithm}
\usepackage[noend]{algpseudocode}
\usepackage{algorithmicx}
\algblock{ParFor}{EndParFor}
\algnewcommand\algorithmicparfor{\textbf{parfor}}
\algnewcommand\algorithmicpardo{\textbf{do}}
\algnewcommand\algorithmicendparfor{\textbf{end\ parfor}}
\algrenewtext{ParFor}[1]{\algorithmicparfor\ #1\ \algorithmicpardo}
\algrenewtext{EndParFor}{\algorithmicendparfor}
\usepackage{graphicx}
\usepackage{subfig}
\usepackage{color}

\definecolor{forestgreen}{rgb}{0, 0.5, 0}

\usepackage{pgfplots}
\usepackage{pgfplotstable}
\definecolor{markercolor}{RGB}{124.9, 255, 160.65}
\pgfplotsset{width=10cm}
\pgfplotsset{
tick label style={font=\small},
label style={font=\small},
legend style={font=\small}
}

\usetikzlibrary{calc}

\newcommand{\logLogSlopeTriangle}[5]
{
    \pgfplotsextra
    {
        \pgfkeysgetvalue{/pgfplots/xmin}{\xmin}
        \pgfkeysgetvalue{/pgfplots/xmax}{\xmax}
        \pgfkeysgetvalue{/pgfplots/ymin}{\ymin}
        \pgfkeysgetvalue{/pgfplots/ymax}{\ymax}

        \pgfmathsetmacro{\xArel}{#1}
        \pgfmathsetmacro{\yArel}{#3}
        \pgfmathsetmacro{\xBrel}{#1-#2}
        \pgfmathsetmacro{\yBrel}{\yArel}
        \pgfmathsetmacro{\xCrel}{\xArel}

        \pgfmathsetmacro{\lnxB}{\xmin*(1-(#1-#2))+\xmax*(#1-#2)} % in [xmin,xmax].
        \pgfmathsetmacro{\lnxA}{\xmin*(1-#1)+\xmax*#1} % in [xmin,xmax].
        \pgfmathsetmacro{\lnyA}{\ymin*(1-#3)+\ymax*#3} % in [ymin,ymax].
        \pgfmathsetmacro{\lnyC}{\lnyA+#4*(\lnxA-\lnxB)}
        \pgfmathsetmacro{\yCrel}{\lnyC-\ymin)/(\ymax-\ymin)} % THE IMPROVED EXPRESSION WITHOUT 'DIMENSION TOO LARGE' ERROR.

        \coordinate (A) at (rel axis cs:\xArel,\yArel);
        \coordinate (B) at (rel axis cs:\xBrel,\yBrel);
        \coordinate (C) at (rel axis cs:\xCrel,\yCrel);

        \draw[#5]   (A)-- node[pos=0.5,anchor=north] {1}
                    (B)-- 
                    (C)-- node[pos=0.5,anchor=west] {#4}
                    cycle;
    }
}
\newcommand{\logLogSlopeTriangleFlip}[5]
{
    \pgfplotsextra
    {
        \pgfkeysgetvalue{/pgfplots/xmin}{\xmin}
        \pgfkeysgetvalue{/pgfplots/xmax}{\xmax}
        \pgfkeysgetvalue{/pgfplots/ymin}{\ymin}
        \pgfkeysgetvalue{/pgfplots/ymax}{\ymax}

        \pgfmathsetmacro{\xBrel}{#1-#2}
        \pgfmathsetmacro{\yBrel}{#3}
        \pgfmathsetmacro{\xCrel}{#1}

        \pgfmathsetmacro{\lnxB}{\xmin*(1-(#1-#2))+\xmax*(#1-#2)} % in [xmin,xmax].
        \pgfmathsetmacro{\lnxA}{\xmin*(1-#1)+\xmax*#1} % in [xmin,xmax].
        \pgfmathsetmacro{\lnyA}{\ymin*(1-#3)+\ymax*#3} % in [ymin,ymax].
        \pgfmathsetmacro{\lnyC}{\lnyA+#4*(\lnxA-\lnxB)}
        \pgfmathsetmacro{\yCrel}{\lnyC-\ymin)/(\ymax-\ymin)} % THE IMPROVED EXPRESSION WITHOUT 'DIMENSION TOO LARGE' ERROR.

	\pgfmathsetmacro{\xArel}{\xBrel}
        \pgfmathsetmacro{\yArel}{\yCrel}

        \coordinate (A) at (rel axis cs:\xArel,\yArel);
        \coordinate (B) at (rel axis cs:\xBrel,\yBrel);
        \coordinate (C) at (rel axis cs:\xCrel,\yCrel);

        \draw[#5]   (A)-- node[pos=0.5,anchor=east] {#4}
                    (B)-- 
                    (C)-- node[pos=0.5,anchor=south] {1}
                    cycle;
    }
}
\usepackage{stmaryrd}

\newcommand{\td}[2]{\frac{{\rm d}#1}{{\rm d}{\rm #2}}}
\newcommand{\pd}[2]{\frac{\partial#1}{\partial#2}}
\newcommand{\pdd}[2]{\frac{\partial^2#1}{\partial#2^2}}
\newcommand{\pdn}[3]{\frac{\partial^{#3}#1}{\partial#2^{#3}}}

\newcommand{\nor}[1]{\left\| #1 \right\|}
\newcommand{\LRp}[1]{\left( #1 \right)}

\newcommand{\LRa}[1]{\left\langle #1 \right\rangle}
\newcommand{\LRb}[1]{\left| #1 \right|}
\newcommand{\LRc}[1]{\left\{ #1 \right\}}

\newcommand{\Grad} {\ensuremath{\nabla}}
\newcommand{\Div} {\ensuremath{\nabla\cdot}}
\newcommand{\jump}[1] {\ensuremath{\llbracket#1\rrbracket}}
\newcommand{\avg}[1] {\ensuremath{\LRc{\!\{#1\}\!}}}

\renewcommand{\L}{L^2\LRp{\Omega}}

\newcommand{\Dhat}{\widehat{D}}

\newcommand{\eval}[2][\right]{\relax
  \ifx#1\right\relax \left.\fi#2#1\rvert}

\newcommand{\reviewerOne}[1]{{{#1}}}
\newcommand{\reviewerTwo}[1]{{{#1}}}

\newcolumntype{C}[1]{>{\centering\let\newline\\\arraybackslash\hspace{0pt}}m{#1}}

\newcommand*\diff[1]{\mathop{}\!{\mathrm{d}#1}}

\makeatletter
\renewcommand\d[1]{\mspace{6mu}\mathrm{d}#1\@ifnextchar\d{\mspace{-3mu}}{}}
\makeatother

\date{}
\begin{document}

%\maketitle
%\tableofcontents
\begin{frontmatter}
\author[rice]{Jesse Chan\corref{cor1}}
\ead{Jesse.Chan@caam.rice.edu}
\cortext[cor1]{Principal Corresponding author}
\author[vt]{John A.\ Evans}
\ead{John.A.Evans@colorado.edu}
\address[rice]{Department of Computational and Applied Mathematics, Rice University, 6100 Main St, Houston, TX, 77005}
\address[vt]{Ann and H.J. Smead Aerospace Engineering Sciences, University of Colorado Boulder, 429 UCB, Boulder, Colorado, 80309}

%\title{Multi-patch discontinuous Galerkin spline finite element methods for time-domain wave propagation}
\title{Multi-patch discontinuous Galerkin isogeometric analysis for wave propagation: explicit time-stepping and efficient mass matrix inversion}

\begin{abstract}
We present a class of spline finite element methods for time-domain wave propagation which are particularly amenable to explicit time-stepping.  The proposed methods utilize a discontinuous Galerkin discretization to enforce continuity of the solution field across geometric patches in a multi-patch setting, which yields a mass matrix with convenient block diagonal structure.  Over each patch, we show how to accurately and efficiently invert mass matrices in the presence of curved geometries by using a weight-adjusted approximation of the mass matrix inverse.  This approximation restores a tensor product structure while retaining provable high order accuracy and semi-discrete energy stability.  We also estimate the maximum stable timestep for spline-based finite elements and show that the use of spline spaces result in less stringent CFL restrictions than equivalent $C^0$ or discontinuous finite element spaces.  Finally, we explore the use of optimal knot vectors based on $L^2$ $n$-widths.  We show how the use of optimal knot vectors can improve both approximation properties and the maximum stable timestep, and present a simple heuristic method for approximating optimal knot positions.  Numerical experiments confirm the accuracy and stability of the proposed methods.  
\end{abstract}
\end{frontmatter}
\section{Introduction}

The key concept behind isogeometric analysis is the integration of computer representations of geometry in design and numerical simulations.  
%computer representations of geometry (i.e.\ CAD) and numerical simulations.  
This integration aims to address the engineering bottleneck in converting engineering designs into a form appropriate for analysis \cite{hughes2005isogeometric} by {representing solutions of PDEs using the spline representations underlying geometric parametrizations}.  While the main advantage of isogeometric analysis is the elimination of errors in the approximation of curved geometries, it has been observed that discretizations using high order B-spline and NURBS possess advantages over traditional high order $C^0$ finite element methods, including greater efficiency per degree of freedom, efficient high order accurate approximations of eigenvalues and eigenfunctions of continuous differential operators, and {larger maximum stable time-step sizes for explicit dynamics} \cite{hughes2008duality,evans2009n, hughes2014finite}.  

In \cite{hughes2008duality}, NURBS-based discretizations were compared to traditional $C^0$ finite element methods for steady-state wave propagation problems such as time-harmonic wave propagation and eigenvalue problems in structural vibrations.  The authors noted that, compared to $C^0$ finite element discretizations, B-spline and NURBS discretizations approximated high frequency modes more accurately relative to the number of degrees of freedom.  These properties of spline and NURBS were also shown to translate into improved numerical simulations of one-dimensional time-dependent wave propagation in \cite{hughes2014finite}.  However, {the application of NURBS-based finite element discretizations to higher dimensional time-dependent wave propagation problems faces several computational challenges, such as the cost of inverting spline mass matrices and dealing with multi-patch geometries}.
%\footnote{Spline collocation methods have been applied to higher dimensional time-dependent wave propagation problems \cite{auricchio2012isogeometric}; however, these are not equivalent in general to finite element discretiations.}  
The goal of this work is to make the advances of isogeometric analysis simpler to realize through the use of multi-patch discontinuous Galerkin discretizations of time-dependent wave propagation problems in two and three dimensions, focusing in particular on explicit time integration.

We consider geometries which are constructed using multiple geometric ``patches'', which are mappings of a reference patch to physical space.  While a single NURBS patch can exactly represent conic sections, \reviewerTwo{geo-metrically} and topologically complex objects require representations in terms of multiple NURBS patches.  We discretize these geometries using a multi-patch discontinuous Galerkin method \cite{langer2015multipatch}, which imposes continuity conditions weakly across patch interfaces through a numerical flux \cite{nguyen2014nitsche}.  These methods are analogous to multi-domain spectral methods \cite{canuto2012spectral}, which combine the accuracy of spectral methods with the use of multiple domains (patches) for geometric flexibility.  %We study multi-patch DG discretizations of model problems in wave propagation in both first and second order form.  

When paired with explicit time integration, multi-patch DG methods require the inversion of a \reviewerTwo{sparse but large} mass matrix over each patch at every timestep.   
%Many isogeometric methods for time-dependent problems have traditionally utilized either implicit time integration or space-time formulations \cite{cottrell2009isogeometric}.  
%The combination of NURBS with explicit time-stepping is uncommon due to the fact that, 
For geometric patches which are affine transformations of reference patches, the tensor product structure of NURBS bases result in mass matrices which can be inverted by applying inverses of one-dimensional matrices.  However, for curved geometries of practical engineering interest, the Kronecker structure of the mass matrix inverse is lost.  We address the issue of mass matrix inversion by using a weight-adjusted approximation to the mass matrix, which restores a tensor product structure to the approximate mass inverse while maintaining provable high order accuracy and energy stability.  %These properties address existing approaches to efficient explicit time-integration, such as isogeometric collocation and preconditioned iterative methods \cite{auricchio2012isogeometric, gao2014fast}.  

This work also studies timestep restrictions for explicit methods under high order NURBS finite element discretizations.  High order $C^0$ and DG finite element discretizations possess maximum stable timesteps which scale as $O(h/p^2)$, where $p$ is the order of approximation { and $h$ is the mesh size} \cite{warburton2008taming}.  By using special filters or discretizations involving dual grids, the maximum stable timestep can be made to scale as $O(h/p)$ instead \cite{warburton2008taming, reyna2014operator}, though these methods can be difficult to generalize to unstructured meshes.  In this work, we derive estimates for the maximum stable timestep based on a modification of arguments in \cite{chan2015gpu}.  These estimates depend on constants in $hp$ inverse and trace inequalities.  
%{ \sout{Multi-patch spline discretizations introduce two meshes: a macro-element mesh consisting of geometric ``patches'' of size $H$, and a micro-element mesh containing $K$ micro-elements { per-side} over which a spline approximation is defined.  While explicit estimates for these constants are not currently available except in special cases, numerical experiments show that these constants (and the resulting spectral radius of the NURBS discretization matrix) grow as $O(hH/p)$ for a micro-element mesh size $h$ which is sufficient small with respect to $p$.} 
While explicit values for these constants are not currently available except in special cases, numerical experiments show that the maximum stable timestep for multi-patch discontinuous Galerkin discretizations scales like $O(h/p)$ provided the number of elements per patch is sufficiently large with respect to $p$.\footnote{We note that the slow growth of the maximum stable timestep with $p$ appears to be correlated with the use of extended-support basis functions.  For example, finite element discretizations using extended-support $C^0$ continuous piecewise polynomial bases also possess a maximum stable timestep which scales more slowly in $p$ than traditional high order $C^0$ and DG finite element methods \cite{banks2016galerkin}.}

Finally, we investigate the use of optimal spline spaces in NURBS discretizations, which are defined using knot locations at the roots of a specific eigenfunction \cite{melkman1978spline}.  These optimal knot positions bear some resemblance to the nonlinear parametrization used in \cite{cottrell2006isogeometric} to eliminate spurious high frequency ``outlier'' modes present in isogeometric discretizations with uniform knots.  Numerical experiments show that, under certain circumstances, spline spaces defined using optimal knot vectors are more efficient than splines defined using uniform knot vectors in terms of degrees of freedom required to reach a certain error.  Moreover, for time-dependent problems, the use of optimal knot vectors increases the maximum stable timestep while retaining reasonable approximation properties.  

The structure of the paper is as follows: Section~\ref{sec:form2} discusses Galerkin formulations of hyperbolic PDEs on multi-patch geometries.  Section~\ref{sec:basis} introduces B-spline bases and optimal knot distributions, and introduces an inexpensive heuristic iteration to approximate optimal knots locations.  Section~\ref{sec:mol} discusses the efficient solution of the semi-discrete finite element system using explicit time integration, including weight-adjusted approximations to the inverse mass matrix and estimates for stable timestep restrictions based on constants in $hp$ inverse and trace inequalities.  Finally, Section~\ref{sec:num} presents numerical experiments illustrating salient properties of the proposed methods, and Section~\ref{sec:conc} concludes with some future directions.  

\section{Multi-patch discontinuous Galerkin formulations}
\label{sec:form2}

In this work, we are interested in the solution of the hyperbolic PDEs of acoustic wave propagation and linear advection over a domain $\Omega \in \mathbb{R}^d$ with boundary $\partial \Omega$.  We begin the model problem of time-dependent advection \reviewerOne{of some concentration $\phi$} (with a divergence-free velocity field $\bm{\beta}$) 
\begin{align*}
\pd{\phi}{t} + \Div \LRp{\bm{\beta}\phi} &= 0, \qquad \Div\bm{\beta} = 0\\
\phi(\bm{x},0) &= \phi_0(\bm{x}).
\end{align*}
Since the advection equation will primarily be used to study dispersive and dissipative properties of the proposed discretizations, periodic boundary conditions will be imposed on the boundary $\partial \Omega$.  

We consider also the acoustic wave equation as a first order pressure-velocity system
\begin{align}
\frac{1}{c^2}\pd{p}{t} + \Div \bm{u} &= f \label{eq:waveeq}\\
\pd{\bm{u}}{t} + \Grad p &= 0. \nonumber
\end{align} 
\reviewerOne{where $p$ and $\bm{u}$ represent the pressure and velocity variables, respectively, and $c$ is the wavespeed.}  
This problem is well-posed when paired with either Dirichlet or Neumann boundary conditions on $\partial \Omega$ and initial conditions 
\[
p(\bm{x},0) = p_0(\bm{x}), \qquad \bm{u}(\bm{x},0) = \bm{u}_0(\bm{x}).
\]
The second order form of the acoustic wave equation can be derived from the first order form.  Taking the derivative in time of the pressure equation and the divergence of the velocity equation and summing the two equations yields
\[
\reviewerOne{\frac{1}{c^2}}\pdd{p}{t} - \Delta p = f,  
\]
\reviewerOne{where $f$ is related to, but not necessarily the same as the $f$ in (\ref{eq:waveeq}).}
The second order formulation is well posed for appropriate boundary conditions and initial conditions 
\[
p(\bm{x},0) = p_0(\bm{x}), \qquad \pd{p}{t}(\bm{x},0) = v_0(\bm{x}).  
\]

%Let $\bm{n}$ be the outward unit normal on $\Omega$; the inflow and outflow boundaries $\Gamma_{\rm in}, \Gamma_{\rm out} \subset \partial \Omega$ are then defined as
%\[
%\Gamma_{\rm in} = \LRc{\bm{x}: \bm{\beta}\cdot \bm{n} < 0}, \qquad \Gamma_{\rm out} = \LRc{\bm{x}: \bm{\beta}\cdot \bm{n} \geq 0}
%\]
%Boundary conditions are imposed on the inflow domain $\left.u\right|_{\Gamma_{\rm in}} = u_0$. 

\begin{figure}
\centering
\includegraphics[width=.4\textwidth]{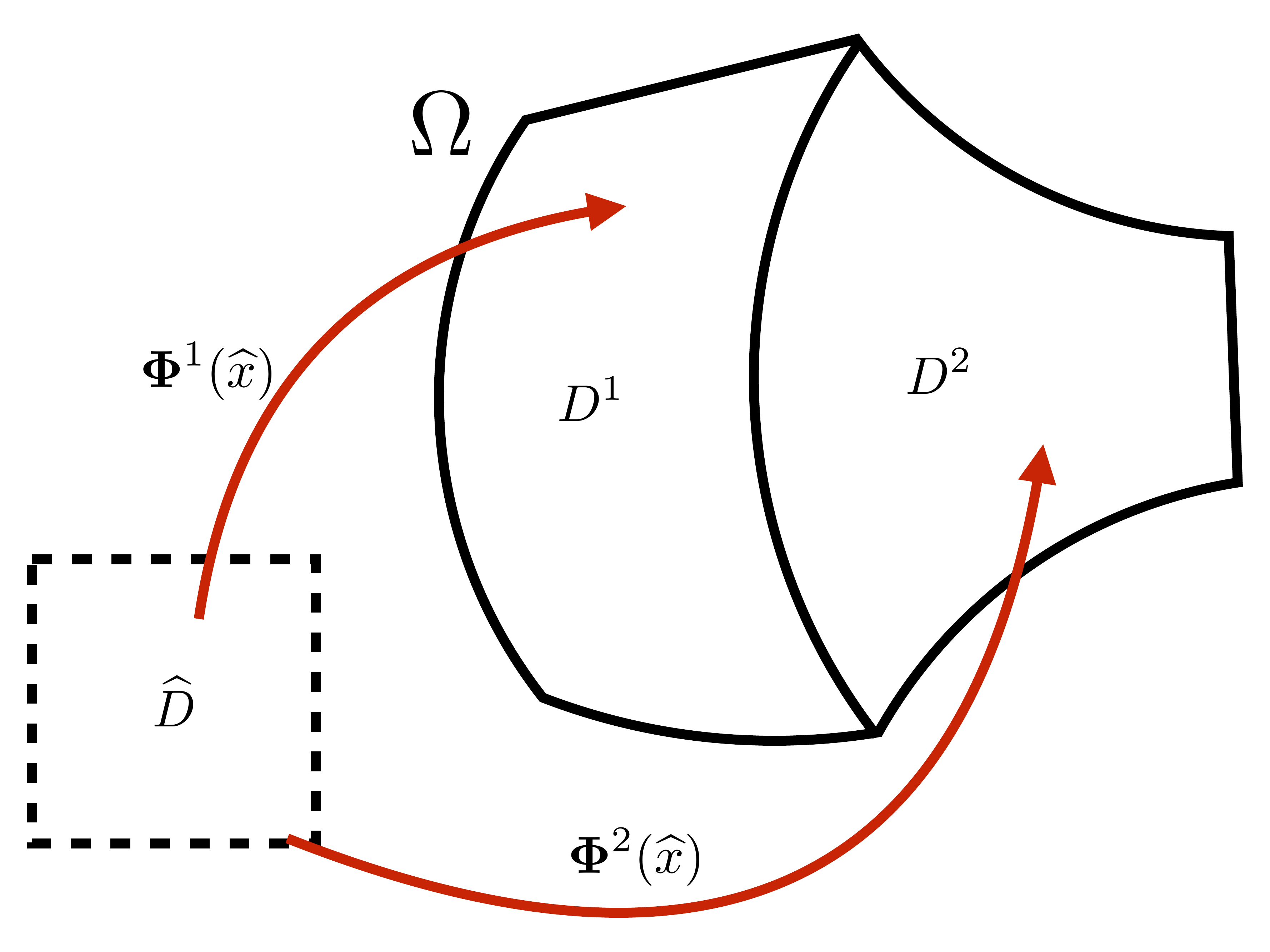}
\caption{An example of a multi-patch geometry.}
\label{fig:multipatch_diagram}
\end{figure}

In practical applications, sufficiently complex geometries are represented using multiple curved geometric patches.  We assume that $\Omega$ is decomposed into non-overlapping patches $D^k$, each of which is the image of the reference patch $\Dhat$ with coordinates $\widehat{\bm{x}}$ under a geometric mapping $\bm{x} = \bm{\Phi}_k(\bm{\widehat{x}})$, as shown in Figure~\ref{fig:multipatch_diagram}.   We define the mesh $\Omega_h = \cup D^k$ and the corresponding global approximation space $V_h(\Omega_h) = \bigoplus V_h\LRp{D^k}$, where $V_h\LRp{D^k}$ is the local approximation space over a single patch.  We further assume that the local approximation space over $D^k$ is the image of a reference approximation space $V_h\LRp{\widehat{D}}$ under $\bm{\Phi}_k$, and  that trace and inverse inequalities hold on $\Dhat$ such that 
\begin{equation}
\nor{u}_{L^2\LRp{\partial \Dhat}} \leq \sqrt{C_T} \nor{u}_{L^2\LRp{\Dhat}}, \qquad \nor{\widehat{\Grad} u}_{L^2\LRp{\Dhat}} \leq C_I \nor{u}_{L^2\LRp{\Dhat}}, \qquad \forall v\in V_h\LRp{\Dhat},
%\nor{u}^2_{L^2\LRp{\partial \Dhat}} \leq C_T \nor{u}^2_{L^2\LRp{\Dhat}}, \qquad \nor{\widehat{\Grad} u}^2_{L^2\LRp{\Dhat}} \leq C_I^2 \nor{u}^2_{L^2\LRp{\Dhat}}, \qquad \forall v\in V_h\LRp{\Dhat}
\label{eq:ineqs}
\end{equation}
where $C_T, C_I$ are constants which depend on the space $V_h\LRp{\Dhat}$.  We will describe specific approximation spaces in more detail in Section~\ref{sec:basis}.  %We note that we have used the norm squared for the trace constant, but only the norm for the inverse inequality constant.  This is to ensure that the scaling
{ We emphasize that the constant $C_I$ appears in the inverse inequality in (\ref{eq:ineqs}), but the constant $\sqrt{C_T}$ appears in the corresponding trace inequality.  The trace inequality $C_T$ in this work is thus the square of the constant typically presented $hp$ finite element trace and inverse inequalities \cite{ciarlet1978finite, warburton2003constants}.  This definition ensures that $C_T, C_I$ both scale similarly (as shown in Section~\ref{sec:numconsts}) and simplifies expressions for bounds involved in estimating stable timestep restrictions in Section~\ref{sec:cfl}.  }

We will discretize the equations of interest based on a discontinuous Galerkin (DG) approach \cite{hesthaven2007nodal,langer2015multipatch}, and will explore first order variational formulations of the acoustic wave equation and linear advection, as well as second order formulations of acoustic wave propagation.  We first introduce the jump and average of discontinuous functions across { patch} interfaces.  Let $D^{k,+}$ denote the neighboring {  patch} \reviewerTwo{sharing} a face $f$ of $D^k$, and let $u^+,u^-$ denote the values of $u$ on $D^{k,+}$ and $D^k$, respectively.  The jump \reviewerTwo{and average} of $u$ across $f$ is then defined as
\[
\jump{u} = u^+ - u^-, \qquad \avg{u} = \frac{u^+ + u^-}{2}.
\]
The jump and average of vector fields are defined component-wise using the jumps and averages of components.  

We note that all formulations presented in the following sections are constructed to be \textit{a-priori} energy stable so long as integrals are explicitly approximated using quadrature (as opposed to quadrature-free implementations).  Additionally, high order accuracy is observed in all numerical experiments where the quadrature is sufficiently accurate \cite{chan2016weight2}.  These formulations are especially important in context of curvilinear meshes, where exact evaluation integrals is inefficient or impossible due to the high order or rational nature of the determinant of the mapping Jacobian $J$.

\subsection{Formulation for linear advection}

We begin with the advection equation, for which we can construct a discontinuous Galerkin formulation based on \cite{kopriva2014energy}: \reviewerTwo{find $\phi \in V_h(\Omega_h)$ such that}
\begin{align}
\sum_k\LRp{\pd{\phi}{t},v}_{L^2\LRp{D^k}} + a(\phi,v) &= 0, \qquad \reviewerTwo{\forall v\in V_h(\Omega_h),}\label{eq:varformadvec}\\
a(\phi,v) &= a_{\rm skew}(\phi,v) + a_{\rm sym}(\phi,v) \nonumber\\
a_{\rm skew}(\phi,v) &= \frac{1}{2}\sum_k\int_{\Omega} \LRp{-\phi\bm{\beta}\cdot \Grad v + \bm{\beta}\cdot \Grad \phi v}\diff{{\bm{x}}}\nonumber\\
a_{\rm sym}(\phi,v) &= \frac{1}{2}\sum_{D^k}\int_{\partial D^k} \LRp{\avg{\bm{\beta}}\cdot\bm{n} - \tau\LRb{\avg{\bm{\beta}}\cdot\bm{n}}} \jump{\phi} v \diff{{\bm{x}}}.\nonumber  
\end{align}
where $0 \leq \tau$ is a parameter which switches between a dissipative flux ($\tau > 0$) and a central flux ($\tau = 0$).  Assuming appropriate boundary conditions (homogeneous or periodic), taking $v = \phi$ shows the formulation to be energy stable, such that
\[
\reviewerOne{\pd{}{t}\sum_{k}\LRp{\phi,\phi}_{D^k} \leq 0.}
\]

\subsection{First order formulation of the acoustic wave equation}

We now consider the acoustic wave equation, and follow approaches outlined in \cite{chan2016weight1, chan2016weight2} to construct an energy-stable discretization using a skew-symmetric formulation and penalty flux.  This formulation is given for $U = \LRc{p,\bm{u}}, V = \LRc{q,\bm{v}}$ as follows: \reviewerTwo{find $U \in V_h(\Omega_h) \times V_h(\Omega_h)^d$ such that}
\begin{align}
\sum_k\LRp{\frac{1}{c^2}\pd{p}{t},q}_{L^2\LRp{D^k}} + \LRp{\pd{\bm{u}}{t},\bm{v}}_{L^2\LRp{D^k}} + a(U,V) &= \sum_k\LRp{f,q}_{L^2\LRp{D^k}}, \qquad \reviewerTwo{V \in V_h(\Omega_h) \times V_h(\Omega_h)^d,} \label{eq:varformwave1}
\end{align}
where we define $a(U,V)$ as
\begin{align*}
a(U,V) &= \sum_k\LRp{-\LRp{\bm{u},\Grad q}_{L^2\LRp{D^k}} + \frac{1}{2}\LRa{2\avg{\bm{u}}\cdot\bm{n} - \tau_p \jump{p},q}_{\partial D^k}} \\
 &+ \sum_k\LRp{\LRp{\Grad p,\bm{v}}_{L^2\LRp{D^k}} + \frac{1}{2}\LRa{\jump{p} - \tau_u \jump{\bm{u}}\cdot\bm{n},\bm{v}\cdot\bm{n}}_{\partial D^k}},
\end{align*}
where $\tau_p, \tau_u$ are positive penalization parameters, and both the pressure $p$ and each component of the velocity $\bm{u}$ are approximated locally from $V_h\LRp{D^k}$.  The variational formulation over each patch can be mapped \reviewerTwo{(or pulled back)} to $\Dhat$ as
\begin{align*}
\int_{\Dhat} \pd{p}{t} q J\diff{\bm{\widehat{x}}} + \int_{\Dhat} \widehat{\Grad}\cdot \LRp{\bm{G}^T\bm{u}}qJ\diff{\bm{\widehat{x}}} + \frac{1}{2}\int_{\partial \widehat{D}}\LRp{2\avg{\bm{u}}\cdot\bm{n} - \tau_p \jump{p}}J^s &= \int_{\Dhat} f q J\diff{\bm{\widehat{x}}}\\
\int_{\Dhat} \pd{\bm{u}}{t} \bm{v} J\diff{\bm{\widehat{x}}} + \int_{\Dhat} \bm{G}\widehat{\Grad} {p}\bm{v}J\diff{\bm{\widehat{x}}} + \frac{1}{2}\int_{\partial \widehat{D}}\LRp{\jump{p} - \tau_u \jump{\bm{u}}\cdot\bm{n}}\bm{v}\cdot\bm{n} J^s &=0, 
\end{align*}
where $\bm{G}_{ij} = \pd{\widehat{\bm{x}}_i}{\bm{x}_j}$ \reviewerTwo{is the inverse of the Jacobian matrix of $\bm{\Phi}_k(\hat{x})$ (also referred to as the metric terms \cite{kopriva2006metric} or geometric factors \cite{hesthaven2007nodal})}, and $J, J^s$ denote the determinants of volume and surface geometric Jacobians, respectively.  

The global formulation results from the summation of local formulations over each {  patch}.  For constant wavespeed $c^2$  and $\tau_p, \tau_u = 1$, the formulation reduces to the standard upwind DG method.  A standard energy argument gives that the method is dissipative for $\tau_p, \tau_u > 0$, in the sense that 
\[
\pd{}{t}\sum_k \int_{D^k}\LRp{\frac{p^2}{c^2} + \LRb{\bm{u}}^2}\diff{{\bm{x}}} \leq -\sum_k \int_{\partial D^k} \LRp{\tau_p\jump{p}^2 + \tau_u\LRb{\jump{\bm{u}}\cdot \bm{n}}^2}\diff{{\bm{x}}}.
\]
The dissipation resulting from the stabilization terms provides a natural way to damp spurious non-conforming components of the solution in time-domain simulations \cite{chan2016short}.  Similar energy stable skew-symmetric formulations can be constructed for electromagnetics and elastodynamics \cite{warburton2013low, chan2017weight}.  

Dirichlet boundary conditions on pressure are enforced in an energy-stable fashion by specifying the jump of $p$ on boundary faces 
\[
\jump{p} = 2 (p_D - p),
\]
where $p_D$ is the prescribed value of $p$ on the boundary.  Similarly, Neumann boundary conditions are enforced by specifying the jump of $\bm{u}$ on boundary faces
\[
\jump{\bm{u}}\cdot \bm{n} = 2 (u_N - \bm{u}\cdot\bm{n}),
\]
where $u_N$ is the prescribed value of the normal flux on the boundary.

\subsection{Second order formulation of the acoustic wave equation}

We use the acoustic wave equation as a prototypical example of a second order formulation.  We discretize the Laplacian using a symmetric interior penalty DG method \cite{arnold1982interior, langer2015multipatch}.  The resulting formulation is given as follows: \reviewerTwo{find $p\in V_h(\Omega_h)$ such that}
\begin{align}
\sum_k\LRp{\pdd{p}{t},v}_{L^2\LRp{D^k}} + a(p,v) &= \sum_k\LRp{f,v}_{L^2\LRp{D^k}}, \qquad \reviewerTwo{\forall v \in V_h(\Omega_h)}\\
a(p,v) &= \sum_k\LRp{\Grad p,\Grad v}_{L^2\LRp{D^k}}  -\LRa{\avg{\Grad p}\cdot\bm{n}, v}_{L^2\LRp{\partial D^k}} \label{eq:varformwave2}\\
&+ \sum_k \LRa{\jump{p}, \Grad v\cdot\bm{n}}_{L^2\LRp{\partial D^k}} -\tau\LRa{\jump{p},v}_{L^2\LRp{\partial D^k}}\nonumber.
\end{align}
where $\tau$ is a penalization parameter.  In contrast with first order systems, the second order system achieves optimal rates of convergence in the $L^2$ norm and is non-dissipative and energy conserving.  % (though dissipative second order formulations have been introduced in \cite{appelonew}).  

We note that, in order to ensure coercivity and stability, the penalization constant $\tau$ must be sufficiently large with respect to discretization parameters.  We follow the approach in \cite{Shahbazi:2005:SNE:1083044.1083046}, where a straightforward generalization of the proof to curved patches yields the following sufficient global condition for coercivity
\[
%\tau \geq \max_{D^k} {C_T^2}{\nor{J^s}_{L^\infty\LRp{\partial D^k}}\nor{\frac{1}{J}}_{L^\infty\LRp{D^k}}}.
\tau \geq \max_{D^k} {C_T}{\nor{J^s}_{L^\infty\LRp{\partial D^k}}\nor{\frac{1}{J}}_{L^\infty\LRp{D^k}}}.
\]
Here, $C_T$ is the trace constant over the reference patch, while $J, J^s$ are Jacobian factors for the volume and surface (face) of a physical patch, respectively.\footnote{One may also define the penalization constant locally over each patch face; this may be useful for local time-stepping, as the maximum stable timestep is dependent on the magnitude of $\tau$.}

\section{B-spline bases}
\label{sec:basis}

We are interested in approximating solutions over each patch $D^k$ from spline spaces, which are piecewise polynomials of degree $p$ with up to $p-1$ continuous global derivatives.  The following sections briefly review the construction of B-spline basis functions, and introduce choices of knot vectors which modify the approximation properties of the corresponding spline space.

Univariate spline spaces of degree $p$ are constructed given an ordered knot vector $\bm{\Xi}$ 
\[
\bm{\Xi} = \LRc{-1 = \xi_1,\ldots,\xi_{{2p}+K+1} = 1}, \qquad \xi_i \leq \xi_{i+1}.  
\]
The unique (i.e. non-repeating) knots in $\bm{\Xi}$ define $K$ non-overlapping sub-intervals (elements) over which the spline space is defined.  B-spline basis functions can then be constructed recursively as follows
\[
B^0_i(\widehat{x}) = \begin{cases}
1, & \xi_i \leq \widehat{x} \leq \xi_{i+1}\\
0, & {\rm otherwise},
\end{cases}
\qquad 
B^k_i(\widehat{x}) = \frac{\widehat{x} - \xi_i}{\xi_{i+p}-\xi_i}B^{k-1}_i(\widehat{x}) + \frac{\xi_{i+p+1}-\widehat{x}}{\xi_{i+p+1}-\xi_{i+1}}B^{k-1}_{i+1}(\widehat{x}).
\]
In the case when $\xi_{i+p}-\xi_i = 0$, the ratio $\frac{\widehat{x} - \xi_i}{\xi_{i+p}-\xi_i}$ is also set to zero.  Similarly, if $\xi_{i+p+1}-\xi_{i+1} = 0$, $ \frac{\xi_{i+p+1}-\widehat{x}}{\xi_{i+p+1}-\xi_{i+1}}$ is also set to zero.  The resulting basis functions $B^p_i(\widehat{x})$ have local support and span a piecewise polynomial space of degree $p$.  In this work, we focus on \textit{open} knot vectors $\bm{\Xi}$, where the first and last knots are repeated $p$ times
\[
\xi_1 = \ldots = \xi_{p+1}, \qquad \xi_{p+1 + K} = \ldots = \xi_{2p+K+1}.
\]
The resulting bases contain $p+K$ basis functions, are interpolatory at the endpoints, and contain the space of degree $p$ polynomials $P^p([-1,1])$ as a proper subset.  For $K=1$, this basis recovers the degree $p$ Bernstein polynomial basis.  Finally, we define the one-dimensional degree $p$ approximation space over the reference interval $V_h([-1,1])$ as the span of B-splines 
\[
V_h([-1,1]) = {\rm span}\LRc{B^p_i(\widehat{x})}_{i=1}^{p+K}.
\]
For the remainder of this work, we focus on spline spaces of maximal continuity, such that the resulting degree $p$ splines possess $p-1$ continuous derivatives.  { Moreover, we focus on uniform knot vectors, that is, knot vectors with sub-integrals of equal length.} %Spline spaces of reduced continuity may also be used, but are not the current focus of this work.   

\reviewerOne{Any function in $V_h([-1,1])$ is represented as a linear combination of B-spline basis functions
\[
u(\hat{x}) = \sum_{j=1}^{p+K} c_j B^p_j(\hat{x}), \qquad \forall u \in V_h([-1,1]),
\]
where $c_j$ are the coefficients or ``control points'' of the curve defined by $u(x)$.  We note that $c_1, c_{p+K}$ interpolate the values of $u(-1)$ and $u(1)$, respectively.  However, since interior basis functions are not interpolatory, the solution does not necessarily coincide with the value of the interior coefficients/control points. }

\subsection{Higher dimensional splines and approximation spaces}

Given a one-dimensional spline basis, we can extend spline bases to $d$ dimensions through a tensor product construction. For example, B-spline basis functions $B^{p}_{ij}(\hat{x},\hat{y})$ in two dimensions and $B^{p}_{ijk}(\hat{x},\hat{y},\hat{z})$ in three dimensions can be constructed as
\[
B^{p}_{ij}(\widehat{x},\widehat{y}) = B^p_i(\widehat{x}) B^p_j(\widehat{y}), \qquad B^{p}_{ijk}(\widehat{x},\widehat{y},\widehat{z}) = B^p_i(\widehat{x}) B^p_j(\widehat{y})B^p_k(\widehat{z}), \qquad 1 \leq i,j,k\leq p+K.
\]
We can then define the local spline approximation space of order $p$ over the reference patch $\widehat{D}$ as 
\[
V_h\LRp{\widehat{D}} = {\rm span}\LRc{B^p_{i}\LRp{\bm{\widehat{x}}}}_{i=1}^{(p+K)^d},
\]
where $B^p_i\LRp{\bm{\widehat{x}}}$ denotes the $i^{\rm th}$ tensor product spline in $d$ dimensions \reviewerTwo{(here, we have switched from a multi-index $ijk$ to a single index $i$).}

\begin{figure}[t]
\centering
\includegraphics[width=.85\textwidth]{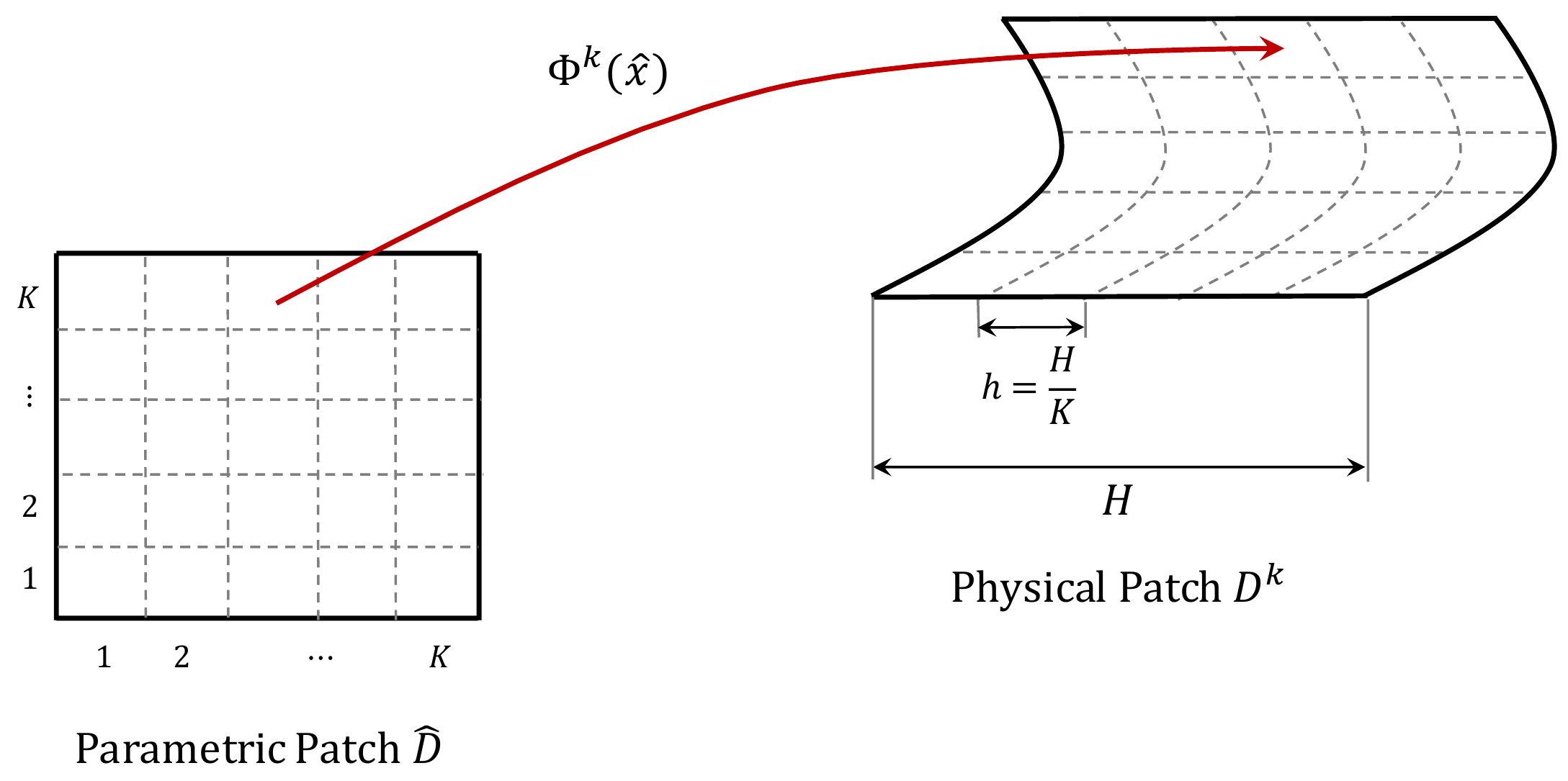}
\caption{{ Illustrations of the patch size $H$, the mesh size $h$, and the number of elements per side $K$ for a particular patch $D^k$.  The dashed lines represented the knot lines both on the parametric patch $\hat{D}$ and the physical patch $D^k$.}}
\label{fig:definitions}
\end{figure}

For isotropic spline approximation spaces, we define the patch size $H$ as the diameter {(which we take to mean the longest side)} of $D^k$.  We will refer to subdivision of patches into patches of size $H/2$ as patch refinement, while mesh refinement will refer to the generation of new knot vectors with $2K$ elements each. Finally, we refer to the physical mesh size as the patch size divided by the number of elements $h = H/K$.  We note that this definition assumes quasi-uniform patches and knot vectors.  { We have graphically depicted all of the aforementioned geometric quantities in Figure \ref{fig:definitions}. Note that when only one element is included per patch, as is the case with the standard discontinuous Galerkin method using discontinuous piecewise polynomial approximations, then $K = 1$ and hence $h = H$.  Thus, the terms element and patch are interchangeable for this setting.}  Because the formulations presented in Section~\ref{sec:form2} are standard Galerkin formulations using spline basis functions, they are consistent and stable in the finite element sense, and the $h$-convergence of these methods on curved geometries can be proven using techniques adapted from \cite{warburton2013low} and interpolation estimates and inverse inequalities established in \cite{bazilevs2006isogeometric, evans2013explicit}.

Finally, we note that in engineering applications, the geometric mapping $\bm{\Phi}_k$ is often defined using a rational NURBS mapping.  An approach popularized in isogeometric analysis is the use of a NURBS space (consisting of rational B-splines) as an approximation space as well.  For the numerical experiments presented in this work, we do not utilize NURBS approximation spaces, and instead map B-splines defined on $\Dhat$ to the physical domain using the mapping $\bm{\Phi}_k$.  \reviewerTwo{However, the extension to NURBS is straightforward, and Section~\ref{sec:wadgnurbs} briefly describes how to extend the techniques of Section~\ref{sec:wadg} to NURBS approximations.}

\subsection{Approximation properties and optimal knot vectors}

In this work, we restrict ourselves to open knot vectors, and seek interior knot locations which improve approximation properties of the underlying spline space.  It is clear that the span of B-splines contains polynomials of degree $p$ over $[-1,1]$; however, the relationship between approximation properties and interior knot locations is more complex.  This relationship has been analyzed using sup-infs and $n$-widths \cite{kolmogoroffOber1936, melkman1978spline, schumaker2007spline, evans2009n}.  Let $W$ be a normed linear space, let $X \subset W$, and $X_n$ be some $n$-dimensional subspace of $X$; then, the sup-inf $E(X;X_n)$ of $X$ relative to $X_n$ is defined as follows: 
\[
E(X;X_n) =  \sup_{x\in X} \inf_{y\in X_n} \nor{x - y},
\]
and can be interpreted as measuring the worst possible best approximation error in approximating $X$ from $X_n$.  The $n$-width of $X$ is then defined to be the infimum of $E(X;X_n)$ over all subspaces $X_n$ with ${\rm dim}\LRp{X_n} = n$, and a specific subspace $X_n$ is referred to as an \textit{optimal} subspace if its sup-inf attains the $n$-width.

Melkman and Michelli showed in \cite{melkman1978spline} that, for $W = L^2([-1,1])$ and 
\[
X = \LRc{f \in W: \pdn{f}{x}{r} \text{ is absolutely continuous on } [-1,1], \quad \nor{f}_{L^2([-1,1])} \leq 1},
\]
an optimal subspace $X_n$ is the spline space of degree $r$ with knot locations at the roots of $y_{n+1,r}(x)$, where $y_{n+1,r}(x)$ is the $(n+1)$st eigenfunction of the boundary value problem
\begin{equation}
(-1)^r\pdn{y}{x}{2r} = \lambda y(x), \qquad \pdn{y}{x}{k}(-1) = \pdn{y}{x}{k}(1) = 0, \quad k = 0, 1, \ldots, r-1.
\label{eq:nwidtheig}
\end{equation}
We note that, for $r=1$, this reduces down the standard Laplacian eigenvalue problem with Dirichlet boundary conditions, whose eigenfunctions $y_{j,1}(x)$ are 
\[
y_{j,1}(x) = \sin\LRp{j \pi \LRp{\frac{1+x}{2}}}, \qquad j = 1,2,\ldots, \qquad x \in (-1,1).
\]
Thus, the optimal $n$-dimensional space for $r = 1$ consists of piecewise linear splines with $n$ equispaced knots along $[-1,1]$, which are simply piecewise linear $C^0$ polynomials on a mesh with $K = n$ elements.  Similarly, for $r = 2$, (\ref{eq:nwidtheig}) reduces to the biharmonic eigenvalue problem with clamped plate conditions \cite{brenner2007mathematical}.  

\begin{figure}
\centering
\subfloat[$r = 2$]{\includegraphics[width=.32\textwidth]{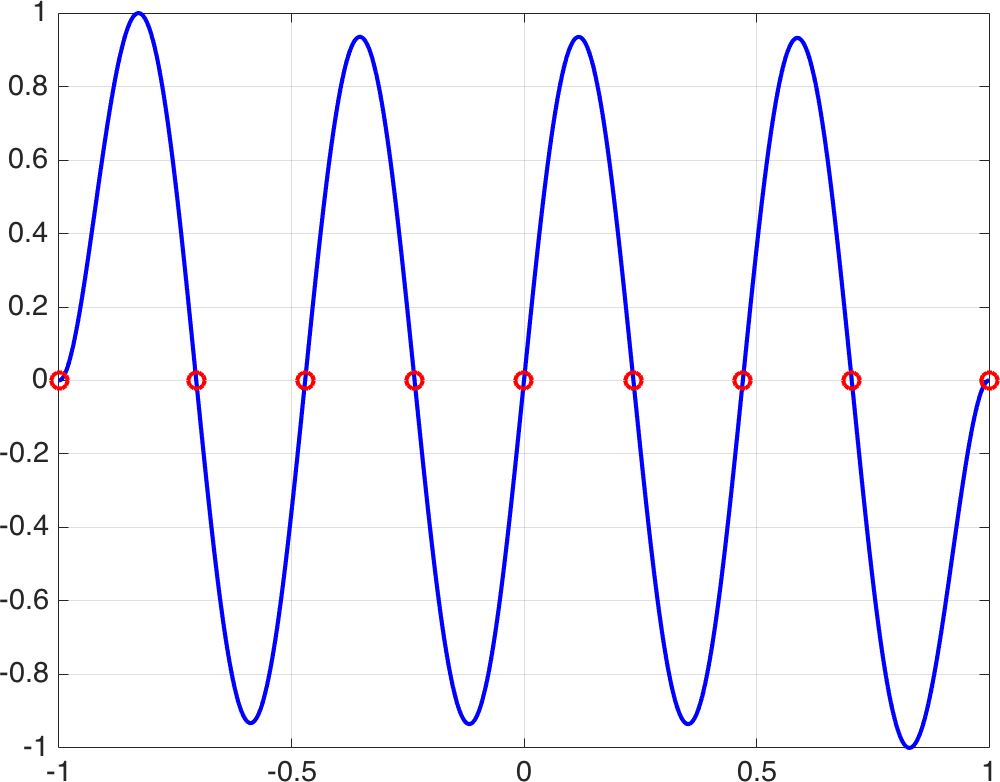}}
\hspace{.5em}
\subfloat[$r = 3$]{\includegraphics[width=.32\textwidth]{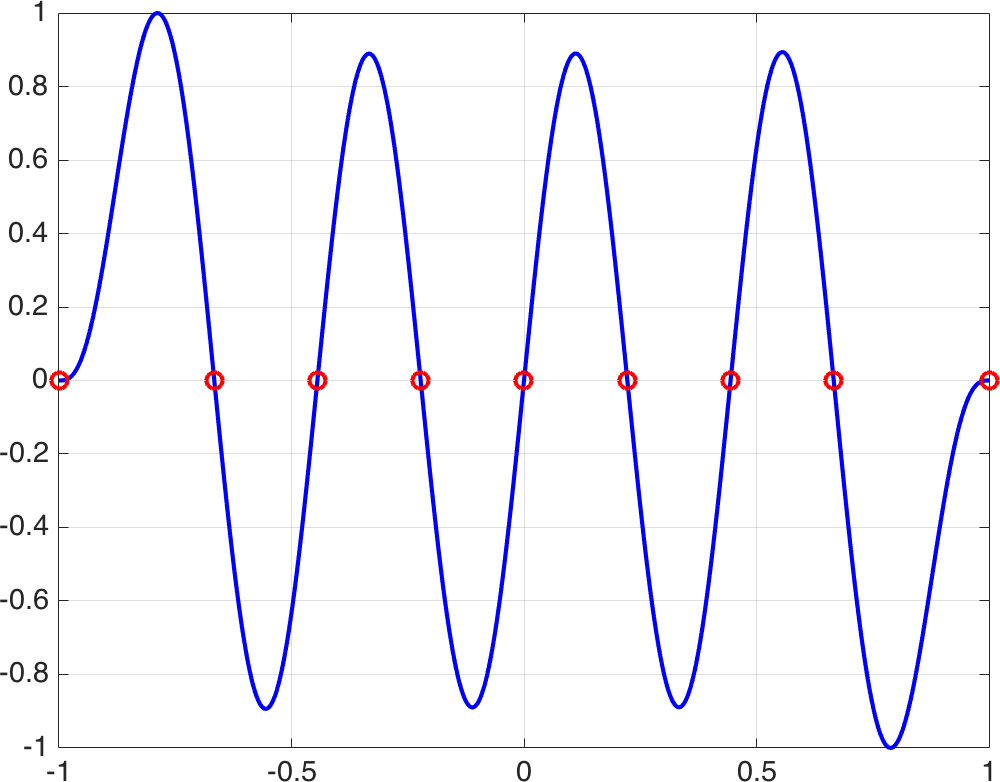}}
\hspace{.5em}
\subfloat[$r = 4$]{\includegraphics[width=.32\textwidth]{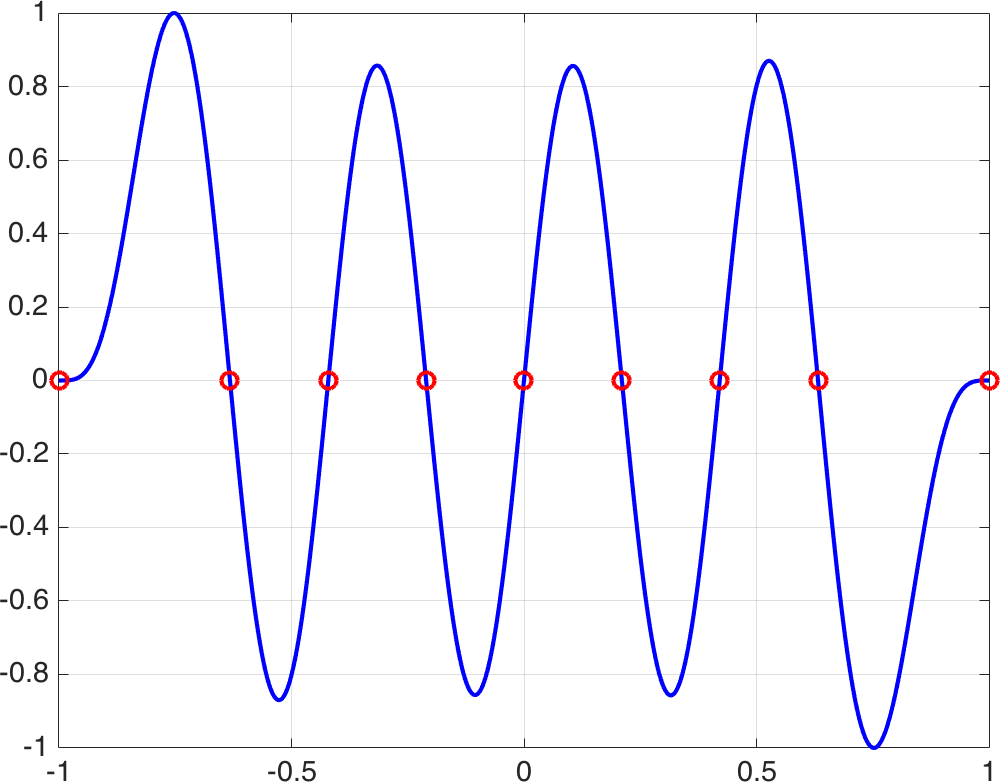}}
\caption{Eigenfunctions $y_{r,n}(x)$ and their roots for $n = 8$ and various $r$.  All eigenfunctions are approximated using a symmetric Galerkin formulation with degree $p = 8$ spline spaces of maximal continuity and $K = 64$ elements. }
\label{fig:eigoptknots}
\end{figure}

For $r > 2$, optimal knot locations are not explicitly known.  However, following \cite{evans2009n}, we can approximate knot locations using a Galerkin approximation of (\ref{eq:nwidtheig}) with $C^{r-1}$ approximation spaces.  Figure~\ref{fig:eigoptknots} shows approximations of eigenfunctions $y_{n+1,r}(x)$ for various $n,r$ using a fine spline space, with knot locations overlaid.  These knot locations are determined by approximating the roots of eigenfunctions numerically up to a tolerance of $10^{-12}$.  We observe that, as $r$ increases, 
%these knot locations cluster towards the interior relative to equispaced knot locations 
\reviewerOne{the gap between each of the boundary knots and the closest interior knot widens}
due to the presence of the homogeneous boundary conditions on derivatives up to degree $(r-1)$.  

\subsection{Knot smoothing}

It is possible to pre-compute and store optimal knot locations for small to moderate $r$ and $n$ for use in practical applications.  However, for large $n$, optimal knot locations become difficult to approximate due to the resolution required to accurately resolve $y_{n+1,r}(x)$.  For large $r$, symmetric Galerkin discretizations of (\ref{eq:nwidtheig}) also become ill-conditioned due to the presence of high order derivatives.  The ill-conditioning of (\ref{eq:nwidtheig}) for $r \gg 1$ may potentially be addressed using a mixed formulation of the eigenvalue problem; however, the construction of stable and convergent mixed methods for higher order elliptic operators is non-trivial \cite{gudi2011interior}. Thus, we seek a heuristic alternative to approximate optimal knot distributions without the computational cost and conditioning issues associated with the discretization and solution of (\ref{eq:nwidtheig}).  We propose an alternative approach to determining a ``good'' knot vector based on an iterative redistribution of knot locations based on the associated Greville abscissae.  

\begin{figure}
\centering
\includegraphics[width=.5\textwidth]{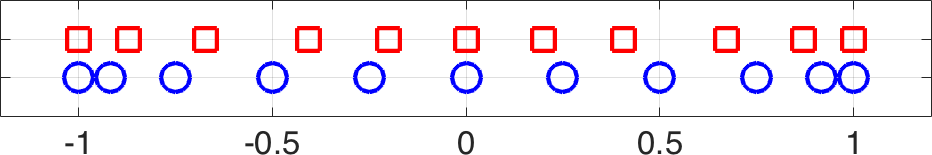}
\caption{Comparison of Greville abscissae for an $p=3, K = 8$ spline space when using uniform knots (shown in \textcolor{blue}{blue}) and knot positions $\tilde{\xi}_i$ defined by equation (\ref{eq:newKnots}) (shown in \textcolor{red}{red}).}
\label{fig:grevilleComparison}
\end{figure}

The Greville abscissae are a collection of points associated with a B-spline knot vector \cite{johnson2005higher}, and are defined as knot averages 
\[
\tau_j = \frac{1}{p}\sum_{i=1}^p \xi_{i+j-1}, \qquad j = 1,\ldots,p.
\]  
It can also be shown that these points satisfy
\[
\sum_{j=1}^{p+K} \tau_j B^p_j(\widehat{x}) = \widehat{x}.
\]
In other words, the Greville abscissae are the coefficients which reproduce the coordinate $\widehat{x}$ under a one-dimensional B-spline basis.  Because the Greville points $\tau_j$ reproduce $\widehat{x}$, the knot locations $\xi_i$ of a knot vector satisfy 
\[
\xi_i = \sum_{j=1}^{p+K} \tau_j B^p_j(\xi_i).  
\]
For an open knot vector, the distribution of the Greville abscissae clusters near the endpoints of the interval (similarly to stable polynomial interpolation nodes, such as Chebyshev or Gauss-Lobatto points \cite{canuto2012spectral}), as shown in Figure~\ref{fig:grevilleComparison}.  

In \cite{hughes2008duality}, the authors introduce a nonlinear parametrization of the interval $[-1,1]$, which is defined implicitly as the mapping which takes Greville abscissae to uniformly spaced points.  Numerical experiments showed that, by using the image of the spline space $V_h$ under this nonlinear mapping in a Galerkin discretization of the Laplacian, the so-called ``outlier frequencies'' of spurious high frequency eigenmodes present in isogeometric discretizations are removed completely.  However, under this mapping, the spline approximation space no longer contains polynomials of degree $p$, and approximation properties are not immediately clear.  Motivated by the results in \cite{hughes2008duality}, we seek instead a heuristic modification of a uniform open knot vector such that the resulting Greville abscissae are closer to equispaced points.  Moreover, we will show that this heuristic modification results in a reasonable approximation of optimal knot locations.   %In contrast, spline spaces defined by smoothed knot vectors are still polynomially complete, and (as demonstrated by the following numerical experiments) appear to produce good approximations of optimal knot locations.  

We wish to define a new knot vector $\tilde{\bm{\Xi}}$ such that the resulting Greville abscissae more closely resemble equispaced points.  We can introduce new knot positions $\tilde{\xi}_i$ by defining a linear parametrization where the Greville abscissae $\tau_j$ are \reviewerOne{arbitrarily} replaced with equispaced points $\widehat{x}_i$ along $[-1,1]$
\begin{equation}
\tilde{\xi}_i = \sum_{j=1}^{p+K} \reviewerOne{\widehat{x}_j} B^p_j(\xi_i).  
\label{eq:newKnots}
\end{equation}
\reviewerOne{One can interpret (\ref{eq:newKnots}) as attempting to reverse-engineer a knot vector for which the Greville abscissae are equispaced.  The goal is to produce new knot locations $\tilde{\xi}_i$ for which the Greville abscissae are less clustered towards the endpoints. }  The positions of the first and last knots produced by (\ref{eq:newKnots}) remains the same because the B-spline basis is interpolatory at the endpoints.  However, the locations of interior knots and Greville abscissae are squeezed closer towards the center of the interval, as shown in Figure~\ref{fig:grevilleComparison}.\footnote{We note that one might seek to further redistribute knots to produce Greville abscissae which are equispaced; however, because Greville abscissae are averages of knot positions and because open knot vectors contain repeating knots at endpoints, producing equispaced Greville points by simply redistributing knot positions is impossible without violating the assumption that knot positions $\xi_i$ must be monotonically increasing.}  

\begin{figure}
\centering
\subfloat[Knot positions]{\includegraphics[width=.4\textwidth]{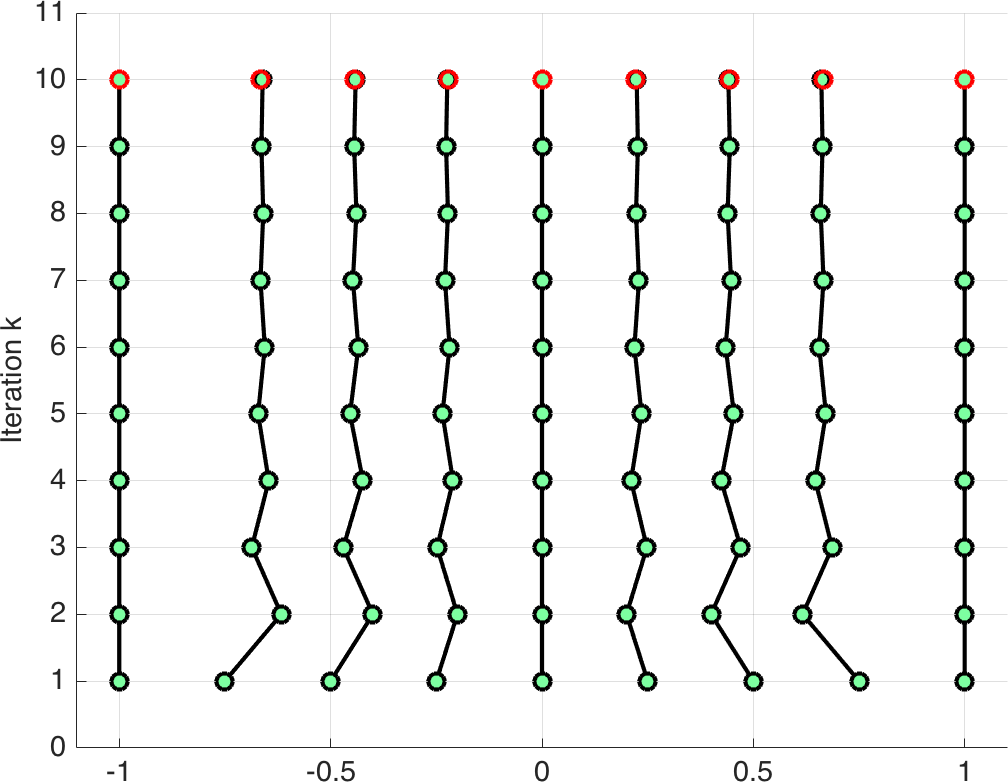}}
\hspace{1em}
\subfloat[Greville abscissae]{\includegraphics[width=.4\textwidth]{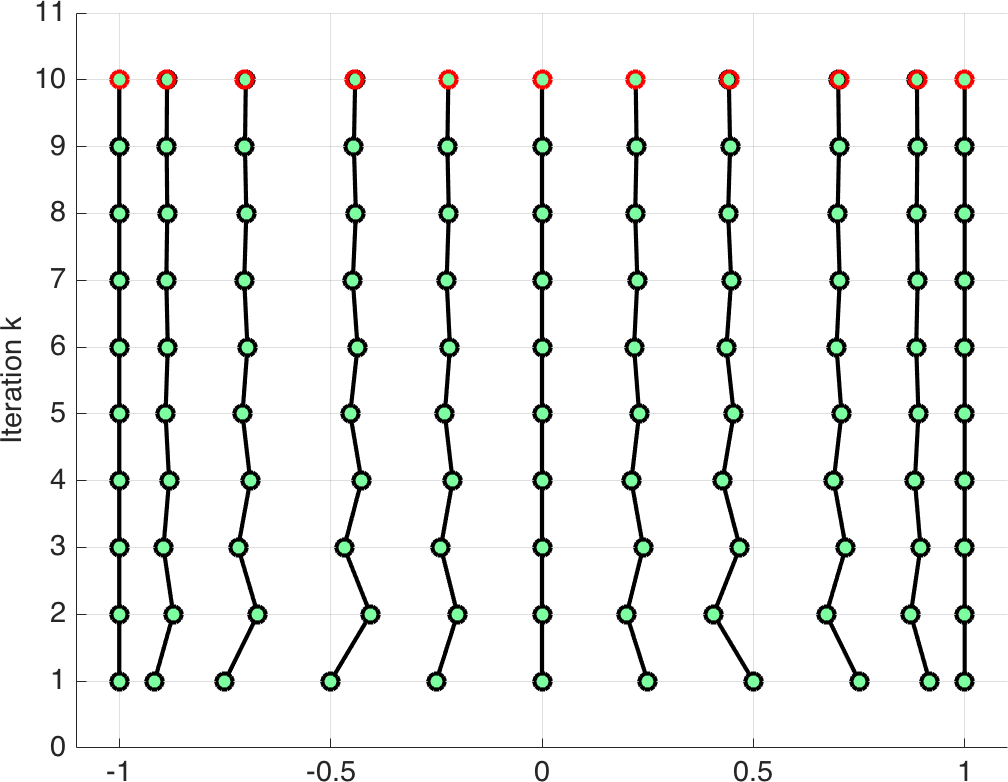}}
\caption{Knot locations for iterations $k = 0, \ldots, 10$ for a spline space with $p = 3, K = 8$.  Computed optimal knot positions and Greville abscissae are overlaid in \textcolor{red}{red}.}
\label{fig:smoothedknots}
\end{figure}

Because the spline basis depends on the knot vector, this process can be repeated 
\begin{equation}
\tilde{\xi}^{k+1}_i = \sum_{j=1}^{p+K} \widehat{x}_j B^p_j(\xi_i;  \tilde{\xi}^k), \qquad \tilde{\xi}^0_i = \xi_i,
\label{eq:knotiter}
\end{equation}
where we have used the notation $B^p_j(\xi_i;  \tilde{\xi}^k)$ to emphasize that the spline basis is evaluated at the original equispaced (i.e. uniform open) knots $\xi_i$, but defined using the knot vector $\tilde{\xi}^k$.\footnote{
We have also tested a knot smoothing iteration based on
\begin{equation}
\tilde{\xi}^{k+1}_i = \sum_{j=1}^{p+K} \widehat{x}_j B^p_j( \tilde{\xi}^k_i;  \tilde{\xi}^k), \qquad \tilde{\xi}^0_i = \xi_i,
\end{equation}
where splines are evaluated at the iterated knots $\tilde{\xi}^k_i$ instead of equispaced knots $\xi_i$.  However, we observe that this iteration produces collapsed knot locations, resulting in spline spaces which are not maximally continuous.  
}
Numerical experiments show that $\tilde{\xi}^k$ converges to an equilibrium distribution as $k\rightarrow \infty$, which can also be determined through the solution of a nonlinear system of equations.  Figure~\ref{fig:smoothedknots} shows the distribution of the knot vectors and Greville abscissae for each iteration $k$, along with optimal knot positions overlaid in red at the final iteration.  Rather than squeezing knots further towards the center, taking $k \rightarrow \infty$ results in knot locations which are close to the optimal knot positions.  

\begin{figure}
\centering
\subfloat[Normalized knot error $\delta_{\rm knot}$ ]{
\begin{tikzpicture}
\begin{axis}[
width=.425\textwidth,
    xlabel={Number of elements $K$},   
    xmin=2, xmax=66,
    ymax=0,ymax=0.25,    
    legend pos=north east, legend cell align=left, legend style={font=\tiny},	
       xmajorgrids=true,  ymajorgrids=true, grid style=dashed,
] 

\addplot[color=blue,mark=*,semithick, mark options={fill=markercolor}]
%coordinates{(4,0.0121)(8,0.0075)(16,0.0046)(32,0.0023)(64,0.0012)};
coordinates{(4,0.210397)(8,0.165792)(16,0.156554)(32,0.15437)(64,0.153806)};

\addplot[color=red,mark=square*,semithick, mark options={fill=markercolor}]
%coordinates{(4,0.009)(8,0.0041)(16,0.0028)(32,0.002)(64,0.0012)};
coordinates{(4,0.0887031)(8,0.0485602)(16,0.0506855)(32,0.0693288)(64,0.0798586)};

\addplot[color=black,mark=triangle*,semithick, mark options={fill=markercolor}]
%coordinates{(4,0.0016)(8,0.0101)(16,0.0074)(32,0.0048)(64,0.0028)};
coordinates{(4,0.0290577)(8,0.0819356)(16,0.099327)(32,0.115575)(64,0.12577)};

\legend{$p = 2$, $p=3$, $p=4$}
\end{axis}
\end{tikzpicture}
}
\hspace{2em}
\subfloat[Normalized Greville error $\delta_{\rm Greville}$]{
\begin{tikzpicture}
\begin{axis}[
width=.425\textwidth,
    xlabel={Number of elements $K$},   
    xmin=2, xmax=66,
    ymax=0,ymax=0.25,    
    legend pos=north east, legend cell align=left, legend style={font=\tiny},	
       xmajorgrids=true,  ymajorgrids=true, grid style=dashed,
] 

\addplot[color=blue,mark=*,semithick, mark options={fill=markercolor}]
%coordinates{(4,0.006)(8,0.0038)(16,0.0024)(32,0.0011)(64,0.00058806)};
coordinates{(4,0.21038)(8,0.100464)(16,0.0851055)(32,0.0805913)(64,0.0788944)};

\addplot[color=red,mark=square*,semithick, mark options={fill=markercolor}]
%coordinates{(4,0.003)(8,0.002)(16,0.0018)(32,0.0012)(64,0.00067857)};
coordinates{(4,0.0887023)(8,0.0365552)(16,0.0390282)(32,0.043861)(64,0.0468188)};

\addplot[color=black,mark=triangle*,semithick, mark options={fill=markercolor}]
%coordinates{(4,0.000401)(8,0.0025)(16,0.0016)(32,0.00089352)(64,0.00063901)};
coordinates{(4,0.0290725)(8,0.040789)(16,0.027371)(32,0.0237124)(64,0.0301967)};

\legend{$p = 2$, $p=3$, $p=4$}
\end{axis}
\end{tikzpicture}
}
\caption{Normalized difference between the locations of smoothed and optimal knots for $p = 2, 3, 4$ and $K = 4, 8, 16, 32, 64$. }
\label{fig:smoothedknoterror}
\end{figure}
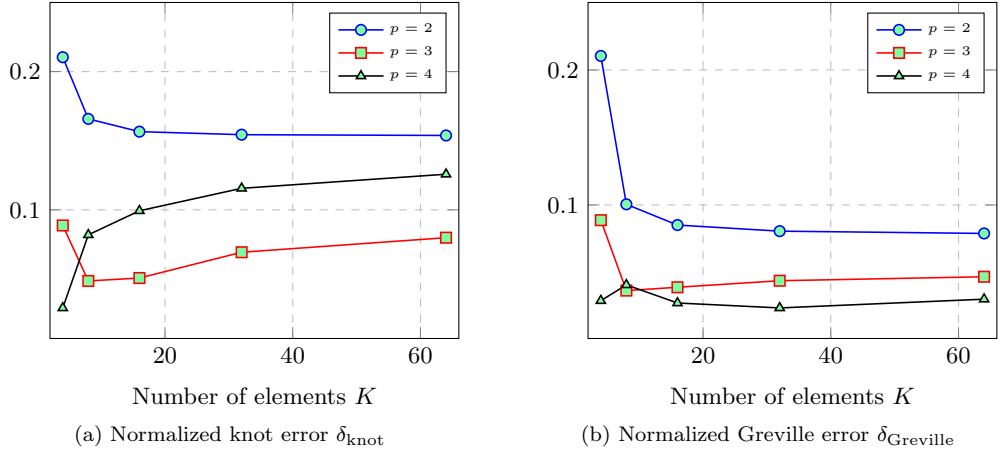

The results in Figure~\ref{fig:smoothedknots} suggest that this simple knot smoothing produces reasonable approximations to optimal knot positions.  Further numerical experiments show that these approximations do not worsen as $K$ increases.  Let $\widehat{\xi}_i$, $\widehat{\tau}_i$ denote knot locations and Greville abscissae for a uniform open knot vector.  Figure~\ref{fig:smoothedknoterror} shows for $p = 2, 3, 4$ the maximum difference between positions of optimal and smoothed knots and Greville abscissae
\[
\delta_{\rm knot} = \frac{\underset{i}{\max} \LRb{\tilde{\xi}_i - \xi^{\rm opt}_i}}{\underset{i}{\max} \LRb{\widehat{\xi}_i - \xi^{\rm opt}_i}}, \qquad \delta_{\rm Greville} = \frac{\underset{i}{\max} \LRb{\tilde{\tau}_i - \tau^{\rm opt}_i}}{\underset{i}{\max} \LRb{\widehat{\tau}_i - \tau^{\rm opt}_i}}, 
\]
where we have normalized by the maximum difference between optimal and uniform knot and Greville positions.  The smoothed knot positions are computed by iterating (\ref{eq:knotiter}) until the change in the knot vector $\nor{\tilde{\bm{\Xi}}^{k} - \tilde{\bm{\Xi}}^{k-1}}_2$ is less than $10^{-8}$.  We observe that the normalized difference between smoothed and optimal knot positions approaches a fixed value around $.1$ as $K$ increases, with the asymptotic value depending irregularly on $p$.  In contrast, the normalized difference between Greville abscissae for smoothed knots and Greville abscissae for optimal knots approaches a $p$-dependent value as $K$ increases, with this asymptotic value decreasing as $p$ increases.  These results quantify how much more accurately smoothed knot vectors approximate optimal knot distributions compared to uniform open knot vectors.  

Finally, we note that while the optimal knot positions minimize the $n$-width of spline spaces over the $L^2$ unit ball of $C^{r-1}$ functions, this does not imply that their approximation properties are uniformly better than those of spline spaces with equispaced knots.  Numerical experiments in Section~\ref{sec:approx} illustrate that spline spaces with optimal knot positions are slightly less accurate than spline spaces with uniform knot vectors in approximating non-oscillatory functions, but are more efficient for the approximation of oscillatory functions.  We will also show that optimal knot positions possess additional advantages over uniform open knot vectors in Section~\ref{sec:cfl}, such as a larger maximum stable timestep and improved approximation properties in the presence of non-affine geometric mappings.   

\section{Method of lines discretization in time}
\label{sec:mol}

\reviewerTwo{In this section, we introduce the semi-discrete system which results from a finite element discretization of time-dependent variational formulations such as (\ref{eq:varformadvec}), (\ref{eq:varformwave1}), and (\ref{eq:varformwave2}).}  Let $\bm{U},\bm{V}$ denote vectors of degrees of freedom for $U,V$, and define the global mass and finite element matrices $\bm{M}_h, \bm{S}_h$ such that 
\[
\bm{V}^T\bm{M}_h\bm{U} = \LRp{ U,V}_{L^2\LRp{\Omega}}, \qquad \bm{V}^T\bm{S}_h\bm{U} = a(U,V),
\]
where $U,V$ are understood to represent either scalar or group variables whose components lie in $V_h$.  
Under a method of lines discretization, a PDE is transformed to a system of ordinary differential equations (ODEs) 
\begin{equation}
\td{\bm{U}}{t} = \bm{A}_h \bm{U}, \qquad \bm{A}_h = -\bm{M}_h^{-1}\bm{S}_h
\label{eq:mol}
\end{equation}
\reviewerTwo{for first-order (in time) formulations and
\begin{equation}
\frac{d^2\bm{U}}{dt^2} = \bm{A}_h \bm{U}, \qquad \bm{A}_h = -\bm{M}_h^{-1}\bm{S}_h
\label{eq:mol_second}
\end{equation} for second-order formulations} by discretizing in space.  The resulting ODE system can then be solved with explicit time stepping for some step size $\Delta t$.  

The following sections discuss efficient strategies for solving (\ref{eq:mol}), focusing on the inversion of $\bm{M}_h^{-1}$ and the choice of timestep $\Delta t$.  Efficient implementations of method of lines solvers typically apply $\bm{A}_h$ in a matrix-free fashion.  However, this requires the computation of $\bm{M}_h^{-1}$ for each application of $\bm{A}_h$.  We discuss in Section~\ref{sec:wadg} the efficient inversion of $\bm{M}_h^{-1}$ using a weight-adjusted approximation to the mass matrix.  Under explicit time stepping, $\Delta t$ must also be chosen such that the eigenvalues of $\Delta t\bm{A}_h$ lie within the stability region of the chosen scheme.  %If an explicit time-stepper is used, the timestep must be sufficiently small to satisfy a stability condition.  
We discuss this condition further in Section~\ref{sec:cfl}.

%The spatial formulations defining $\bm{S}_h$ are constructed such that energy stability holds under inexact quadrature.  As such, the eigenvalues of $\bm{A}_h$ lie in the left half plane

%In order to ensure the eigenvalues of $\bm{A}_h$ lie in the left half plane under inexact quadrature, we adopt a discretization which is \textit{a-priori} energy stable for any arbitrary quadrature with positive weights.  The matrix $\bm{S}_h$ can be written as the sum of skew-symmetric and symmetric parts.   If the symmetric part $\frac{1}{2}\LRp{\bm{S}_h+\bm{S}_h^T}$ is positive semi-definite, all eigenvalues lie in the left half plane and the semi-discrete problem is energy stable \cite{chan2016short}.  We construct a formulation which is stable under inexact numerical integration by splitting $a(u,v)$ into symmetric positive semi-definite and skew-symmetric parts
%\[
%a(u,v) = a_{\rm sym}(u,v) + a_{\rm skew}(u,v), \qquad a_{\rm sym}(u,v) = a_{\rm sym}(v,u), \qquad a_{\rm skew}(u,v) = -a_{\rm skew}(v,u).
%\]
%For advection, integration by parts gives that 
%\[
%a_{\rm sym}(u,v) = \frac{1}{2}\int_{\partial \Omega} \LRb{\bm{\beta}_n} uv\diff{\bm{{x}}}, \qquad {a_{\rm skew}(u,v)} = \frac{1}{2}\int_{\Omega} \LRp{-u\bm{\beta}\cdot \Grad v + \bm{\beta}\cdot \Grad u v}\diff{\bm{{x}}}.
%\]
%Similar approaches have been applied to ensure discrete energy stability for wave propagation problems in a variety of different settings \cite{warburton2013low, chan2015gpu, banks2016galerkin, chan2016weight2, chan2017weight}.  

\subsection{Efficient inversion of mass matrices}
\label{sec:wadg}

Recall that $\bm{\Phi}_k: \Dhat \rightarrow D^k$ denotes the mapping from the reference patch $\Dhat$ to the physical patch $D^k$, such that $\bm{\Phi}_k(\widehat{\bm{x}}) = \bm{x} \in D^k$. Then, because solutions are approximated independently over each patch, the global mass matrix $\bm{M}_h$ is block diagonal, with each block corresponding to a patch mass matrix.  The entries of each patch mass matrix are then 
\[
\LRp{\bm{M}}_{ij} =\reviewerTwo{\int_{\Omega} B^p_i(\bm{\Phi}^{-1}_k\LRp{\bm{x}})B^p_j(\bm{\Phi}^{-1}_k\LRp{\bm{x}})\diff{{\bm{x}}}} = \int_{\Dhat} B^p_i(\bm{\widehat{x}})B^p_j(\bm{\widehat{x}}) J \diff{\widehat{\bm{x}}}.  
\]
If $\bm{\Phi}_k$ is affine, then $J$  (the determinant of the Jacobian of the mapping) is constant, and (assuming a tensor product basis) the mass matrix satisfies
\begin{align*}
\bm{M} &= \widehat{\bm{M}} J = J \LRp{\widehat{\bm{M}}_{\rm 1D} \otimes \ldots \otimes \widehat{\bm{M}}_{\rm 1D}}\\
\LRp{\widehat{\bm{M}}_{\rm 1D}}_{ij} &= \int_{-1}^1 B^p_i({\widehat{x}})B^p_j({\widehat{x}}) \diff{\widehat{{x}}}.  
\end{align*}
Thus, for affine transformations, the inverse of the mass matrix possesses a Kronecker product structure
\[
\bm{M}^{-1} = \LRp{\widehat{\bm{M}}J}^{-1}  = \frac{1}{J} \LRp{\widehat{\bm{M}}^{-1}_{\rm 1D} \otimes \ldots \otimes \widehat{\bm{M}}^{-1}_{\rm 1D}},
\]
which greatly reduces the cost of application in higher dimensions.  However, for non-affine mappings, $J$ is not constant, and the Kronecker product structure of the inverse mass matrix is lost.  

Several approaches have been proposed to address this difficulty in the inversion of NURBS mass matrices.  Spline collocation methods for explicit dynamics circumvent the construction of weighted $L^2$ mass matrices over each patch \cite{auricchio2012isogeometric}.  These formulations discretize the underlying PDE in its strong form, and require only the inversion of an un-weighted generalized Vandermonde matrix, which always admits a Kronecker product structure.  However, in contrast to methods based the Galerkin framework, it can be difficult to prove energy stability of the resulting discretizations.  

Another approach is to replace the exact inversion of the mass matrix with an iterative solver, using a preconditioner possessing a tensor product structure \cite{gao2014fast}.  However, this adds significant complexity and computational cost to solvers based on explicit timestepping, and the convergence of the iteration can depend on the specific geometry.  A more subtle concern associated with preconditioned iterative solvers is energy stability.  The use of an exact mass matrix induces a semi-discrete $L^2$ energy stability, and replacing the mass matrix with a symmetric positive-definite approximation maintains energy stability but changes the norm in which energy is measured.  As noted in \cite{wathen2009chebyshev}, Krylov methods result in a non-linear approximation of the mass matrix inverse, which depends on the initial guess and right hand side vector.  Consequentially, a discrete energy stability no longer holds, as an iterative scheme effectively inverts a slightly different mass matrix each time-step.  

A related iterative approach is the combination of lumped spline mass matrices with predictor multicorrector iterations \cite{auricchio2012isogeometric}.  $C^0$ finite element discretizations can employ ``lumped'' diagonal approximations to the mass matrix, giving rise to the spectral element method when lumping is performed through collocation at Gauss-Lobatto quadrature nodes \cite{patera1984spectral}.  %The relationship of collocation nodes to quadrature also provides a way to incorporate curved geometries into spectral element discretizations in an provably energy stable manner \cite{kopriva2006metric, kopriva2010quadrature}.  
However, the extension of mass lumping to NURBS is not straightforward, as collocation nodes (Greville abscissae) can not be used to construct a sufficiently accurate quadrature for spline spaces, and alternative row-sum mass lumping techniques have been observed to reduce accuracy to second order \cite{cottrell2006isogeometric}.  However, high order accuracy can be restored at the cost of additional corrector iterations \cite{auricchio2012isogeometric}.  These predictor-multicorrector iterations were recently shown to be equivalent to a linear operator, thus avoiding problems of non-linearity present for Krylov iterations \cite{evans2017explicit}.  

We seek an alternative approach to restoring the Kronecker product structure by replacing the weighted $L^2$ mass matrix with a \textit{weight-adjusted} approximation.  Let $B^p_i(\widehat{\bm{x}})$ denote a basis for the $d$-dimensional spline space over the reference patch.  Then, the mass matrix may be approximated by a weight-adjusted mass matrix
\[
\bm{M} \approx \widehat{\bm{M}} \bm{M}^{-1}_{1/J} \widehat{\bm{M}}, \qquad \reviewerTwo{ \LRp{\bm{M}_{1/J}}_{ij} }= \int_{\Dhat} B^p_j(\widehat{\bm{x}})B^p_i(\widehat{\bm{x}}) \frac{1}{J} \diff{\widehat{\bm{x}}},
\]
such that the inverse of the weight-adjusted mass matrix is 
\begin{align*}
\bm{M}^{-1} \approx \widehat{\bm{M}}^{-1} \bm{M}_{1/J} \widehat{\bm{M}}^{-1}.
\end{align*}
This approximation requires only the application of a weighted mass matrix $\bm{M}_{1/J}$ and the inversion of reference mass matrices $\widehat{\bm{M}}^{-1}$.  
\reviewerTwo{Assuming a fixed geometry, the matrix $\bm{M}_{1/J}$ is constant in time, and can thus be pre-computed and reused at each time-step.  
%Moreover, assuming a Cholesky factorization of $\widehat{\bm{M}} = \bm{L}\bm{L}^T$, the inverse of the reference mass matrix is $\widehat{\bm{M}}^{-1} = \bm{L}^{-T}\bm{L}^{-1}$, and the application of the weight-adjusted mass matrix inverse can be written as 
%\[
%\widehat{\bm{M}}^{-1} \bm{M}_{1/J} \widehat{\bm{M}}^{-1} = \bm{L}^{-T}\bm{L}^{-1}\bm{M}_{1/J}\bm{L}^{-T}\bm{L}^{-1} = \bm{L}^{-T}\bm{L}^{-1}\bm{M}_{1/J}\bm{L}^{-T}\bm{L}^{-1}.
%\]
}
Alternatively, the application of $\bm{M}_{1/J}$ can be performed efficiently in a matrix-free fashion using quadrature and the tensor-product structure of the spline space, while the inversion of $\widehat{\bm{M}}$ can be made efficient by taking advantage of the Kronecker structure present on the reference patch
\[
\widehat{\bm{M}}^{-1} = \widehat{\bm{M}}_{\rm 1D}^{-1} \otimes \ldots  \widehat{\bm{M}}_{\rm 1D}^{-1}.  
\]
The computational cost of this approach can be further reduced by fusing products of reference element matrices to yield an implementation in terms of ``projection'' and ``interpolation'' matrices, as described \cite{chan2016weight1, chan2016weight2}.  
%\reviewerTwo{We note that the relative efficiency of each strategy will likely depend on the choice of computational architecture used.  For example, the low storage cost of the matrix-free application is attractive for accelerator architecture such as Graphics Processing Units, which  \cite{chan2017weight} }

For sufficiently regular $J$, this approximation is high order accurate \cite{chan2016weight1,chan2016weight2}.  Additionally, since the weight-adjusted approximation is symmetric-positive definite, it preserves the energy stability of the underlying Galerkin formulation.  Moreover, the weighed and weight-adjusted mass matrices induce discrete $L^2$-equivalent norms with the same equivalence constants \cite{chan2016weight1}.  As a result, the substitution of the weight adjusted mass matrix does not decrease the maximum stable time-step.  

The matrix-free application of the weight-adjusted mass matrix inverse also requires the specification of a sufficiently accurate quadrature.  In this work, we have used a $(p+1)$ point Gaussian quadrature over each knot interval, resulting in composite quadrature rules with a relatively large number of points.  Future work will investigate optimal and reduced quadrature rules tailored to splines \cite{hughes2010efficient, auricchio2012simple, schillinger2014reduced, hiemstra2017optimal}, which greatly reduce the number of quadrature points relative to the number of degrees of freedom in a single patch.  Another method to reduce costs is the use of \reviewerTwo{``weighted quadrature'' methods for applying NURBS mass matrices \cite{calabro2016fast}.  However, we note that energy stability can be lost under this approach, since the curvilinear mass matrices computed using weighted quadrature are not necessarily symmetric}.  

Finally, we note that, while the efficient computational implementation of the proposed methods is outside of the scope of this work, different computational implementations may be preferable depending on computational parameters.  For example, the inverse mass matrix may be pre-multiplied into derivative matrices, resulting in dense size $p+K$ operators in each dimension.  This approach may be useful if $K = O(p)$, and is trivially parallelizable down to individual degrees of freedom.  However, for $K \gg p$, the matrices $\bm{S}_i$ and $\bm{M}$ possess a sparse banded structure, which can be exploited to improve computational efficiency on some architectures.  A Cholesky factorization may be used to apply $\bm{M}^{-1}$ in linear complexity with respect to the size of $\bm{M}$, though the constant will depend on $p$ through the matrix bandwidth.  A disadvantage of this approach is that the triangular solve necessary in Cholesky is largely sequential over each patch.\footnote{\reviewerTwo{We note that, in higher dimensions, applying $\widehat{\bm{M}}^{-1}$ in a Kronecker product fashion involves multiple one-dimensional mass matrix inversions.  It is possible to offset the sequential nature of the Cholesky backsolve by parallelizing over each of these one-dimensional mass inversions.  However, this approach also requires significantly more storage, which can be problematic on many-core architectures such as Graphics Processing Units \cite{chan2017weight}.}}

\subsubsection{NURBS basis functions and mass matrices}
\label{sec:wadgnurbs}

\reviewerTwo{The use of rational NURBS basis functions has been popularized in isogeometric analysis \cite{hughes2005isogeometric}.  While the approximate weight-adjusted mass matrix inverse described in Section~\ref{sec:wadg} is formulated for mapped B-spline basis functions, it is straightforward to extend the approach to approximate the inverse of NURBS mass matrices.  

NURBS basis functions are rational B-splines, defined as follows 
\begin{equation*}
R^p_i(\widehat{\bm{x}}) = \frac{B^p_i(\widehat{\bm{x}}) w_i}{\sum_{j=1}^{p+K}B^p_j(\widehat{\bm{x}}) w_j} = \frac{B^p_i(\widehat{\bm{x}}) w_i}{w_R(\widehat{\bm{x}})},
\end{equation*}
where we have introduced the rational weight $w_R(\widehat{\bm{x}}) = \sum_{j=1}^{p+K}B^p_j(\widehat{\bm{x}}) w_j$.  We note that the weight $w_R$ varies from patch to patch.  The NURBS mass matrix $\bm{M}_R$ for a curvilinear geometric patch $D^k$ is computed by mapping/pulling back to the reference element
\[
(\bm{M}_R)_{ij} = \int_{\Dhat} R^p_i(\widehat{\bm{x}}) R^p_j(\widehat{\bm{x}}) J\diff{\widehat{\bm{x}}} = \int_{\Dhat} B^p_i(\widehat{\bm{x}}) B^p_j(\widehat{\bm{x}}) \frac{J}{w_R^2}\diff{\widehat{\bm{x}}}.  
\]
Thus, the NURBS mass matrix $\bm{M}_R$ can be interpreted as a weighted B-spline mass matrix with weight $J/w_R^2$, and can be approximated by a weight-adjusted mass matrix inverse 
\[
\bm{M}_R^{-1} \approx \widehat{\bm{M}}^{-1} \bm{M}_{w_R^2/J}  \widehat{\bm{M}}^{-1}, \qquad  \LRp{\bm{M}_{w_R^2/J}}_{ij} = \int_{\Dhat} B^p_i(\widehat{\bm{x}})B^p_j(\widehat{\bm{x}}) \frac{w_R^2}{J}\diff{\widehat{\bm{x}}}.
\]
Approximation results in \cite{bazilevs2006isogeometric, evans2013explicit} imply that the weight-adjusted approximation to the NURBS mass matrix will be high order accurate if $w_R^2/J$ is sufficiently regular.  For example, the weight-adjusted approximation of $\bm{M}_R^{-1}$ is expected to maintain high order accuracy if $J, w_R$ vary smoothly and $w_R$ is bounded away from zero.  
}

%\note{Equivalence of WADG with weighted $L^2$ norm:
%\[
%C_1 \nor{u}_{T^{-1}_{1/w}} \leq \nor{u}_w \leq C_2 \nor{u}_{T^{-1}_{1/w}}. 
%\]
%Can I prove that $C_1 = 1$ and bound $C_2$ based on $w(x)$?
%
%We first note that $C_1, C_2$ are extrema of generalized Rayleigh quotients, or minimum and maximum generalized eigenvalues.  Assume $w(x) \in P^M$, $M < N$, or let integrals be computed with quadrature exact for degree $2N+M$ polynomials (these conditions are equivalent), then $N-M$ eigenvalues are unity - note that for $wu \in P^N$, 
%\[
%\nor{u}_{T^{-1}_{1/w}}^2 = \LRp{T_{1/w}^{-1} u, u} = \LRp{\frac{1}{w}T_{1/w}^{-1} u, wu} = \LRp{u,wu} = \nor{u}_w^2.
%\]
%This implies $u \in P^{N-M}$, from which we can construct $N-M$ weighted-orthogonal polynomials.  These are generalized eigenfunctions.  
%%Max extremal eigenvalues are $\Pi_N\LRp{P_N / w}$ where $P_N$ are $L^2$ orthogonal bases.  
%}

\subsection{Stable timestep restrictions}
\label{sec:cfl}
%Let $\bm{U},\bm{V}$ denote vectors of degrees of freedom for $u,v \in V_h$, and define the mass and advection matrices $\bm{M}, \bm{S}_h$ such that $\bm{V}^T\bm{M}\bm{U} = \int_\Omega u v$ and $\bm{V}^T\bm{S}\bm{U} = a(u,v)$.  Under the method of lines, a PDE is transformed to a system of ordinary differential equations (ODEs) 
%\[
%\td{\bm{U}}{t} = \bm{A}_h \bm{U}, \qquad \bm{A}_h = \bm{M}^{-1}\bm{S}_h
%\]
%by discretizing in space.  This system may then be solved using explicit time stepping with some timestep $\Delta t$, such that the eigenvalues of $\Delta t\bm{A}_h$ lie within the stability region of the chosen scheme.  
%If an explicit time-stepper is used, a stability restriction on the timestep must also be imposed.  

\reviewerTwo{For first-order systems,} a rough estimate of the maximum stable timestep is that it is inversely proportional to the spectral radius $\rho(\bm{A}_h)$, which can be estimated by bounding the real and imaginary parts of eigenvalues of $\bm{A}_h$.  The real and imaginary parts of the spectra of $\bm{A}_h$ can be estimated using the symmetric and skew-symmetric parts of $\bm{A}_h$ and theorems from Bendixson and Hirsch \cite{marcus1992survey}, as is done in \cite{chan2015gpu}.  
\[
\max_j \reviewerOne{\LRb{{\rm Re}\LRp{\lambda_j}}} \leq \rho\LRp{\frac{\bm{A}_h+\bm{A}_h^T}{2}}, \qquad \max_j \reviewerOne{\LRb{{\rm Im}\LRp{\lambda_j}}} \leq \rho\LRp{\frac{\bm{A}_h-\bm{A}_h^T}{2}}
\]
%Let $\bm{U},\bm{V}$ denote vectors of degrees of freedom for $u,v \in V_h$.  Since $\bm{A}_h = \bm{M}^{-1}\bm{S}_h$ where $\bm{V}^T\bm{S}\bm{U} = a(u,v)$, 
We estimate the spectral radius of the symmetric part of $\bm{A}_h$ through a generalized Rayleigh quotient
\[
\max_j \reviewerOne{\LRb{{\rm Re}\LRp{\lambda_j}}} \leq \sup_{u \in V_h} \frac{\LRb{a_{\rm sym}(u,u)}}{\nor{u}^2_{\L}}, \qquad a_{\rm sym}(u,v) = \frac{a(u,v) + a(v,u)}{2}.  
\]
The spectral radius of the skew-symmetric part \reviewerTwo{$\rho\LRp{\frac{\bm{A}_h-\bm{A}_h^T}{2}}$} can be bounded using properties of normal matrices \cite{chan2015gpu}; however, it is more straightforward to estimate this term by noting that \reviewerTwo{$\frac{i}{2} \LRp{\bm{A}_h-\bm{A}_h^T}$} is Hermitian and using the generalized Rayleigh quotient.   

\reviewerOne{For illustration, we will derive these bounds for scalar advection.}  We define the skew-symmetric part $a_{\rm skew}(u,v)$ to be sesquilinear such that 
\[
a_{\rm skew}(u,v) = \frac{1}{2}\int_{\Omega} \LRp{-u\bm{\beta}\cdot \overline{\Grad v} + \bm{\beta}\cdot \Grad u \overline{v}}\diff{{\bm{x}}}, \qquad a_{\rm skew}(u,vi) = \bm{V}^T\frac{i}{2} \LRp{\bm{A}_h-\bm{A}_h}\bm{U} = a_{\rm skew}(v,ui).
\]
Then, the Bendixson-Hirsch theorems give that
\[
\max_j \reviewerOne{\LRb{{\rm Im}\LRp{\lambda_j}}} \leq \sup_{u \in V_h} \frac{\LRb{a_{\rm skew}(u,ui)}}{\nor{u}^2_{\L}}, \qquad \LRb{a_{\rm skew}(u,ui)} = \int_{\Omega} \LRp{u\bm{\beta}\cdot \Grad u}\diff{{\bm{x}}} \leq \nor{\bm{\beta}}_{L^\infty}\nor{u}_{\L}\nor{ \Grad u}_{\L}.
\]
%For quasi-uniform patches, 
%\[
%\nor{ \Grad u}_{\L} \leq \nor{\bm{G}}_{L^{\infty}} \nor{\widehat{\Grad}u}_{\L} \leq C \frac{\nor{\widehat{\Grad}u}_{\L}}{H} \leq C\nor{\widehat{\Grad}u}_{L^2\LRp{\widehat{D}}} \leq C C_I \nor{u}_{L^2\LRp{\widehat{D}}} \leq C C_I \nor{u}_{\L}.
%\]
Mapping the formulation for advection to and from the reference patch $\Dhat$ and applying inverse and trace inequalities then yields
\begin{align*}
\LRb{a_{\rm sym}(u,u)} 
&\leq \max_k \LRp{\frac{\nor{\bm{\beta}}_{L^\infty}}{2} C_T \nor{J^s}_{L^\infty}\nor{J^{-1}}_{L^\infty} \nor{u}^2_{L^2\LRp{D^k}}} \\
\LRb{a_{\rm skew}(u,u)} &\leq \max_k \LRp{{\nor{\bm{\beta}}_{L^\infty}} C_I \nor{J\bm{G}}_{L^\infty}\nor{J^{-1}}_{L^\infty} \nor{u}^2_{L^2\LRp{D^k}}}.
\end{align*}
These bounds together imply that 
\begin{align*}
\rho\LRp{\bm{A}_h} &\leq \max_k \LRp{\nor{\bm{\beta}}_{L^\infty}\max\LRc{\frac{1}{2} C_T \nor{J^s}_{L^\infty}, C_I \nor{J\bm{G}}_{L^\infty}}\nor{J^{-1}}_{L^\infty}}.
%\max_j \reviewerOne{\LRb{{\rm Im}\LRp{\lambda_j}}} &\leq \frac{\nor{\bm{\beta}}_{L^\infty}}{2} .
\end{align*}

In one dimension, on a uniform domain subdivided into patches of size \reviewerOne{$H$}, this reduces to the expression
\[
\rho\LRp{\bm{A}_h} \leq \max_k  \LRp{\frac{\nor{\bm{\beta}}_{L^\infty}}{\reviewerOne{H}} C_p}, \qquad \qquad \reviewerOne{C_p = \max\LRc{\frac{1}{2} C_T, C_I}},
\]
which can be used to determine a stable timestep by setting $\Delta t = 1/\rho\LRp{\bm{A}_h}$.  

For the acoustic wave equation, a similar argument can be applied to bound the spectral radius via
\[
\rho\LRp{\bm{A}_h} \leq \max_k\LRp{ \frac{\nor{c}_{L^\infty}}{\reviewerOne{H}} C_p},
\]
where $c$ is the wavespeed.  Setting $\Delta t = O(1/\rho\LRp{\bm{A}_h})$ then yields an estimate of the maximum stable timestep 
\[
\Delta t \leq C_t\frac{\reviewerOne{H}}{C_p\nor{c}_{L^\infty} },
\]
\reviewerOne{where $C_t$ is some constant which depends on the stability region of the explicit time-stepper.}
Thus, for high order finite element discretizations, the stable timestep is determined not only by the \reviewerOne{patch} size $\reviewerOne{H}$ and maximum wavespeed ($\nor{\bm{\beta}}_{L^{\infty}}$ or $\nor{c}_{L^\infty}$, but also by trace and inverse inequality constants $C_T, C_I$ \reviewerOne{which depend on the number of elements $K$ per patch and the order $p$}.

\reviewerTwo{For second-order systems, a rough estimate of the maximum stable timestep is that it is inversely proportional to the square root of the spectral radius $\rho(\bm{A}_h)$.  For the second-order order formulation of the acoustic wave equation, the spectral radius $\rho(\bm{A}_h)$ is equivalent to the Rayleigh quotient
\begin{align*}
\rho(\bm{A}_h) = \sup_{u \in V_h} \frac{a(u,u)}{\nor{u}^2_{\L}}
\end{align*}
provided the constant $\tau$ is sufficiently large to guarantee positive-definiteness of $\bm{A}_h$.  If $\tau$ scales linearly with the trace constant $C_T$, then an application of inverse and trace inequalities yields
\begin{align*}
\rho(\bm{A}_h) \leq C_{\tau} \frac{\left( C_T^2 + C_I^2 \right) \nor{c}^2_{L^\infty}}{H^2}
\end{align*}
where $C_{\tau}$ is a constant of proportionality between $\tau$ and $C_T$.  Thus, we arrive at a maximum stable timestep estimate of
\[
\Delta t \leq C_t C_\tau\frac{H}{ \sqrt{C_T^2 + C_I^2} \nor{c}_{L^\infty} },
\]
where again $C_t$ is some constant which depends on the stability region of the explicit time-stepper.}

\subsubsection{Constants in trace and inverse inequalities for spline spaces}

It is possible to derive explicit formulas for the constants in trace inequalities for polynomials \cite{warburton2003constants, hillewaert2011sharp, chan2015hp}.  While one can repeat this process for inverse inequalities or spline spaces, it is more difficult to derive tight, explicit expressions for such constants in terms of $p$.  However, because these constants are given as extrema of generalized Rayleigh quotients, one can also compute constants in trace or inverse inequality through the solution of a generalized eigenvalue problem over the reference patch.  In one dimension, these problems are given as
\begin{equation}
\widehat{\bm{M}}^fu = \lambda_T\widehat{\bm{M}}u, \qquad \widehat{\bm{K}}u = \lambda_I\widehat{\bm{M}}u,
\label{eq:eigconst}
\end{equation}
where $\widehat{\bm{M}}^f$ is the reference boundary mass matrix
\[
\LRp{\widehat{\bm{M}}^f}_{ij} = B^p_i(-1) B^p_j(-1) + B^p_i(1) B^p_j(1),
\]
and $\widehat{\bm{K}}$ is the reference stiffness matrix
\[
\LRp{\widehat{\bm{K}}}_{ij} = \int_{-1}^1 \pd{B^p_i}{\widehat{x}}\pd{B^p_j}{\widehat{x}} \diff{\widehat{x}}.
\]
{ Then, the constant $C_T$ is the largest generalized eigenvalue of the boundary mass matrix, and the constant $C_I$ is the square root of the largest generalized eigenvalue of the stiffness matrix.}

Assuming a tensor product construction of $V_h$, constants in higher dimensional trace inequalities can be determined in terms of the constants in one-dimensional trace inequalities.  For example, in two dimensions, the mass, face mass, and stiffness matrices are given as
\begin{align*}
\widehat{\bm{M}} &= \widehat{\bm{M}}_{\rm 1D}\otimes \widehat{\bm{M}}_{\rm 1D}, \qquad \widehat{\bm{M}}^f = \widehat{\bm{M}}_{\rm 1D}\otimes \widehat{\bm{M}}^f_{\rm 1D}\\
\widehat{\bm{K}} &= \widehat{\bm{M}}_{\rm 1D}\otimes \widehat{\bm{K}}_{\rm 1D} + \widehat{\bm{K}}_{\rm 1D} \otimes \widehat{\bm{M}}_{\rm 1D}.
\end{align*}
Substituting these relations into the generalized eigenvalue problem (\ref{eq:eigconst}) for the trace constant and reducing gives
\[
\bm{I}\otimes \widehat{\bm{M}}^f_{\rm 1D} u = \lambda_T\LRp{\bm{I}\otimes \widehat{\bm{M}}}u.  
\]
Because the spectra of $\bm{I}\otimes \widehat{\bm{M}}^f_{\rm 1D}$ consists only of multiples of one-dimensional eigenvalues $\lambda_T$, constants in higher dimensional trace inequalities are identical to the constants in one spatial dimension.  

In comparison, constants in $d$-dimensional inverse inequalities scale with $\sqrt{d}$, which can be seen by diagonalizing (\ref{eq:eigconst}).  Let $\bm{W}_{\rm 1D}$ denote the coefficient matrix for eigenvectors of (\ref{eq:eigconst}) in one dimension, such that $\bm{W}_{\rm 1D}^* \widehat{\bm{K}}_{\rm 1D} \bm{W}_{\rm 1D} = \bm{\Lambda}_{\rm 1D}$ and $\bm{W}_{\rm 1D}^* \widehat{\bm{M}}_{\rm 1D} \bm{W}_{\rm 1D} = \bm{I}$.  Then, the change of basis $\bm{W}_{\rm 1D}\otimes \bm{W}_{\rm 1D}$ reduces the eigenvalue problem for constants in the two-dimensional inverse inequality to 
\[
\LRp{\bm{I}_{\rm 1D}\otimes \bm{\Lambda}_{\rm 1D} + \bm{\Lambda}_{\rm 1D}\otimes\bm{I}_{\rm 1D}} u = \lambda u,
\]
Let $\mu_j$ denote one-dimensional generalized eigenvalues for the stiffness matrix, and let the two dimensional eigenvalues be indexed by $ij$.  Then, $\lambda_{ij} = \mu_i + \mu_j$ and $\max_{ij} \lambda_{ij} = 2\max_j \mu_j$.  Assuming a indexing by $i_1, \ldots, i_d$ in $d$ dimensions, the result generalizes to higher dimensions as $\max_{i_1,\ldots,i_d} \lambda_{i_1,\ldots,i_d} = d\max_j \mu_j$.  

For polynomials, it can be shown that the constants $C_I,C_T$ depend quadratically on the degree of approximation $p$ (see, for example \cite{warburton2003constants, ozisik2010constants,hillewaert2011sharp}).  For standard $C^0$-continuous { and discontinuous Galerkin discretizations} of the advection equation { using a polynomial approximation over each patch}, it follows that $\rho(\bm{A}_h) = O(p^2 / h)$ { since $H = h$}.  { In contrast, for spline spaces of maximal smoothness, numerical experiments in the following sections and analytic results for special cases \cite{takacs2016approximation} suggest that the constants $C_I, C_T$ grow as $O(pK)$ when the number elements $K$ in a patch is sufficiently large with respect to $p$, and thus $\rho(\bm{A}_h)$ grows as $O(p / h)$.}  In the special case of periodic splines, the stronger result can be shown that $C_I, C_T$ are independent of $p$ \cite{takacs2016approximation}.  

\subsubsection{Numerically computed constants}
\label{sec:numconsts}

We consider spline bases with uniform, optimal, and smoothed knot vectors, though we limit our experiments using optimal knot vectors to $p \leq 5$ due to issues of numerical stability \reviewerOne{in the computation of optimal knot locations}.  For each basis, we compute constants $C_I,C_T$ in trace and inverse inequalities through the solution of a generalized eigenvalue problem.  The number of elements is increased with $p$ as $K = \lceil p/2 \rceil, p, 2p$.  %We consider two scenarios for the number of elements $K$: in the first, we fix the number of elements $K = 4, 8, 16$ while increasing the degree $p$.  In the second scenario, we increase the number of elements with $p$ as $K = \lceil p/2 \rceil, p, 2p$.  

\begin{figure}
\centering
\subfloat[$K = \lceil \frac{p}{2}\rceil$]{
\begin{tikzpicture}
\begin{axis}[
width=.35\textwidth,
    xlabel={Degree of approximation $p$},   
        ylabel={$C_I / K$},
    xmin=.5, xmax=12.5,
    ymax=0,ymax=90,        
    legend pos=north west, legend cell align=left, legend style={font=\tiny},	
    xmajorgrids=true,  ymajorgrids=true, grid style=dashed,
] 

\addplot[color=blue,mark=*,semithick, mark options={fill=markercolor}]
coordinates{(1,2)(2,4.5)(3,8)(4,12.5)(5,18)(6,24.5)(7,32)(8,40.5)(9,50)(10,60.5)(11,72)(12,84.5)};
\addplot[color=red,mark=square*,semithick, mark options={fill=markercolor}]
coordinates{(1,3)(2,6)(3,6.5)(4,9.5)(5,11.0172)(6,14.178)(7,16.171)(8,19.4951)(9,21.7992)(10,25.2711)(11,27.8076)(12,31.4113)};
\addplot[color=black,mark=triangle*,semithick, mark options={fill=markercolor}]
coordinates{(1,3)(2,6)(3,6.5)(4,9.5)(5,10.5251)(6,13.6497)(7,14.7755)(8,17.9829)(9,19.1754)(10,22.4404)(11,23.6889)(12,26.9988)};

\legend{Polynomial $C_I$, Uniform $C_I/K$, Smoothed $C_I/K$ }
\end{axis}
\end{tikzpicture}
}
\subfloat[$K = p$]{
\begin{tikzpicture}
\begin{axis}[
width=.35\textwidth,
    xlabel={Degree of approximation $p$},   
    xmin=.5, xmax=12.5,
    ymax=0,ymax=90,    
    legend pos=north west, legend cell align=left, legend style={font=\tiny},	
    xmajorgrids=true,  ymajorgrids=true, grid style=dashed,
] 

\addplot[color=blue,mark=*,semithick, mark options={fill=markercolor}]
coordinates{(1,2)(2,4.5)(3,8)(4,12.5)(5,18)(6,24.5)(7,32)(8,40.5)(9,50)(10,60.5)(11,72)(12,84.5)};
\addplot[color=red,mark=square*,semithick, mark options={fill=markercolor}]
coordinates{(1,3)(2,4)(3,5.70747)(4,7.72711)(5,9.95857)(6,12.3561)(7,14.8924)(8,17.5488)(9,20.3115)(10,23.1696)(11,26.1146)(12,29.1391)};
\addplot[color=black,mark=triangle*,semithick, mark options={fill=markercolor}]
coordinates{(1,3)(2,4)(3,5.34993)(4,6.8709)(5,8.45922)(6,10.0982)(7,11.7767)(8,13.4862)(9,15.2211)(10,16.9775)(11,18.7525)(12,20.5439)};

\legend{Polynomial $C_I$, Uniform $C_I/K$, Smoothed $C_I/K$ }
\end{axis}
\end{tikzpicture}
}
\subfloat[$K = 2p$]{
\begin{tikzpicture}
\begin{axis}[
width=.35\textwidth,
    xlabel={Degree of approximation $p$},   
    xmin=.5, xmax=12.5,
    ymax=0,ymax=90,    
    legend pos=north west, legend cell align=left, legend style={font=\tiny},	
    xmajorgrids=true,  ymajorgrids=true, grid style=dashed,
] 

\addplot[color=blue,mark=*,semithick, mark options={fill=markercolor}]
coordinates{(1,2)(2,4.5)(3,8)(4,12.5)(5,18)(6,24.5)(7,32)(8,40.5)(9,50)(10,60.5)(11,72)(12,84.5)};
\addplot[color=red,mark=square*,semithick, mark options={fill=markercolor}]
coordinates{(1,2)(2,3.42857)(3,5.27022)(4,7.3491)(5,9.61122)(6,12.0258)(7,14.5716)(8,17.233)(9,19.9977)(10,22.8559)(11,25.7995)(12,28.8217)};
\addplot[color=black,mark=triangle*,semithick, mark options={fill=markercolor}]
coordinates{(1,2)(2,3.15356)(3,4.41024)(4,5.72791)(5,7.09105)(6,8.48919)(7,9.91299)(8,11.3583)(9,12.8221)(10,14.3019)(11,15.796)(12,17.3028)};

\legend{Polynomial $C_I$, Uniform $C_I/K$, Smoothed $C_I/K$ }
\end{axis}
\end{tikzpicture}
}
\caption{Scaled constants in inverse inequalities for splines with various choices of $K$, along with corresponding constants for polynomials. }
\label{fig:iconsts}
\end{figure}
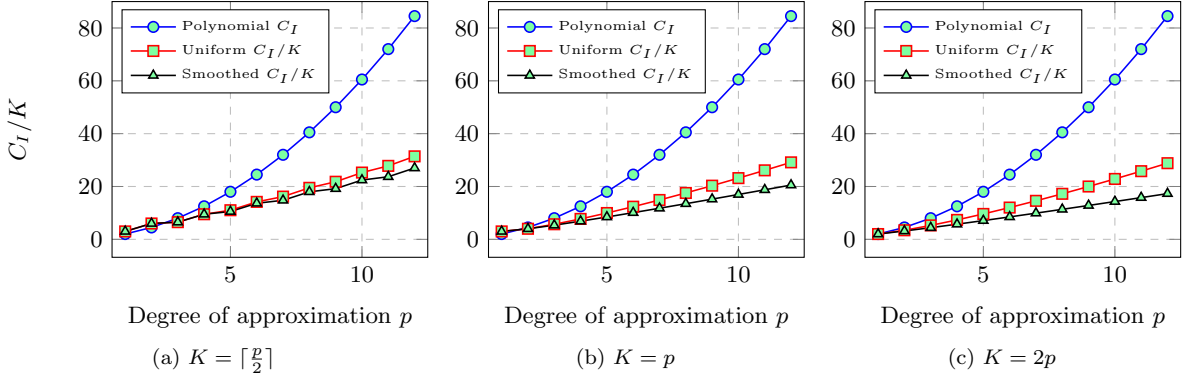

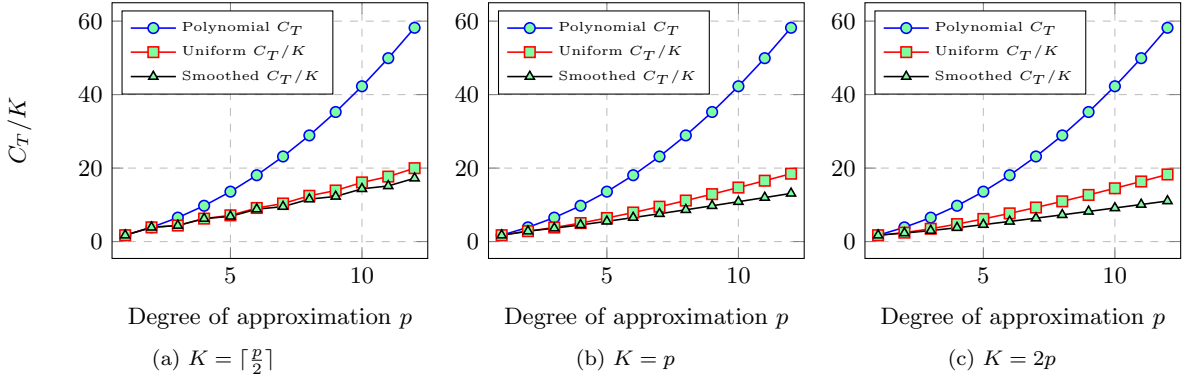
\begin{figure}
\centering
\subfloat[$K = \lceil \frac{p}{2}\rceil$]{
\begin{tikzpicture}
\begin{axis}[
width=.35\textwidth,
    xlabel={Degree of approximation $p$},   
        ylabel={$C_T / K$},
    xmin=.5, xmax=12.5,
    ymax=0,ymax=65,        
    legend pos=north west, legend cell align=left, legend style={font=\tiny},	
    xmajorgrids=true,  ymajorgrids=true, grid style=dashed,
] 

\addplot[color=blue,mark=*,semithick, mark options={fill=markercolor}]
coordinates{(1,1.73205)(2,3.87298)(3,6.5216)(4,9.74981)(5,13.5914)(6,18.0596)(7,23.1598)(8,28.894)(9,35.2634)(10,42.2685)(11,49.9096)(12,58.187)};
\addplot[color=red,mark=square*,semithick, mark options={fill=markercolor}]
coordinates{(1,1.73205)(2,3.87298)(3,4.37379)(4,6.25368)(5,7.1462)(6,9.14314)(7,10.3601)(8,12.4668)(9,13.8967)(10,16.0995)(11,17.6856)(12,19.973)};
\addplot[color=black,mark=triangle*,semithick, mark options={fill=markercolor}]
coordinates{(1,1.73205)(2,3.87298)(3,4.37379)(4,6.25368)(5,6.85572)(6,8.8263)(7,9.51337)(8,11.5444)(9,12.2862)(10,14.3577)(11,15.1409)(12,17.2432)};

\legend{Polynomial $C_T$, Uniform $C_T/K$, Smoothed $C_T/K$ }
\end{axis}
\end{tikzpicture}
}
\subfloat[$K = p$]{
\begin{tikzpicture}
\begin{axis}[
width=.35\textwidth,
    xlabel={Degree of approximation $p$},   
    xmin=.5, xmax=12.5,
    ymax=0,ymax=65,    
    legend pos=north west, legend cell align=left, legend style={font=\tiny},	
    xmajorgrids=true,  ymajorgrids=true, grid style=dashed,
] 

\addplot[color=blue,mark=*,semithick, mark options={fill=markercolor}]
coordinates{(1,1.73205)(2,3.87298)(3,6.5216)(4,9.74981)(5,13.5914)(6,18.0596)(7,23.1598)(8,28.894)(9,35.2634)(10,42.2685)(11,49.9096)(12,58.187)};
\addplot[color=red,mark=square*,semithick, mark options={fill=markercolor}]
coordinates{(1,1.73205)(2,2.8364)(3,3.82963)(4,5.04001)(5,6.41327)(6,7.90642)(7,9.49532)(8,11.165)(9,12.9051)(10,14.7078)(11,16.5672)(12,18.4783)};
\addplot[color=black,mark=triangle*,semithick, mark options={fill=markercolor}]
coordinates{(1,1.73205)(2,2.8364)(3,3.6403)(4,4.54195)(5,5.51171)(6,6.52839)(7,7.57837)(8,8.65295)(9,9.74686)(10,10.8566)(11,11.9797)(12,13.1144)};

\legend{Polynomial $C_T$, Uniform $C_T/K$, Smoothed $C_T/K$ }
\end{axis}
\end{tikzpicture}
}
\subfloat[$K = 2p$]{
\begin{tikzpicture}
\begin{axis}[
width=.35\textwidth,
    xlabel={Degree of approximation $p$},   
    xmin=.5, xmax=12.5,
    ymax=0,ymax=65,    
    legend pos=north west, legend cell align=left, legend style={font=\tiny},	
    xmajorgrids=true,  ymajorgrids=true, grid style=dashed,
] 

\addplot[color=blue,mark=*,semithick, mark options={fill=markercolor}]
coordinates{(1,1.73205)(2,3.87298)(3,6.5216)(4,9.74981)(5,13.5914)(6,18.0596)(7,23.1598)(8,28.894)(9,35.2634)(10,42.2685)(11,49.9096)(12,58.187)};
\addplot[color=red,mark=square*,semithick, mark options={fill=markercolor}]
coordinates{(1,1.73205)(2,2.43668)(3,3.48657)(4,4.75661)(5,6.16327)(6,7.6751)(7,9.27477)(8,10.9506)(9,12.694)(10,14.4983)(11,16.3578)(12,18.2682)};
\addplot[color=black,mark=triangle*,semithick, mark options={fill=markercolor}]
coordinates{(1,1.73205)(2,2.32019)(3,2.9916)(4,3.77273)(5,4.61001)(6,5.4806)(7,6.37331)(8,7.28317)(9,8.20701)(10,9.14265)(11,10.0885)(12,11.0432)};

\legend{Polynomial $C_T$, Uniform $C_T/K$, Smoothed $C_T/K$ }
\end{axis}
\end{tikzpicture}
}
\caption{Scaled constants in trace inequalities for splines with various choices of $K$, along with corresponding constants for polynomials. }
\label{fig:tconsts}
\end{figure}

Figures~\ref{fig:iconsts} and \ref{fig:tconsts} show the growth of $C_I /K, C_T/K$ for one dimensional spline bases, where we have scaled the computed trace and inverse inequality constants by $1/K$ to isolate the growth of $C_I, C_T$ with respect to $p$.  We observe that, for $K \geq p$, the growth of $C_I, C_T$ for all spline spaces is roughly linear in { both} $p$ { and $K$}.  Additionally, while knot smoothing does not decrease the asymptotic rate of growth of $C_I, C_T$, we observe that it results in noticeable reductions in each constant.  Moreover, the reduction in $C_I, C_T$ from knot smoothing increases the larger $K$ is relative to $p$.

For a uniform open knot vector, the distribution of the Greville abscissae clusters near the endpoints of the interval.  Because knot smoothing reduces the clustering of Greville abscissae near the boundaries, the reduction in $C_I, C_T$ is to be expected.  For example, the clustering of stable interpolation points has been identified as contributing to the $O(1/p^2)$ CFL condition in high order pseudospectral methods: for Gauss-Legendre-Lobatto or Chebyshev nodes, the minimum spacing between nodes decreases as $O(1/p^2)$.  Kosloff and Tal-Ezer proposed a mapping of the spatial domain \cite{kosloff1993modified} to address this timestep limitation, and showed that the mapping improved the CFL condition from $O(1/p^2)$ to $O(1/p)$.  However, this mapping also reduces accuracy (especially at lower degrees of approximation \cite{atcheson2013explicit}).  \reviewerOne{By contrast, knot smoothing does not significantly decrease accuracy.  In fact, the approximation power of spline spaces under knot smoothing is actually improved for certain functions.    }

%% trace inequalities

Finally, we examine the difference between $C_I, C_T$ for smoothed and optimal knot vectors.  Table~\ref{tab:consts} shows constants $C_I, C_T$ in inverse and trace inequalities for smoothed and optimal spline spaces of degree $p = 2,\ldots, 5$ with $K = p, 2p$ elements.  We observe that the values of $C_I, C_T$ do not differ greatly between smoothed and optimal knot vectors.  

\begin{table}
\centering
\subfloat[Scaled trace constants \reviewerOne{$C_T/K$}]{
\begin{tabular}{|c||c|c|c|c|c|}
\hline
 & $p = 2$ & $p = 3$ & $p = 4$ & $p = 5$ \\
\hhline{|=|=|=|=|=|=|}
%$K = p$ (uniform knots) &  2.8364  &  3.8296&    5.0400&    6.4133\\
\hline
$K = p$ (smoothed knots) & 2.8364 &   3.6403&    4.5420&    5.5117\\
\hline
$K = p$ (optimal knots) & 2.8364  &  3.6604&    4.5303&    5.4323\\
\hhline{|=|=|=|=|=|=|}
%$K = 2p$ (uniform knots) &  2.4367    &3.4866    &4.7566  &  6.1633\\
\hline
$K = 2p$ (smoothed knots) &  2.3202   & 2.9916&    3.7727&    4.6100\\
\hline
$K = 2p$  (optimal knots)  & 2.3215   & 3.0011   & 3.7515  &  4.5251\\
\hline
\end{tabular}
}\\
\subfloat[Scaled inverse inequality constants \reviewerOne{$C_I/K$}]{
\begin{tabular}{|c||c|c|c|c|c|}
\hline
 & $p = 2$ & $p = 3$ & $p = 4$ & $p = 5$ \\
\hhline{|=|=|=|=|=|=|}
%$K = p$ (uniform knots) & 4.0000   & 5.7075&    7.7271&    9.9586\\
\hline
$K = p$ (smoothed knots) &     4.0000    &5.3499    &6.8709&    8.4592\\
\hline
$K = p$ (optimal knots) & 4.0000    &5.3907&    6.8514&    8.3260\\
\hhline{|=|=|=|=|=|=|}
%$K = 2p$ (uniform knots) & 3.4286    &5.2702    &7.3491    &9.6112 \\
\hline
$K = 2p$ (smoothed knots) & 3.1536   & 4.4102&    5.7279&    7.0910\\
\hline
$K = 2p$ (optimal knots) & 3.1768   & 4.4311   & 5.6895   & 6.9474\\
\hline
\end{tabular}
}
\caption{Comparison of constants in spline trace and inverse inequalities using smoothed and optimal knot vectors.}
\label{tab:consts}
\end{table}

\subsection{On the CFL restriction for piecewise polynomial and spline finite elements: $h$ vs $p$ methods}

\reviewerOne{As mentioned previously, $C^0$ and discontinuous finite element methods result in a stable timestep restriction of $O(h/p^2)$ for a given $h$ and $p$.  Assuming an {$O(pK)$} scaling of inverse and trace constants $C_I, C_T$ for spline approximations, spline discretizations yield a stable time-step restriction of $O(h/p)$.  This $O(h/p)$ stable timestep restriction for splines applies when interpreting spline finite elements as $h$-methods, where $h, K$ are assumed to be independent of $p$ with $K \gg p$.  

However, when interpreting spline discretizations as $p$-methods \cite{babuska1981p}, resolution is achieved by increasing the order rather than by mesh refinement, and the assumption $K \gg p$ may not be appropriate.  Consider { then the case of} a single patch of size $H=O(1)$ { with $K = O(p)$.  Since $h = O(1/K)$ and $C_I, C_T$ scale like $O(pK)$ when $K = O(p)$}, we have that the stable timestep restriction for splines on a single patch is $O(h/p) = O(1/p^2)$.  This analysis shows {spline discretizations with $K = O(p)$ exhibit stable timestep sizes with the same asymptotic dependence on $p$ as those of $C^0$ and discontinuous finite element methods}.  However, we note that this does not account for the fact that for the same stable timestep, spline discretizations achieve an extra factor $O(1/K) = O(1/p)$ of additional mesh resolution compared to a polynomial discretization.  Thus splines can achieve the same stable timestep restriction as $C^0$ finite elements while providing a higher resolution.  

An alternative way to compare the stable timestep restriction for splines and $C^0$ finite elements is in terms of degrees of freedom.  Up to this point, our discussion has related the time-step size to the mesh size $h$ and the polynomial degree $p$, but for the same mesh size $h$ and polynomial degree $p$, a piecewise polynomial finite element basis has a larger number of degrees of freedom $N_{\rm dofs}$ than a corresponding spline basis.  Thus, it may seem appropriate instead to examine the maximum stable timestep in terms of the number of degrees of freedom.  In this direction, we note that $N_{\rm dofs} \sim p/h$ for $C^0$ or discontinuous piecewise polynomial finite elements in the one-dimensional setting, so the maximum stable timestep scales like $1/(pN_{\rm dofs})$.  For one-dimensional spline finite elements, $N_{\rm dofs} \sim 1/h$, so the maximum stable timestep also scales like $1/(pN_{\rm dofs})$.  This analysis seems to suggest spline finite elements do not exhibit any improvement over $C^0$ or discontinuous piecewise polynomial finite elements, at least with respect to timestep size.  However, for the same polynomial degree $p$ and number of degrees of freedom $N_{\rm dofs}$, a spline basis is generally more accurate than a piecewise polynomial finite element basis \cite{evans2009n}.  Thus, to attain the same level of accuracy for a fixed $p$, one can utilize a larger timestep size with splines than with $C^0$ or discontinuous finite elements.  %In fact, the same level of accuracy is typically attained when the mesh sizes are comparable, so an $O(p)$ times larger timestep may be employed for splines versus $C^0$ or discontinuous finite elements.  
}

\section{Numerical experiments}
\label{sec:num}

In the following sections, we present numerical experiments which verify properties of the proposed approximation spaces and discretizations.  All computations use Gaussian quadratures of degree $p$ over each knot interval, and all time-dependent simulations use a 4th order 5-stage low-storage Runge-Kutta method for time integration \cite{carpenter1994fourth} (with the second order formulation rewritten as a first order system).  All experiments use first or second order multi-patch DG formulations for the acoustic wave equation, except for the computation of discrete dispersion relations (which use the multi-patch DG formulation for the advection equation).  

\subsection{One dimensional experiments}

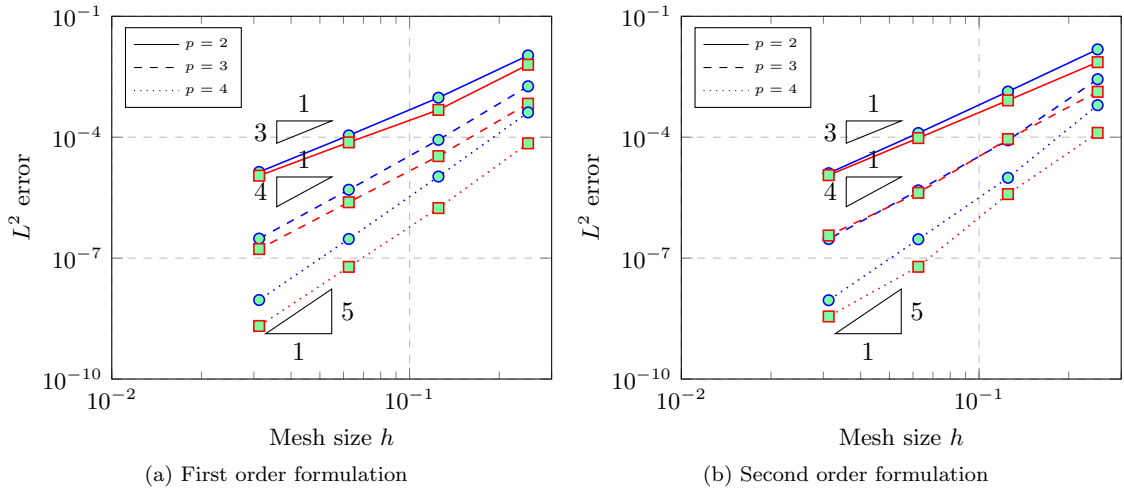
\begin{figure}
\centering
\subfloat[First order formulation]{
\begin{tikzpicture}
\begin{loglogaxis}[
width=.45\textwidth,
    xlabel={Mesh size $h$},
    ylabel={$L^2$ error},
    xmin=1e-2, xmax=.3,
    ymin=1e-10, ymax=.1,    
    legend pos=north west, legend cell align=left, legend style={font=\tiny},	
       xmajorgrids=true,  ymajorgrids=true, grid style=dashed,
    legend entries={$p = 2$, $p=3$, $p=4$},
] 
\pgfplotsset{
cycle list={
{blue, mark=*}, {red,mark=square*},
{blue, mark=*,dashed}, {red,mark=square*,dashed},
{blue, mark=*,dotted}, {red,mark=square*,dotted},
legend style={mark options={no markers}}
}}
\addlegendimage{no markers,black}
\addlegendimage{no markers,black,dashed}
\addlegendimage{no markers,black,dotted}

\addplot+[semithick, mark options={solid,fill=markercolor}]
coordinates{(0.25,0.0106673)(0.125,0.000954815)(0.0625,0.000110851)(0.03125,1.37364e-05)};
\addplot+[semithick, mark options={solid,fill=markercolor}]
coordinates{(0.25,0.00629266)(0.125,0.000476247)(0.0625,7.45318e-05)(0.03125,1.11051e-05)};
\logLogSlopeTriangleFlip{.5}{.125}{.65}{3}{}

    % #1. Relative offset in x direction.
    % #2. Width in x direction, so xA-xB.
    % #3. Relative offset in y direction.
    % #4. Slope d(y)/d(log10(x)).

%    4.0252
\addplot+[semithick, mark options={solid,fill=markercolor}]
coordinates{(0.25,0.00183485)(0.125,8.58892e-05)(0.0625,4.97075e-06)(0.03125,3.0529e-07)};
%    3.8576
\addplot+[semithick, mark options={solid,fill=markercolor}]
coordinates{(0.25,0.000679534)(0.125,3.39478e-05)(0.0625,2.43527e-06)(0.03125,1.67989e-07)};
\logLogSlopeTriangleFlip{.5}{.125}{.475}{4}{}

%    5.0305
\addplot+[semithick, mark options={solid,fill=markercolor}]
coordinates{(0.25,0.000412175)(0.125,1.05591e-05)(0.0625,2.98645e-07)(0.03125,9.13748e-09)};
%    4.8741
\addplot+[semithick, mark options={solid,fill=markercolor}]
coordinates{(0.25,7.04497e-05)(0.125,1.74612e-06)(0.0625,6.0791e-08)(0.03125,2.07289e-09)};
\logLogSlopeTriangle{.5}{.15}{.125}{5}{}

\end{loglogaxis}
\end{tikzpicture}
}
\subfloat[Second order formulation]{
\begin{tikzpicture}
\begin{loglogaxis}[
width=.45\textwidth,
    xlabel={Mesh size $h$}, 
    ylabel={$L^2$ error},
    xmin=1e-2, xmax=.3,
    ymin=1e-10, ymax=.1,    
    legend pos=north west, legend cell align=left, legend style={font=\tiny},	
       xmajorgrids=true,  ymajorgrids=true, grid style=dashed,
    legend entries={$p = 2$, $p=3$, $p=4$},
] 
\pgfplotsset{
cycle list={
{blue, mark=*}, {red,mark=square*},
{blue, mark=*,dashed}, {red,mark=square*,dashed},
{blue, mark=*,dotted}, {red,mark=square*,dotted},
legend style={mark options={no markers}}
}}
\addlegendimage{no markers,black}
\addlegendimage{no markers,black,dashed}
\addlegendimage{no markers,black,dotted}

\addplot+[semithick, mark options={solid,fill=markercolor}]
coordinates{(0.25,0.0150916)(0.125,0.00135624)(0.0625,0.000126749)(0.03125,1.27992e-05)};
\addplot+[semithick, mark options={solid,fill=markercolor}]
coordinates{(0.25,0.00727599)(0.125,0.000820116)(0.0625,9.52161e-05)(0.03125,1.14411e-05)};
\logLogSlopeTriangleFlip{.5}{.125}{.65}{3}{}

    % #1. Relative offset in x direction.
    % #2. Width in x direction, so xA-xB.
    % #3. Relative offset in y direction.
    % #4. Slope d(y)/d(log10(x)).

\addplot+[semithick, mark options={solid,fill=markercolor}]
coordinates{(0.25,0.00274207)(0.125,8.37621e-05)(0.0625,4.79288e-06)(0.03125,2.97136e-07)};
\addplot+[semithick, mark options={solid,fill=markercolor}]
coordinates{(0.25,0.001341)(0.125,8.99204e-05)(0.0625,4.17057e-06)(0.03125,3.6443e-07)};
\logLogSlopeTriangleFlip{.5}{.125}{.475}{4}{}

\addplot+[semithick, mark options={solid,fill=markercolor}]
coordinates{(0.25,0.000623337)(0.125,9.84177e-06)(0.0625,2.93926e-07)(0.03125,9.019e-09)};
\addplot+[semithick, mark options={solid,fill=markercolor}]
coordinates{(0.25,0.000128353)(0.125,3.87921e-06)(0.0625,6.10617e-08)(0.03125,3.57278e-09)};
\logLogSlopeTriangle{.5}{.15}{.125}{5}{}

\end{loglogaxis}
\end{tikzpicture}
}
\caption{One-dimensional convergence of multi-patch DG methods under knot insertion and patch refinement.  Uniform knots are used for all experiments.  Results for refinement by knot insertion are denoted by circles, while results for patch refinement are denoted by squares.}
\label{fig:hconvergence1D}
\end{figure}

First, we compute convergence rates of $L^2$ errors under mesh refinement.  We begin with the acoustic wave equation in 1D, using the exact solution $p(x,t) = \cos\LRp{\frac{3\pi x}{2}}\cos\LRp{\frac{3\pi t}{2}}$ evaluated at final time $T=1/2$..  Figure~\ref{fig:hconvergence1D} shows the convergence of one-dimensional multi-patch discretizations.  Optimal $O(h^{p+1})$ rates of convergence are observed under both mesh refinement (knot insertion) and patch refinement (subdivision of patches) for first and second order formulations.  We also observe that, while optimal rates are still observed, $L^2$ errors for knot insertion are larger than $L^2$ errors for patch refinement, and that this effect grows as the degree $p$ increases.  However, we will demonstrate that, as suggested in the previous section, spline spaces achieve a lower $L^2$ error per degree of freedom.  We also observe that $L^2$ errors are larger when using smoothed knots, though (as indicated by previous experiments) this may be due to the fact that the solution is relatively smooth and non-oscillatory.  

We now compare per-dof efficiency between uniform/smoothed knot vectors and mesh/patch refinement.  Figure~\ref{fig:convergence1D} shows the convergence of $L^2$ errors for uniform and smoothed knots at various orders of approximation.  For uniform knots, mesh refinement is more efficient than patch refinement in terms of the number of degrees of freedom required to reach a certain $L^2$ error.  The same behavior is observed for smoothed knots at higher orders ($p>3$), though mesh refinement becomes less efficient than patch refinement at lower orders ($p=2$) for the first order formulation of DG.  This is consistent with observations in Section~\ref{sec:approx} that  (in one dimension) smoothed knots provide an advantage over uniform knot vectors only in the pre-asymptotic under-resolved regime.  For both the first and second order formulations, $L^2$ errors are slightly higher when using smoothed knots, though the difference between uniform and smoothed knot vectors becomes smaller as the number of patches increases.  

\begin{figure}
\centering
\subfloat[First order formulation, uniform knots]{
\begin{tikzpicture}
\begin{loglogaxis}[
width=.45\textwidth,
    xlabel={Number of degrees of freedom},   
    ylabel={$L^2$ error},
    xmin=0, xmax=150,
    ymin=1e-10, ymax=.1,    
    legend pos=north east, legend cell align=left, legend style={font=\tiny},	
       xmajorgrids=true,  ymajorgrids=true, grid style=dashed,
    legend entries={$p = 2$, $p=3$, $p=4$, $p=5$},
] 
\pgfplotsset{
cycle list={
{blue, mark=*}, {red,mark=square*}, {black,mark=triangle*},
{blue, mark=*,dashed}, {red,mark=square*,dashed}, {black,mark=triangle*,dashed},
{blue, mark=*,dotted}, {red,mark=square*,dotted}, {black,mark=triangle*,dotted},
{blue, mark=*,dashdotted}, {red,mark=square*,dashdotted}, {black,mark=triangle*,dashdotted}},
legend style={mark options={no markers}}
}
\addlegendimage{no markers,black}
\addlegendimage{no markers,black,dashed}
\addlegendimage{no markers,black,dotted}
\addlegendimage{no markers,black,dashdotted}

\addplot+[semithick, mark options={solid,fill=markercolor}]
coordinates{(12,0.0106631)(20,0.00095481)(36,0.000110851)(68,1.37364e-05)};
\addplot+[semithick, mark options={solid,fill=markercolor}]
coordinates{(24,0.000855511)(40,0.00010137)(72,1.2405e-05)};
\addplot+[semithick, mark options={solid,fill=markercolor}]
coordinates{(48,7.68373e-05)(80,1.02342e-05)};

\addplot+[semithick, mark options={solid,fill=markercolor}]
coordinates{(14,0.00183495)(22,8.58897e-05)(38,4.97075e-06)(70,3.0529e-07)};
\addplot+[semithick, mark options={solid,fill=markercolor}]
coordinates{(28,7.94928e-05)(44,4.72941e-06)(76,2.92657e-07)};
\addplot+[semithick,  mark options={solid,fill=markercolor}]
coordinates{(56,4.30145e-06)(88,2.845e-07)};

\addplot+[semithick, mark options={solid,fill=markercolor}]
coordinates{(16,0.00041217)(24,1.05595e-05)(40,2.98648e-07)(72,9.13749e-09)};
\addplot+[semithick, mark options={solid,fill=markercolor}]
coordinates{(32,8.10195e-06)(48,2.39798e-07)(80,7.44693e-09)};
\addplot+[semithick,  mark options={solid,fill=markercolor}]
coordinates{(64,2.19131e-07)(96,7.01346e-09)};

\addplot+[semithick, mark options={solid,fill=markercolor}]
coordinates{(18,7.95525e-05)(26,1.04357e-06)(42,1.41959e-08)(74,2.14126e-10)};
\addplot+[semithick, mark options={solid,fill=markercolor}]
coordinates{(36,7.59324e-07)(52,1.17281e-08)(84,1.84885e-10)};
\addplot+[semithick,  mark options={solid,fill=markercolor}]
coordinates{(72,1.08769e-08)(104,2.53553e-10)};
\end{loglogaxis}
\end{tikzpicture}
}
\subfloat[First order formulation, smoothed knots]{
\begin{tikzpicture}
\begin{loglogaxis}[
width=.45\textwidth,
    xlabel={Number of degrees of freedom},   
    ylabel={$L^2$ error},
    xmin=0, xmax=150,
    ymin=1e-10, ymax=.1,    
    legend pos=north east, legend cell align=left, legend style={font=\tiny},	
       xmajorgrids=true,  ymajorgrids=true, grid style=dashed,
    legend entries={$p = 2$, $p=3$, $p=4$, $p=5$},
] 
\pgfplotsset{
cycle list={
{blue, mark=*}, {red,mark=square*}, {black,mark=triangle*},
{blue, mark=*,dashed}, {red,mark=square*,dashed}, {black,mark=triangle*,dashed},
{blue, mark=*,dotted}, {red,mark=square*,dotted}, {black,mark=triangle*,dotted},
{blue, mark=*,dashdotted}, {red,mark=square*,dashdotted}, {black,mark=triangle*,dashdotted}},
legend style={mark options={no markers}}
}
\addlegendimage{no markers,black}
\addlegendimage{no markers,black,dashed}
\addlegendimage{no markers,black,dotted}
\addlegendimage{no markers,black,dashdotted}

\addplot+[semithick, mark options={solid,fill=markercolor}]
coordinates{(12,0.0210966)(20,0.00277327)(36,0.000377052)(68,4.89405e-05)};
\addplot+[semithick, mark options={solid,fill=markercolor}]
coordinates{(24,0.00122621)(40,0.000191822)(72,2.82418e-05)};
\addplot+[semithick, mark options={solid,fill=markercolor}]
coordinates{(48,8.23467e-05)(80,1.233e-05)};

\addplot+[semithick, mark options={solid,fill=markercolor}]
coordinates{(14,0.00350747)(22,0.000280473)(38,2.07193e-05)(70,1.3648e-06)};
\addplot+[semithick, mark options={solid,fill=markercolor}]
coordinates{(28,9.7771e-05)(44,9.29359e-06)(76,7.0952e-07)};
\addplot+[semithick,  mark options={solid,fill=markercolor}]
coordinates{(56,4.41167e-06)(88,3.45481e-07)};

\addplot+[semithick, mark options={solid,fill=markercolor}]
coordinates{(16,0.000444765)(24,2.14376e-05)(40,8.86832e-07)(72,3.23075e-08)};
\addplot+[semithick, mark options={solid,fill=markercolor}]
coordinates{(32,6.77359e-06)(48,3.72033e-07)(80,1.57506e-08)};
\addplot+[semithick,  mark options={solid,fill=markercolor}]
coordinates{(64,1.72723e-07)(96,7.66514e-09)};

\addplot+[semithick, mark options={solid,fill=markercolor}]
coordinates{(18,4.94071e-05)(26,1.49166e-06)(42,3.28811e-08)(74,5.88541e-10)};
\addplot+[semithick, mark options={solid,fill=markercolor}]
coordinates{(36,4.05142e-07)(52,1.42633e-08)(84,3.36761e-10)};
\addplot+[semithick,  mark options={solid,fill=markercolor}]
coordinates{(72,6.53937e-09)(104,2.53727e-10)};
\end{loglogaxis}
\end{tikzpicture}
}\\
\subfloat[Second order formulation, uniform knots]{
\begin{tikzpicture}
\begin{loglogaxis}[
width=.45\textwidth,
    xlabel={Number of degrees of freedom},   
    ylabel={$L^2$ error},
    xmin=0, xmax=150,
    ymin=1e-10, ymax=.1,    
    legend pos=north east, legend cell align=left, legend style={font=\tiny},	
       xmajorgrids=true,  ymajorgrids=true, grid style=dashed,
    legend entries={$p = 2$, $p=3$, $p=4$, $p=5$},
] 
\pgfplotsset{
cycle list={
{blue, mark=*}, {red,mark=square*}, {black,mark=triangle*},
{blue, mark=*,dashed}, {red,mark=square*,dashed}, {black,mark=triangle*,dashed},
{blue, mark=*,dotted}, {red,mark=square*,dotted}, {black,mark=triangle*,dotted},
{blue, mark=*,dashdotted}, {red,mark=square*,dashdotted}, {black,mark=triangle*,dashdotted}},
legend style={mark options={no markers}}
}
\addlegendimage{no markers,black}
\addlegendimage{no markers,black,dashed}
\addlegendimage{no markers,black,dotted}
\addlegendimage{no markers,black,dashdotted}

\addplot+[semithick, mark options={solid,fill=markercolor}]
coordinates{(12,0.0162082)(20,0.00141823)(36,0.000125808)(68,1.26424e-05)};
\addplot+[semithick, mark options={solid,fill=markercolor}]
coordinates{(24,0.000703955)(40,0.0001003)(72,1.20226e-05)};
\addplot+[semithick, mark options={solid,fill=markercolor}]
coordinates{(48,0.000125593)(80,1.30507e-05)};

\addplot+[semithick, mark options={solid,fill=markercolor}]
coordinates{(14,0.00271749)(22,8.35594e-05)(38,4.79948e-06)(70,2.97146e-07)};
\addplot+[semithick, mark options={solid,fill=markercolor}]
coordinates{(28,7.90339e-05)(44,4.66277e-06)(76,2.96322e-07)};
\addplot+[semithick,  mark options={solid,fill=markercolor}]
coordinates{(56,4.37416e-06)(88,2.86356e-07)};

\addplot+[semithick, mark options={solid,fill=markercolor}]
coordinates{(16,0.000623315)(24,9.82538e-06)(40,2.93292e-07)(72,9.02485e-09)};
\addplot+[semithick, mark options={solid,fill=markercolor}]
coordinates{(32,9.57229e-06)(48,2.79274e-07)(80,8.95717e-09)};
\addplot+[semithick,  mark options={solid,fill=markercolor}]
coordinates{(64,2.29205e-07)(96,8.78945e-09)};

\addplot+[semithick, mark options={solid,fill=markercolor}]
coordinates{(18,0.000116261)(26,9.33596e-07)(42,1.13206e-08)(74,1.79528e-10)};
\addplot+[semithick, mark options={solid,fill=markercolor}]
coordinates{(36,9.69517e-07)(52,1.10416e-08)(84,1.89631e-10)};
\addplot+[semithick,  mark options={solid,fill=markercolor}]
coordinates{(72,1.35402e-08)(104,1.86444e-10)};
\end{loglogaxis}
\end{tikzpicture}
}
\subfloat[Second order formulation, smoothed knots]{
\begin{tikzpicture}
\begin{loglogaxis}[
width=.45\textwidth,
    xlabel={Number of degrees of freedom},   
    ylabel={$L^2$ error},
    xmin=0, xmax=150,
    ymin=1e-10, ymax=.1,    
    legend pos=north east, legend cell align=left, legend style={font=\tiny},	
       xmajorgrids=true,  ymajorgrids=true, grid style=dashed,
    legend entries={$p = 2$, $p=3$, $p=4$, $p=5$},
] 
\pgfplotsset{
cycle list={
{blue, mark=*}, {red,mark=square*}, {black,mark=triangle*},
{blue, mark=*,dashed}, {red,mark=square*,dashed}, {black,mark=triangle*,dashed},
{blue, mark=*,dotted}, {red,mark=square*,dotted}, {black,mark=triangle*,dotted},
{blue, mark=*,dashdotted}, {red,mark=square*,dashdotted}, {black,mark=triangle*,dashdotted}},
legend style={mark options={no markers}}
}
\addlegendimage{no markers,black}
\addlegendimage{no markers,black,dashed}
\addlegendimage{no markers,black,dotted}
\addlegendimage{no markers,black,dashdotted}

\addplot+[semithick, mark options={solid,fill=markercolor}]
coordinates{(12,0.0100998)(20,0.00183953)(36,0.000334866)(68,3.38972e-05)};
\addplot+[semithick, mark options={solid,fill=markercolor}]
coordinates{(24,0.00148053)(40,0.000351277)(72,2.61946e-05)};
\addplot+[semithick, mark options={solid,fill=markercolor}]
coordinates{(48,0.000243683)(80,4.22802e-05)};

\addplot+[semithick, mark options={solid,fill=markercolor}]
coordinates{(14,0.00341013)(22,0.000116041)(38,8.20099e-06)(70,9.12238e-07)};
\addplot+[semithick, mark options={solid,fill=markercolor}]
coordinates{(28,0.000166115)(44,1.59664e-05)(76,6.36059e-07)};
\addplot+[semithick,  mark options={solid,fill=markercolor}]
coordinates{(56,1.36033e-05)(88,9.49495e-07)};

\addplot+[semithick, mark options={solid,fill=markercolor}]
coordinates{(16,0.00036462)(24,9.40051e-06)(40,1.13518e-06)(72,3.45277e-08)};
\addplot+[semithick, mark options={solid,fill=markercolor}]
coordinates{(32,1.08149e-05)(48,5.82473e-07)(80,1.9993e-08)};
\addplot+[semithick,  mark options={solid,fill=markercolor}]
coordinates{(64,5.05283e-07)(96,2.09237e-08)};

\addplot+[semithick, mark options={solid,fill=markercolor}]
coordinates{(18,4.74162e-05)(26,1.99624e-06)(42,3.33867e-08)(74,3.99503e-10)};
\addplot+[semithick, mark options={solid,fill=markercolor}]
coordinates{(36,9.26336e-07)(52,2.22652e-08)(84,3.71367e-10)};
\addplot+[semithick,  mark options={solid,fill=markercolor}]
coordinates{(72,8.35937e-09)(104,2.87484e-10)};

\end{loglogaxis}
\end{tikzpicture}
}
\caption{Convergence under mesh refinement (knot insertion) of first and second order formulations of the 1D acoustic wave equation using multi-patch DG for $p = 2,\ldots,5$.  Different markers correspond to patches refinement: two patches (circle), four patches (square), eight patches (triangle).  }
\label{fig:convergence1D}
\end{figure}

\subsection{Approximation of oscillatory functions using uniform, optimal, and smoothed knot vectors}
\label{sec:approx}

We now compare the approximation properties of splines under uniform, optimal, and smoothed knot vectors.  Recall that optimal knot vectors are defined such that the corresponding spline space minimizes the ``worst best approximation error'' in the sense of $n$-widths.  This implies that optimal knot vectors (and their approximations using smoothed knot vectors) will be less accurate than spline spaces using uniform knot vectors for certain functions.  

\begin{figure}
\centering
\subfloat[$L^2$ error vs wavelengths ($k$) per dof]{\includegraphics[width=.4675\textwidth]{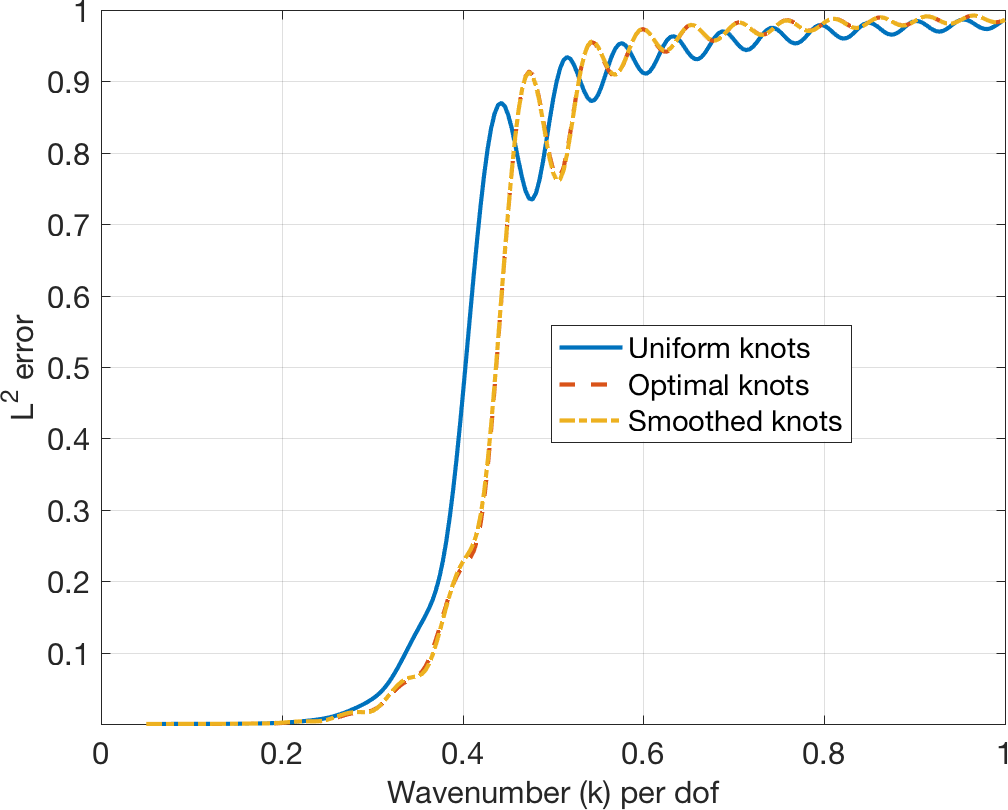}}
\hspace{.5em}
\subfloat[Semi-log scale]{\includegraphics[width=.475\textwidth]{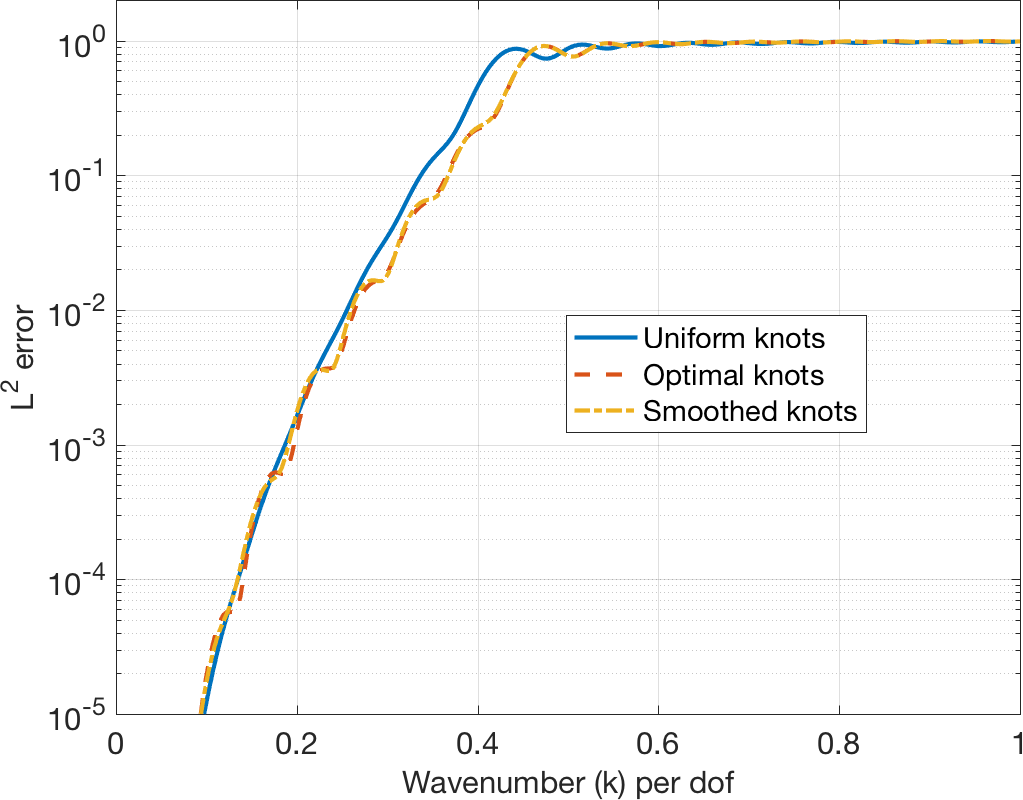}}
\caption{$L^2$ error in approximating $u(x) = \cos\LRp{\frac{(2k-1)\pi x}{2}}$ at various $k$ with $p = 4, K = 16$.}
\label{fig:ppw}
\end{figure}

Figure~\ref{fig:ppw} shows $L^2$ errors in approximating $\cos\LRp{\frac{(2k-1)\pi x}{2}}$ for various $k$ under uniform, optimal, and smoothed knots over a single patch.  \reviewerOne{Here, we approximate the function using an $L^2$ projection, using a composite Gaussian quadrature rule which is exact for the product of two degree $p$ B-splines.}  The number of wavelengths is taken to be $k = 1/{N_{\rm dofs}},\ldots, N_{\rm dofs}/2$ (where $N_{\rm dofs}$ denotes the dimension of the spline space).  We observe that, while $L^2$ errors grow with $k$ for all knot distributions, this growth of error is delayed until a slightly higher $k$ when using optimal and smoothed knot distributions (compared to uniform knot distributions).  These results also indicate that smoothed and optimal knot distributions behave similarly to uniform knots for smaller values of $k$.  

\begin{figure}
\centering
\subfloat[First order formulation]{
\begin{tikzpicture}
\begin{semilogyaxis}[
width=.45\textwidth,
    xlabel={Wavelengths $k$ per dof},   
    ylabel={$L^2$ error},
    ymin=1e-5, ymax=2,    
    legend pos=south east, legend cell align=left, legend style={font=\tiny},	
       xmajorgrids=true,  ymajorgrids=true, grid style=dashed,
] 
\addplot[color=blue,mark=*, semithick, mark options={solid,fill=markercolor}]
%coordinates{(1,1.86794e-05)(3,0.00199192)(5,0.0632619)(7,0.616786)(9,0.712371)(11,0.768705)};
coordinates{(1/11,1.86794e-05)(3/11,0.00199192)(5/11,0.0632619)(7/11,0.616786)(9/11,0.712371)(11/11,0.768705)};
\addplot[color=red,mark=*, dashed, semithick, mark options={solid,fill=markercolor}]
%coordinates{(1,5.56475e-05)(3,0.00612773)(5,0.0351819)(7,0.405041)(9,0.772862)(11,0.731632)};
coordinates{(1/11,5.56475e-05)(3/11,0.00612773)(5/11,0.0351819)(7/11,0.405041)(9/11,0.772862)(11/11,0.731632)};
\legend{Uniform knots, Smoothed knots}
\end{semilogyaxis}
\end{tikzpicture}
}
\subfloat[Second order formulation]{
\begin{tikzpicture}
\begin{semilogyaxis}[
width=.45\textwidth,
    xlabel={Wavelengths $k$ per dof},   
    ylabel={$L^2$ error},
    ymin=1e-5, ymax=2,    
    legend pos=south east, legend cell align=left, legend style={font=\tiny},	
       xmajorgrids=true,  ymajorgrids=true, grid style=dashed,
] 
\addplot[color=blue,mark=*, semithick, mark options={solid,fill=markercolor}]
%coordinates{(1,1.77384e-05)(3,0.00244169)(5,0.0592546)(7,1.07406)(9,0.662767)(11,0.872701)};
coordinates{(1/11,1.77384e-05)(3/11,0.00244169)(5/11,0.0592546)(7/11,1.07406)(9/11,0.662767)(11/11,0.872701)};
\addplot[color=red,mark=*, dashed, semithick, mark options={solid,fill=markercolor}]
%coordinates{(1,1.78285e-05)(3,0.00465007)(5,0.107561)(7,0.751711)(9,0.864029)(11,0.773654)};
coordinates{(1/11,1.78285e-05)(3/11,0.00465007)(5/11,0.107561)(7/11,0.751711)(9/11,0.864029)(11/11,0.773654)};

\legend{Uniform knots, Smoothed knots}
\end{semilogyaxis}
\end{tikzpicture}
}
\caption{Comparison of $L^2$ errors for uniform and smooth knots as a function of wavenumber $k$ using $p = 3, K= 8$ splines.}
\label{fig:wavelengthconverge}
\end{figure}
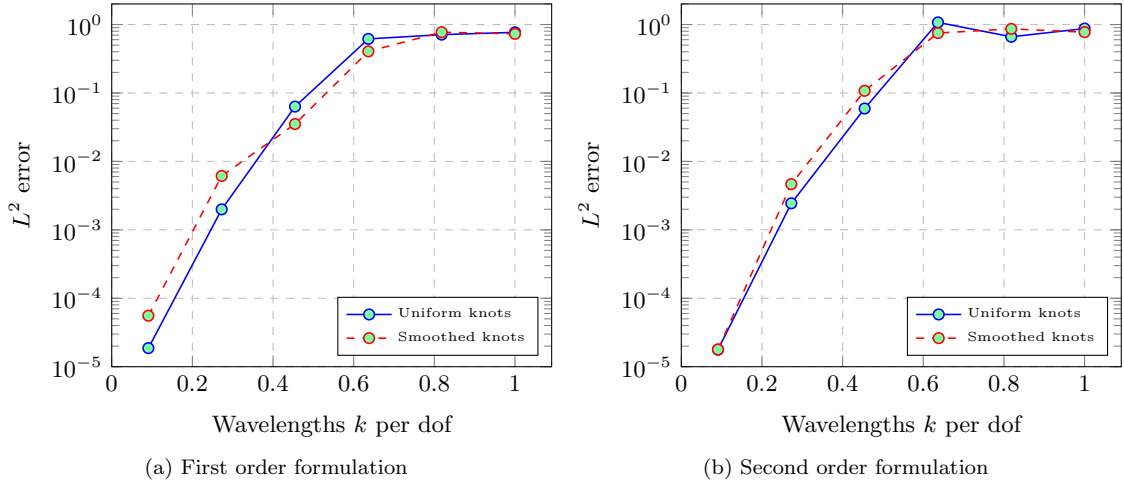

Finally, we compare DG solution errors for both uniform and smoothed knot vectors.  We use the exact pressure solution $p(x,t) = \cos\LRp{\frac{2k-1}{2}\pi x}\cos\LRp{\frac{2k-1}{2}\pi t}$  at final time $T=1/2$, and compute $L^2$ errors on a single patch domain.  Errors are computed for various wavenumbers $k$ but for a fixed order $p$ and number of elements $K$.  The $L^2$ errors (shown in Figure~\ref{fig:wavelengthconverge}) show qualitatively similarity to the approximation results of Figure~\ref{fig:ppw}, with uniform knots performing better at small $k$.  Smoothed knots result in lower errors for higher $k$, though this effect is more pronounced for first order formulations of the acoustic wave equation and less noticeable for the second order formulation.  However, in all cases, the $L^2$ errors for both uniform and smoothed knots are roughly of the same order of magnitude.  We note that the value of $k$ at which smoothed knots deliver lower error than uniform knots depends on $p$ and $K$; however, we observe in numerical experiments that this ``crossover point'' grows roughly proportionally to $K$ but depends only weakly on $p$.  

\subsection{Two-dimensional experiments}
\label{sec:twodimexp}
\begin{figure}
\centering
\subfloat[Warped mesh, $\alpha = 1/8$]{\includegraphics[width=.275\textwidth]{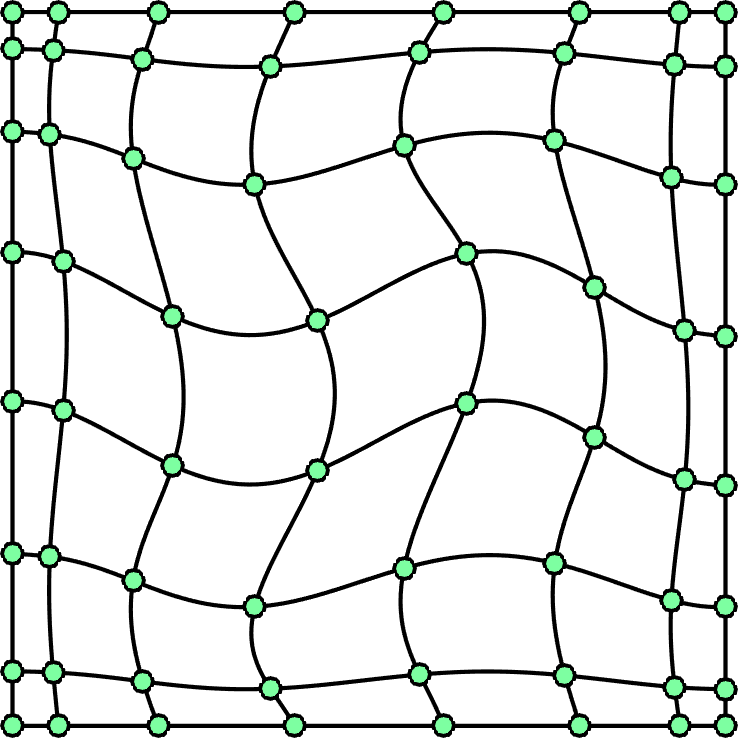}\label{subfig:mapped}}
\subfloat[$L^2$ errors for $k = 1$]{
\begin{tikzpicture}
\begin{loglogaxis}[
    legend cell align=left,
    legend style={legend pos=south east, font=\tiny},
    width=.36\textwidth,
    xlabel={Mesh size $h$},
    ylabel={$L^2$ error}, 
             ymin=1e-10, ymax=2,    
    grid style=dashed,
] 
\addplot[color=blue,mark=*,semithick, mark options={solid,fill=markercolor}]
coordinates{(0.5,0.0110908)(0.25,0.00138298)(0.125,2.29018e-05)(0.0625,3.18561e-07)};
\addplot[color=red,mark=x,semithick, mark options={solid,fill=markercolor}]
coordinates{(0.5,0.0111346)(0.25,0.00138582)(0.125,2.30061e-05)(0.0625,3.18565e-07)};
\addplot[color=black,mark=triangle*,semithick, mark options={solid,fill=markercolor}]
coordinates{(0.5,0.000986363)(0.25,8.87776e-05)(0.125,2.18816e-06)(0.0625,1.65013e-09)};

\legend{$L^2$ projection,Weight-adjusted,Difference}
%\legend{Uniform, Optimal, Smoothed}
\end{loglogaxis}
\end{tikzpicture}
}
\subfloat[$L^2$ errors for $k = 10$]{
\begin{tikzpicture}
\begin{loglogaxis}[
    legend cell align=left,
    legend style={legend pos=south east, font=\tiny},
    width=.36\textwidth,
         ymin=1e-10, ymax=2,    
    xlabel={Mesh size $h$},
    ylabel={$L^2$ error}, 
    grid style=dashed,
] 
\addplot[color=blue,mark=*,semithick, mark options={solid,fill=markercolor}]
coordinates{(0.5,0.999801)(0.25,0.903719)(0.125,0.218432)(0.0625,0.00375629)};
\addplot[color=red,mark=x,semithick, mark options={solid,fill=markercolor}]
coordinates{(0.5,0.999803)(0.25,0.90373)(0.125,0.218442)(0.0625,0.00375629)};
\addplot[color=black,mark=triangle*,semithick, mark options={solid,fill=markercolor}]
coordinates{(0.5,0.00163151)(0.25,0.00435668)(0.125,0.00205162)(0.0625,5.51212e-06)};
\legend{$L^2$ projection, Weight-adjusted, Difference}
%\legend{Uniform, Optimal, Smoothed}
\end{loglogaxis}
\end{tikzpicture}
}
\caption{Errors for $L^2$ and weight-adjusted projection when approximating $\cos\LRp{\frac{k\pi x}{2}}\cos\LRp{\frac{k\pi y}{2}}$ on a warped curvilinear mesh with $p = 4$ splines for $k = 1$ and $k=10$. }
%\caption{$L^2$ vs weight-adjusted projection. }
\label{fig:l2vswadg}
\end{figure}
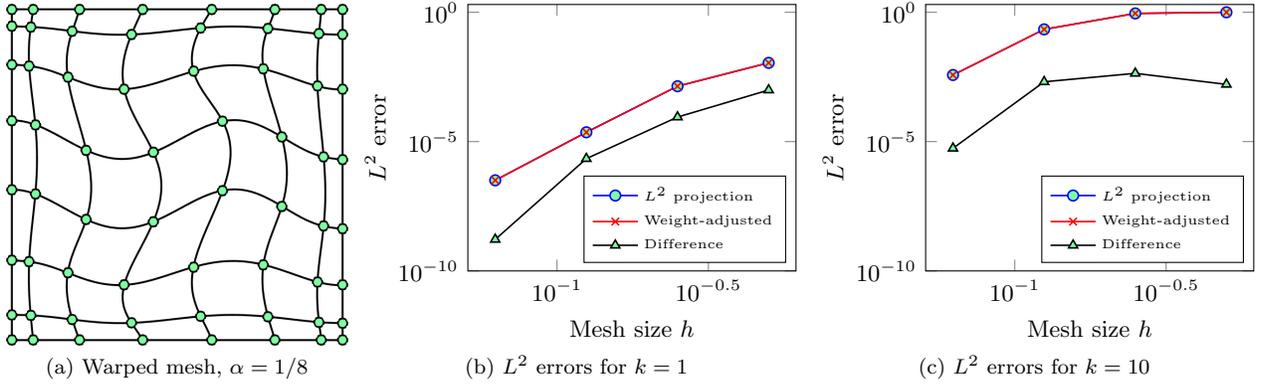

%In higher spatial dimensions, we observe similar approximation errors as reported in Figures~\ref{fig:href1D}, \ref{fig:pref1D}, and \ref{fig:ppw}.  
In higher spatial dimensions, we observe similar approximation behavior as reported in Figure~\ref{fig:ppw}.  We consider the curvilinear mapping $(\tilde{x},\tilde{y})$ of a bi-unit square $[-1,1]^2$ shown in Figure~\ref{subfig:mapped} and defined by 
\[
\tilde{x} = x + \alpha\cos\LRp{\frac{3 \pi y}{2}} \cos\LRp{\frac{\pi x}{2}}, \qquad \tilde{y} = y + \alpha \sin\LRp{\frac{3 \pi x}{2}}\cos\LRp{\frac{\pi y}{2}}.
\]
We begin by verifying that the weight-adjusted approximation to the mass matrix inverse is sufficiently accurate approximation.  Figure~\ref{fig:l2vswadg} shows $L^2$ errors in approximating the smooth function $\cos\LRp{\frac{k\pi x}{2}}\cos\LRp{\frac{k\pi y}{2}}$ from a uniform $p=4$ spline space.  The errors for the $L^2$ projection computed using the true inverse mass matrix and using the weight-adjusted mass matrix inverse are nearly identical, and are indistinguishable visually for both $k=1$ (corresponding to a smooth function) and for $k=10$ (corresponding to an oscillatory function).  

We refer to the $L^2$ projection computed using the weight adjusted mass matrix inverse as a weight-adjusted projection.  igure~\ref{fig:l2vswadg} also shows the difference between the $L^2$ and weight-adjusted projections, which we observe to be at least one order of magnitude smaller than the approximation error for both $k=1$ and $k=10$.  Additionally, this difference appears to converge more rapidly than the error under mesh refinement.  

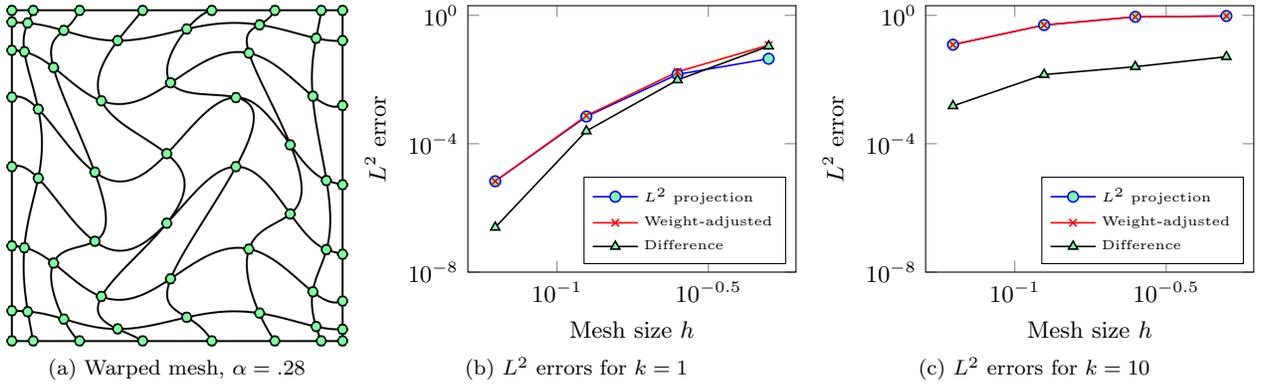
\begin{figure}
\centering
\subfloat[Warped mesh, $\alpha = .28$]{\includegraphics[width=.275\textwidth]{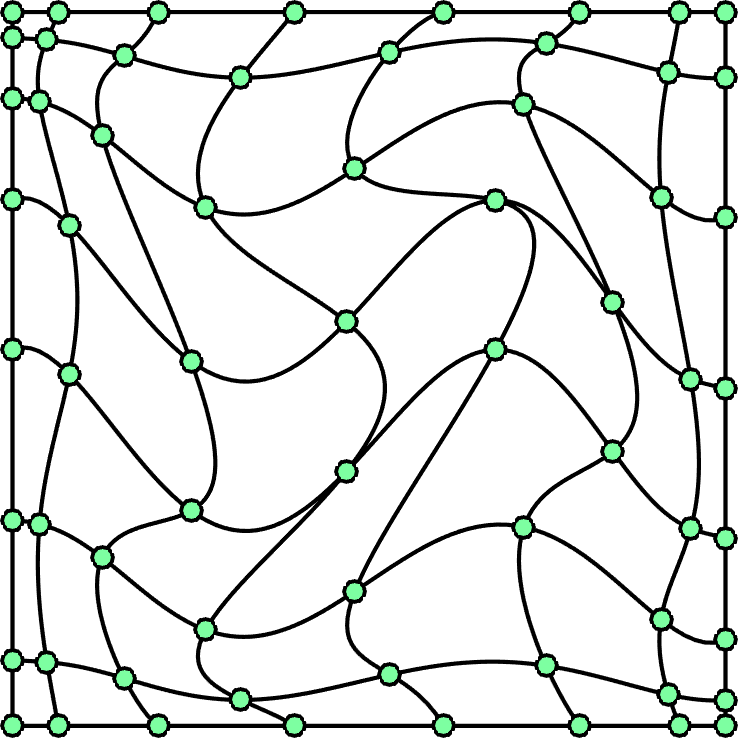}\label{subfig:mapped2}}
\subfloat[$L^2$ errors for $k = 1$]{
\begin{tikzpicture}
\begin{loglogaxis}[
    legend cell align=left,
    legend style={legend pos=south east, font=\tiny},
    width=.36\textwidth,
     ymin=1e-8, ymax=2,    
    xlabel={Mesh size $h$},
    ylabel={$L^2$ error}, 
    grid style=dashed,
] 
\addplot[color=blue,mark=*,semithick, mark options={solid,fill=markercolor}]
coordinates{(0.5,0.0444376)(0.25,0.0144725)(0.125,0.000702933)(0.0625,6.79061e-06)};
\addplot[color=red,mark=x,semithick, mark options={solid,fill=markercolor}]
coordinates{(0.5,0.1182)(0.25,0.017419)(0.125,0.000745669)(0.0625,6.79513e-06)};
\addplot[color=black,mark=triangle*,semithick, mark options={solid,fill=markercolor}]
coordinates{(0.5,0.109529)(0.25,0.00969371)(0.125,0.000248811)(0.0625,2.47881e-07)};

\legend{$L^2$ projection,Weight-adjusted,Difference}
%\legend{Uniform, Optimal, Smoothed}
\end{loglogaxis}
\end{tikzpicture}
}
\subfloat[$L^2$ errors for $k = 10$]{
\begin{tikzpicture}
\begin{loglogaxis}[
    legend cell align=left,
    legend style={legend pos=south east, font=\tiny},
    width=.36\textwidth,
         ymin=1e-8, ymax=2,    
    xlabel={Mesh size $h$},
    ylabel={$L^2$ error}, 
    grid style=dashed,
] 
\addplot[color=blue,mark=*,semithick, mark options={solid,fill=markercolor}]
coordinates{(0.5,0.952758)(0.25,0.89781)(0.125,0.4989)(0.0625,0.121744)};
\addplot[color=red,mark=x,semithick, mark options={solid,fill=markercolor}]
coordinates{(0.5,0.954139)(0.25,0.898161)(0.125,0.499104)(0.0625,0.121753)};
\addplot[color=black,mark=triangle*,semithick, mark options={solid,fill=markercolor}]
coordinates{(0.5,0.0513167)(0.25,0.0251125)(0.125,0.0142617)(0.0625,0.00150984)};

\legend{$L^2$ projection, Weight-adjusted, Difference}
%\legend{Uniform, Optimal, Smoothed}
\end{loglogaxis}
\end{tikzpicture}
}
\caption{Errors for $L^2$ and weight-adjusted projection when approximating $\cos\LRp{\frac{k\pi x}{2}}\cos\LRp{\frac{k\pi y}{2}}$ on a more aggressively warped curvilinear mesh with $p = 4$ splines for $k = 1$ and $k=10$. }
%\caption{$L^2$ vs weight-adjusted projection. }
\label{fig:l2vswadg2}
\end{figure}

\reviewerTwo{Theory in \cite{chan2016weight1} suggests that the weight-adjusted approximation becomes less accurate as the geometric mapping becomes less regular.  We test this by considering a more aggressively warped mapping with $\alpha = .28$.  This value is chosen to maximize the irregularity of $J$ without losing invertibility of the geometric mapping (for a slightly larger $\alpha = .285$, negative values of $J$ are detected, indicating that the mapping becomes non-invertible).  Figure~\ref{fig:l2vswadg2} shows a visualization of the warped mesh, along with  $L^2$ errors for both the  $L^2$ and weight-adjusted projections and the $L^2$ difference between the $L^2$ and weight-adjusted projection.   

We consider again the function $\cos\LRp{\frac{k\pi x}{2}}\cos\LRp{\frac{k\pi y}{2}}$ for $k=1$ and $k = 10$.  For $k=10$, the approximated function is oscillatory and less well-resolved, and we do not observe any significant differences between the $L^2$ and weight-adjusted projections.  However, for $k=1$, the approximated function is smooth, and we observe a noticeable difference between the $L^2$ error for the weight-adjusted and $L^2$ projections on the coarsest mesh.  This confirms the fact that the weight-adjusted projection becomes less accurate as the geometric mapping approaches non-invertibility.  However, even for this more irregular mapping, the difference between the weight-adjusted and $L^2$ projection decreases rapidly under mesh refinement.  

Finally, we note that it was shown in \cite{chan2016weight2} that computing the weight-adjusted projection of some function $u$ is equivalent to computing 
\[
\widehat{\Pi}_h\LRp{\frac{1}{J}\widehat{\Pi}_h \LRp{uJ}},
\]
where $\widehat{\Pi}_h$ denotes the $L^2$ projection operator onto $V_h\LRp{\widehat{D}}$.  This suggests that the weight-adjusted inverse mass matrix is a highly accurate approximation to the inverse curvilinear mass matrix, so long as the geometric mapping is sufficiently well-resolved approximated by the spline approximation space on the reference element.
%These experiments suggest that under mesh refinement, 
}

\subsubsection{Approximation errors under degree refinement}

\begin{figure}
\centering
%\subfloat[Warped mesh, $\alpha = 1/8$]{\includegraphics[width=.5\textwidth]{figs/mapped.png}}
\subfloat[$L^2$ errors ]{
\begin{tikzpicture}
\begin{semilogyaxis}[
    legend cell align=left,
    legend style={legend pos=north east, font=\tiny},
    width=.5\textwidth,
    xlabel={Degrees of freedom},
    ylabel={$L^2$ error}, 
    grid style=dashed,
] 

\addplot[color=blue,mark=*,semithick, mark options={solid,fill=markercolor}]
coordinates{(16,0.0345598)(36,0.00291687)(64,0.000910808)(100,0.000228283)(144,4.68154e-05)(196,1.08049e-05)(256,2.35056e-06)}
[yshift=4pt] node[above, pos=.8, color=black] {$\alpha = 1/64$};
\addplot[color=red,mark=square*,semithick, mark options={solid,fill=markercolor}]
coordinates{(16,0.0345598)(36,0.0018151)(64,0.000489718)(100,8.29769e-05)(144,1.29416e-05)};
\addplot[color=black,mark=triangle*,semithick, mark options={solid,fill=markercolor}]
coordinates{(16,0.0345598)(36,0.00184734)(64,0.000495297)(100,8.68887e-05)(144,1.63983e-05)(196,1.69855e-06)(256,4.46827e-07)};
\addplot[color=violet,mark=diamond*,semithick, mark options={solid,fill=markercolor}]
coordinates{(9,0.0353096)(16,0.0346236)(25,0.00282808)(36,0.00233183)(49,0.00142429)(64,0.000710684)(81,0.000305903)(100,0.000123394)(121,4.1829e-05)(144,2.40683e-05)(169,6.2327e-06)(196,5.06423e-06)(225,1.05795e-06)(256,8.08652e-07)};

\addplot[color=blue,mark=*,dashed,semithick, mark options={solid,fill=markercolor}]
coordinates{(16,0.0367809)(36,0.0169805)(64,0.0111257)(100,0.00520972)(144,0.00196454)(196,0.000609694)(256,0.000194924)}
[yshift=4pt] node[above, pos=.85,color=black] {$\alpha = 1/8$};
\addplot[color=red,mark=square*,dashed,semithick, mark options={solid,fill=markercolor}]
coordinates{(16,0.0367809)(36,0.014593)(64,0.00828577)(100,0.00316138)(144,0.000836099)};
\addplot[color=black,mark=triangle*,dashed,semithick, mark options={solid,fill=markercolor}]
coordinates{(16,0.0367809)(36,0.014438)(64,0.008333)(100,0.00330281)(144,0.000899406)(196,0.000319152)(256,0.000106414)};
\addplot[color=violet,mark=diamond*,dashed,semithick, mark options={solid,fill=markercolor}]
coordinates{(9,0.0744743)(16,0.0402554)(25,0.0332601)(36,0.0185583)(49,0.01937)(64,0.00949402)(81,0.00704509)(100,0.00431925)(121,0.00311497)(144,0.00145279)(169,0.00121632)(196,0.000388163)(225,0.00042613)};

\legend{Uniform, Optimal, Smoothed,Polynomial}
%\legend{Uniform, Optimal, Smoothed}
\end{semilogyaxis}
\end{tikzpicture}
}
\caption{$L^2$ best approximation errors for $\cos\LRp{\frac{k\pi x}{2}}\cos\LRp{\frac{k\pi y}{2}}, k = 1$ under different mesh warpings and degree refinement $p = 2,\ldots, 8$.  All results assume that spline spaces contain $K = p$ elements, \reviewerTwo{except for the ``Polynomial'' curve, which assumes degree $p = 2,\ldots, 12$ polynomial approximations over a single element}. }
\label{fig:mapped}
\end{figure}
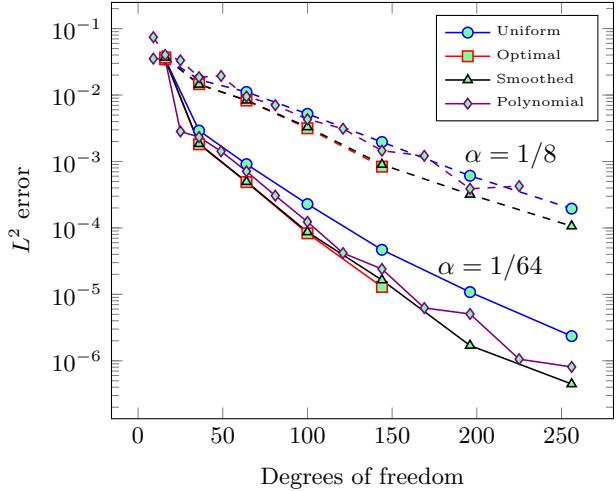

We now investigate approximation errors under degree refinement on curved meshes.  These experiments suggest that splines under optimal and smoothed knot vectors possess some advantages over uniform knot distributions in the presence of non-affine mappings.  Figure~\ref{fig:mapped} shows computed $L^2$ errors in approximating $\cos\LRp{\frac{3\pi x}{2}}\cos\LRp{\frac{3\pi y}{2}}$ on a curvilinear mesh.  We examine errors under degree refinement for both tensor-product polynomial and spline spaces (where we fix the number of  elements $K=p$).  \reviewerOne{We compute approximations using quadrature-based weight-adjusted projection, where the inverse curvilinear mass matrix is approximated using a weight-adjusted approximation.}  

In contrast to the affine case, we observe that both polynomials and spline spaces with uniform knot vectors produce $L^2$ errors of similar magnitude per degree of freedom, even for low values of $\alpha$ (lightly warped curvilinear meshes).  We also observe that splines spaces using optimal and smoothed knot vectors produce lower $L^2$ errors than both polynomials and uniform splines.  These numerical results suggest that $L^2$ approximation errors under spline spaces are comparable to polynomials not only for oscillatory functions (as was observed in 1D), but also for curvilinear mappings.  Additionally, spline spaces with optimal or smoothed knot vectors produce better approximations for both oscillatory functions and curved mappings.  

\subsubsection{$L^2$ errors for multi-patch DG solutions }

Next, we examine the behavior of multi-patch DG methods for the acoustic wave equation in two dimensions using both uniform and smoothed knot vectors.  We compute $L^2$ errors for the exact solution $p(x,y,t) = \cos\LRp{\frac{3\pi x}{2}}\cos\LRp{\frac{3\pi y}{2}}\cos\LRp{\frac{3\sqrt{2} \pi t}{2}}$ at $T = 1/2$.  For affine patches, the convergence of $L^2$ errors is similar to the 1D case; thus, we consider a curved interior warping of a square (Figure~\ref{fig:mapped}) in order to account for the effect of geometric mappings.  

\begin{figure}
\subfloat[First order formulation, $\alpha = 1/8$]{
\begin{tikzpicture}
\begin{loglogaxis}[
    legend cell align=left,
    legend style={legend pos=south east, font=\tiny},
    width=.4\textwidth,
     ymin=1e-17, ymax=2,    
    xlabel={Mesh size $h$},
    ylabel={$L^2$ error}, 
    grid style=dashed,
] 
\addplot[color=blue,mark=*,semithick, mark options={solid,fill=markercolor}]
coordinates{(0.5,0.154255)(0.25,0.0280139)(0.125,0.000775805)(0.0625,7.51259e-06)};
\addplot[color=red,mark=x,semithick, mark options={solid,fill=markercolor}]
coordinates{(0.5,0.154523)(0.25,0.0280342)(0.125,0.000775107)(0.0625,7.51217e-06)};
\addplot[color=black,mark=triangle*,semithick, mark options={solid,fill=markercolor}]
coordinates{(0.5,7.56052e-06)(0.25,1.32834e-07)(0.125,1.51285e-11)(0.0625,1.69165e-17)};

\legend{$L^2$ projection,Weight-adjusted,Difference}
%\legend{Uniform, Optimal, Smoothed}
\end{loglogaxis}
\end{tikzpicture}
}
\subfloat[Second order formulation, $\alpha = 1/8$]{
\begin{tikzpicture}
\begin{loglogaxis}[
    legend cell align=left,
    legend style={legend pos=south east, font=\tiny},
    width=.4\textwidth,
     ymin=1e-17, ymax=2,    
    xlabel={Mesh size $h$},
    ylabel={$L^2$ error}, 
    grid style=dashed,
] 
\addplot[color=blue,mark=*,semithick, mark options={solid,fill=markercolor}]
coordinates{(0.5,0.270044)(0.25,0.0483258)(0.125,0.00134115)(0.0625,1.00287e-05)};
\addplot[color=red,mark=x,semithick, mark options={solid,fill=markercolor}]
coordinates{(0.5,0.272282)(0.25,0.0483238)(0.125,0.00133867)(0.0625,1.0026e-05)};
\addplot[color=black,mark=triangle*,semithick, mark options={solid,fill=markercolor}]
coordinates{(0.5,4.77833e-05)(0.25,1.1668e-06)(0.125,3.6017e-11)(0.0625,2.65088e-16)};

\legend{$L^2$ projection,Weight-adjusted,Difference}
%\legend{Uniform, Optimal, Smoothed}
\end{loglogaxis}
\end{tikzpicture}
}\\
\subfloat[First order formulation, $\alpha = .28$]{
\begin{tikzpicture}
\begin{loglogaxis}[
    legend cell align=left,
    legend style={legend pos=south east, font=\tiny},
    width=.4\textwidth,
     ymin=1e-17, ymax=2,    
    xlabel={Mesh size $h$},
    ylabel={$L^2$ error}, 
    grid style=dashed,
] 
\addplot[color=blue,mark=*,semithick, mark options={solid,fill=markercolor}]
coordinates{(0.5,0.338603)(0.25,0.143169)(0.125,0.0234825)(0.0625,0.00033695)};
\addplot[color=red,mark=x,semithick, mark options={solid,fill=markercolor}]
coordinates{(0.5,0.388935)(0.25,0.146233)(0.125,0.0235623)(0.0625,0.000336824)};
\addplot[color=black,mark=triangle*,semithick, mark options={solid,fill=markercolor}]
coordinates{(0.5,0.00977253)(0.25,0.000187168)(0.125,3.16763e-07)(0.0625,9.44863e-12)};

\legend{$L^2$ projection,Weight-adjusted,Difference}
%\legend{Uniform, Optimal, Smoothed}
\end{loglogaxis}
\end{tikzpicture}
}
\subfloat[Second order formulation, $\alpha = .28$]{
\begin{tikzpicture}
\begin{loglogaxis}[
    legend cell align=left,
    legend style={legend pos=south east, font=\tiny},
    width=.4\textwidth,
     ymin=1e-17, ymax=2,  
           xlabel={Mesh size $h$},
    ylabel={$L^2$ error}, 
    grid style=dashed,
] 
\addplot[color=blue,mark=*,semithick, mark options={solid,fill=markercolor}]
coordinates{(0.5,0.660221)(0.25,0.272201)(0.125,0.0336296)(0.0625,0.000337592)};
\addplot[color=red,mark=x,semithick, mark options={solid,fill=markercolor}]
coordinates{(0.5,0.68491)(0.25,0.280896)(0.125,0.0335946)(0.0625,0.000337711)};
\addplot[color=black,mark=triangle*,semithick, mark options={solid,fill=markercolor}]
coordinates{(0.5,0.0580108)(0.25,0.000653681)(0.125,7.91095e-07)(0.0625,1.2372e-11)};

\legend{$L^2$ projection,Weight-adjusted,Difference}
%\legend{Uniform, Optimal, Smoothed}
\end{loglogaxis}
\end{tikzpicture}
}

\caption{$L^2$ errors for the acoustic wave equation when inverting the curvilinear mass matrix and the weight-adjusted mass matrix.   Uniform $p = 4$ splines are used, and the meshes are the two warped curvilinear meshes shown in Figures~\ref{fig:l2vswadg} and \ref{fig:l2vswadg2}. }
\label{fig:l2vswadgerror}
\end{figure}
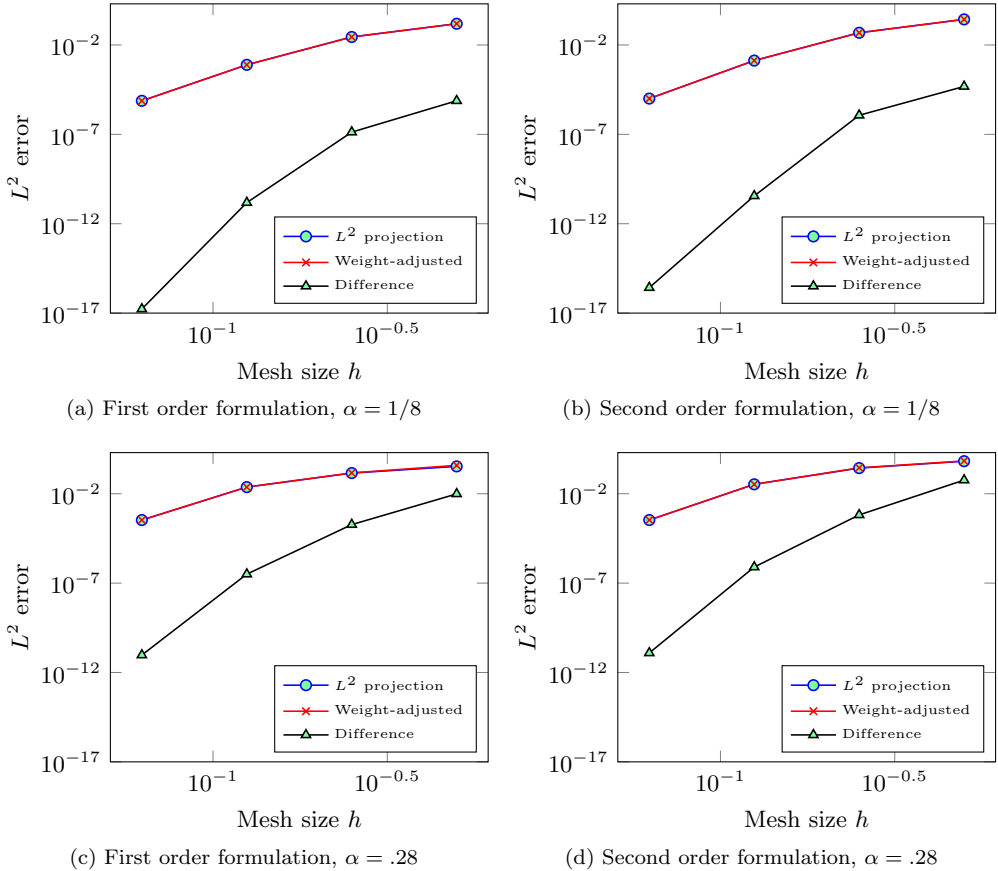

\begin{figure}
\centering
\subfloat[First order formulation, uniform knots]{
\begin{tikzpicture}
\begin{loglogaxis}[
width=.45\textwidth,
    xlabel={Number of degrees of freedom},   
    ylabel={$L^2$ error},
    xmin=0, xmax=15000,
    ymin=1e-8, ymax=1,    
    legend pos=north east, legend cell align=left, legend style={font=\tiny},	
       xmajorgrids=true,  ymajorgrids=true, grid style=dashed,
    legend entries={$p = 2$, $p=5$},
] 
\pgfplotsset{
cycle list={
{blue, mark=*}, {red,mark=square*}, {black,mark=triangle*},
{blue, mark=*,dashed}, {red,mark=square*,dashed}, {black,mark=triangle*,dashed},
{blue, mark=*,dotted}, {red,mark=square*,dotted}, {black,mark=triangle*,dotted},
{blue, mark=*,dashdotted}, {red,mark=square*,dashdotted}, {black,mark=triangle*,dashdotted}},
legend style={mark options={no markers}}
}
\addlegendimage{no markers,black}
\addlegendimage{no markers,black,dashed}
\addlegendimage{no markers,black,dotted}
\addlegendimage{no markers,black,dashdotted}

\addplot+[semithick, mark options={solid,fill=markercolor}]
coordinates{(64,0.242695)(144,0.0715275)(400,0.00658214)(1296,0.000525301)(4624,6.03423e-05)};
\addplot+[semithick, mark options={solid,fill=markercolor}]
coordinates{(256,0.0357164)(576,0.00686446)(1600,0.000545882)(5184,6.20872e-05)};
\addplot+[semithick, mark options={solid,fill=markercolor}]
coordinates{(1024,0.00435721)(2304,0.000488199)(6400,5.83311e-05)};

%\addplot+[semithick, mark options={solid,fill=markercolor}]
%coordinates{(100,0.132165)(196,0.0334468)(484,0.00222776)(1444,6.0903e-05)(4900,3.07891e-06)};
%\addplot+[semithick, mark options={solid,fill=markercolor}]
%coordinates{(400,0.00927583)(784,0.00166601)(1936,5.76982e-05)(5776,3.01869e-06)};
%\addplot+[semithick,  mark options={solid,fill=markercolor}]
%coordinates{(1600,0.000542582)(3136,5.61886e-05)(7744,2.98486e-06)};
%
%\addplot+[semithick, mark options={solid,fill=markercolor}]
%coordinates{(144,0.0619962)(256,0.017935)(576,0.000867235)(1600,9.30881e-06)(5184,2.06241e-07)};
%\addplot+[semithick, mark options={solid,fill=markercolor}]
%coordinates{(576,0.00259989)(1024,0.000491411)(2304,9.22672e-06)(6400,2.10835e-07)};
%\addplot+[semithick,  mark options={solid,fill=markercolor}]
%coordinates{(2304,7.22065e-05)(4096,8.75141e-06)(9216,1.9849e-07)};

\addplot+[semithick, mark options={solid,fill=markercolor}]
coordinates{(196,0.0288428)(324,0.007685)(676,0.000482686)(1764,1.86168e-06)(5476,1.65548e-08)};
\addplot+[semithick, mark options={solid,fill=markercolor}]
coordinates{(784,0.000691368)(1296,0.000135871)(2704,1.79719e-06)(7056,1.65087e-08)};
\addplot+[semithick,  mark options={solid,fill=markercolor}]
coordinates{(3136,9.71515e-06)(5184,1.35491e-06)(10816,1.51573e-08)};

\end{loglogaxis}
\end{tikzpicture}
}
\subfloat[First order formulation, smoothed knots]{
\begin{tikzpicture}
\begin{loglogaxis}[
width=.45\textwidth,
    xlabel={Number of degrees of freedom},   
    ylabel={$L^2$ error},
    xmin=0, xmax=15000,
    ymin=1e-8, ymax=1,    
    legend pos=north east, legend cell align=left, legend style={font=\tiny},	
       xmajorgrids=true,  ymajorgrids=true, grid style=dashed,
        legend entries={$p = 2$, $p=5$},
] 
\pgfplotsset{
cycle list={
{blue, mark=*}, {red,mark=square*}, {black,mark=triangle*},
{blue, mark=*,dashed}, {red,mark=square*,dashed}, {black,mark=triangle*,dashed},
{blue, mark=*,dotted}, {red,mark=square*,dotted}, {black,mark=triangle*,dotted},
{blue, mark=*,dashdotted}, {red,mark=square*,dashdotted}, {black,mark=triangle*,dashdotted}},
legend style={mark options={no markers}}
}
\addlegendimage{no markers,black}
\addlegendimage{no markers,black,dashed}
\addlegendimage{no markers,black,dotted}
\addlegendimage{no markers,black,dashdotted}

\addplot+[semithick, mark options={solid,fill=markercolor}]
coordinates{(64,0.153639)(144,0.0418493)(400,0.00442582)(1296,0.000564782)(4624,7.79315e-05)};
\addplot+[semithick, mark options={solid,fill=markercolor}]
coordinates{(256,0.0256049)(576,0.00420403)(1600,0.000586667)(5184,8.83518e-05)};
\addplot+[semithick, mark options={solid,fill=markercolor}]
coordinates{(1024,0.00388923)(2304,0.000639539)(6400,9.87198e-05)};

%\addplot+[semithick, mark options={solid,fill=markercolor}]
%coordinates{(100,0.0930522)(196,0.0171325)(484,0.00138355)(1444,8.64099e-05)(4900,4.83299e-06)};
%\addplot+[semithick, mark options={solid,fill=markercolor}]
%coordinates{(400,0.0072846)(784,0.000906647)(1936,8.32599e-05)(5776,5.96369e-06)};
%\addplot+[semithick,  mark options={solid,fill=markercolor}]
%coordinates{(1600,0.000530081)(3136,6.26226e-05)(7744,4.79183e-06)};
%
%\addplot+[semithick, mark options={solid,fill=markercolor}]
%coordinates{(144,0.0438962)(256,0.00927686)(576,0.000397798)(1600,1.40392e-05)(5184,3.60713e-07)};
%\addplot+[semithick, mark options={solid,fill=markercolor}]
%coordinates{(576,0.002197)(1024,0.000194099)(2304,1.36422e-05)(6400,4.07817e-07)};
%\addplot+[semithick,  mark options={solid,fill=markercolor}]
%coordinates{(2304,6.84186e-05)(4096,5.84491e-06)(9216,2.17164e-07)};

\addplot+[semithick, mark options={solid,fill=markercolor}]
coordinates{(196,0.0209849)(324,0.00383302)(676,0.00020681)(1764,2.96835e-06)(5476,3.21437e-08)};
\addplot+[semithick, mark options={solid,fill=markercolor}]
coordinates{(784,0.000588702)(1296,5.06718e-05)(2704,1.96211e-06)(7056,2.59168e-08)};
\addplot+[semithick,  mark options={solid,fill=markercolor}]
coordinates{(3136,8.87494e-06)(5184,4.8802e-07)(10816,1.50029e-08)};

\end{loglogaxis}
\end{tikzpicture}
}\\
\subfloat[Second order formulation, uniform knots]{
\begin{tikzpicture}
\begin{loglogaxis}[
width=.45\textwidth,
    xlabel={Number of degrees of freedom},   
    ylabel={$L^2$ error},
    xmin=0, xmax=15000,
    ymin=1e-8, ymax=1,    
    legend pos=north east, legend cell align=left, legend style={font=\tiny},	
       xmajorgrids=true,  ymajorgrids=true, grid style=dashed,
    legend entries={$p = 2$, $p=5$},
] 
\pgfplotsset{
cycle list={
{blue, mark=*}, {red,mark=square*}, {black,mark=triangle*},
{blue, mark=*,dashed}, {red,mark=square*,dashed}, {black,mark=triangle*,dashed},
{blue, mark=*,dotted}, {red,mark=square*,dotted}, {black,mark=triangle*,dotted},
{blue, mark=*,dashdotted}, {red,mark=square*,dashdotted}, {black,mark=triangle*,dashdotted}},
legend style={mark options={no markers}}
}
\addlegendimage{no markers,black}
\addlegendimage{no markers,black,dashed}
\addlegendimage{no markers,black,dotted}
\addlegendimage{no markers,black,dashdotted}

\addplot+[semithick, mark options={solid,fill=markercolor}]
coordinates{(64,0.359903)(144,0.0988039)(400,0.00700231)(1296,0.000596158)(4624,6.44021e-05)};
\addplot+[semithick, mark options={solid,fill=markercolor}]
coordinates{(256,0.0647089)(576,0.00662372)(1600,0.000637467)(5184,6.86903e-05)};
\addplot+[semithick, mark options={solid,fill=markercolor}]
coordinates{(1024,0.00635845)(2304,0.000602189)(6400,6.75732e-05)};

%\addplot+[semithick, mark options={solid,fill=markercolor}]
%coordinates{(100,0.190103)(196,0.0415714)(484,0.00204923)(1444,6.0709e-05)(4900,3.1327e-06)};
%\addplot+[semithick, mark options={solid,fill=markercolor}]
%coordinates{(400,0.0181315)(784,0.0017007)(1936,5.89272e-05)(5776,3.07284e-06)};
%\addplot+[semithick,  mark options={solid,fill=markercolor}]
%coordinates{(1600,0.00077351)(3136,5.6308e-05)(7744,3.01529e-06)};
%
%\addplot+[semithick, mark options={solid,fill=markercolor}]
%coordinates{(144,0.113336)(256,0.0242542)(576,0.000803146)(1600,1.00835e-05)(5184,2.17046e-07)};
%\addplot+[semithick, mark options={solid,fill=markercolor}]
%coordinates{(576,0.00480314)(1024,0.000546509)(2304,9.80952e-06)(6400,2.23621e-07)};
%\addplot+[semithick,  mark options={solid,fill=markercolor}]
%coordinates{(2304,0.000104349)(4096,9.11874e-06)(9216,2.13889e-07)};

\addplot+[semithick, mark options={solid,fill=markercolor}]
coordinates{(196,0.0419686)(324,0.0108523)(676,0.00046736)(1764,1.88401e-06)(5476,1.74519e-08)};
\addplot+[semithick, mark options={solid,fill=markercolor}]
coordinates{(784,0.00112145)(1296,0.000171121)(2704,1.86063e-06)(7056,1.87817e-08)};
\addplot+[semithick,  mark options={solid,fill=markercolor}]
coordinates{(3136,1.45546e-05)(5184,1.52434e-06)(10816,1.79508e-08)};
\end{loglogaxis}
\end{tikzpicture}
}
\subfloat[Second order formulation, smoothed knots]{
\begin{tikzpicture}
\begin{loglogaxis}[
width=.45\textwidth,
    xlabel={Number of degrees of freedom},   
    ylabel={$L^2$ error},
    xmin=0, xmax=15000,
    ymin=1e-8, ymax=1,    
    legend pos=north east, legend cell align=left, legend style={font=\tiny},	
       xmajorgrids=true,  ymajorgrids=true, grid style=dashed,
        legend entries={$p = 2$, $p=5$},
] 
\pgfplotsset{
cycle list={
{blue, mark=*}, {red,mark=square*}, {black,mark=triangle*},
{blue, mark=*,dashed}, {red,mark=square*,dashed}, {black,mark=triangle*,dashed},
{blue, mark=*,dotted}, {red,mark=square*,dotted}, {black,mark=triangle*,dotted},
{blue, mark=*,dashdotted}, {red,mark=square*,dashdotted}, {black,mark=triangle*,dashdotted}},
legend style={mark options={no markers}}
}
\addlegendimage{no markers,black}
\addlegendimage{no markers,black,dashed}
\addlegendimage{no markers,black,dotted}
\addlegendimage{no markers,black,dashdotted}

\addplot+[semithick, mark options={solid,fill=markercolor}]
coordinates{(64,0.354708)(144,0.0897704)(400,0.00913605)(1296,0.0010847)(4624,0.000122353)};
\addplot+[semithick, mark options={solid,fill=markercolor}]
coordinates{(256,0.0635767)(576,0.00668419)(1600,0.00106006)(5184,0.000137968)};
\addplot+[semithick, mark options={solid,fill=markercolor}]
coordinates{(1024,0.00603353)(2304,0.000890171)(6400,0.000122191)};

%\addplot+[semithick, mark options={solid,fill=markercolor}]
%coordinates{(100,0.189564)(196,0.0341518)(484,0.00200879)(1444,0.000127476)(4900,9.41605e-06)};
%\addplot+[semithick, mark options={solid,fill=markercolor}]
%coordinates{(400,0.0180879)(784,0.0014787)(1936,0.000118006)(5776,1.07647e-05)};
%\addplot+[semithick,  mark options={solid,fill=markercolor}]
%coordinates{(1600,0.00075651)(3136,8.13628e-05)(7744,8.37923e-06)};
%
%\addplot+[semithick, mark options={solid,fill=markercolor}]
%coordinates{(144,0.113363)(256,0.0163238)(576,0.00056338)(1600,2.35036e-05)(5184,5.80131e-07)};
%\addplot+[semithick, mark options={solid,fill=markercolor}]
%coordinates{(576,0.0047926)(1024,0.000296857)(2304,1.63539e-05)(6400,5.98945e-07)};
%\addplot+[semithick,  mark options={solid,fill=markercolor}]
%coordinates{(2304,0.000103714)(4096,8.61181e-06)(9216,4.79104e-07)};

\addplot+[semithick, mark options={solid,fill=markercolor}]
coordinates{(196,0.0420621)(324,0.00671924)(676,0.000284211)(1764,3.67017e-06)(5476,4.51605e-08)};
\addplot+[semithick, mark options={solid,fill=markercolor}]
coordinates{(784,0.00111957)(1296,6.0986e-05)(2704,2.161e-06)(7056,5.66052e-08)};
\addplot+[semithick,  mark options={solid,fill=markercolor}]
coordinates{(3136,1.45279e-05)(5184,8.23008e-07)(10816,9.14345e-08)};

\end{loglogaxis}
\end{tikzpicture}
}
\caption{Convergence under mesh refinement (knot insertion) of first and second order formulations of the 2D acoustic wave equation using multi-patch DG for $p = 2$ and $p=5$ on \reviewerOne{a warped curved square domain with $\alpha = 1/8$}.  Different markers correspond to patches refinement: two patches (circle), four patches (square), eight patches (triangle).  }
\label{fig:convergence2D}
\end{figure}
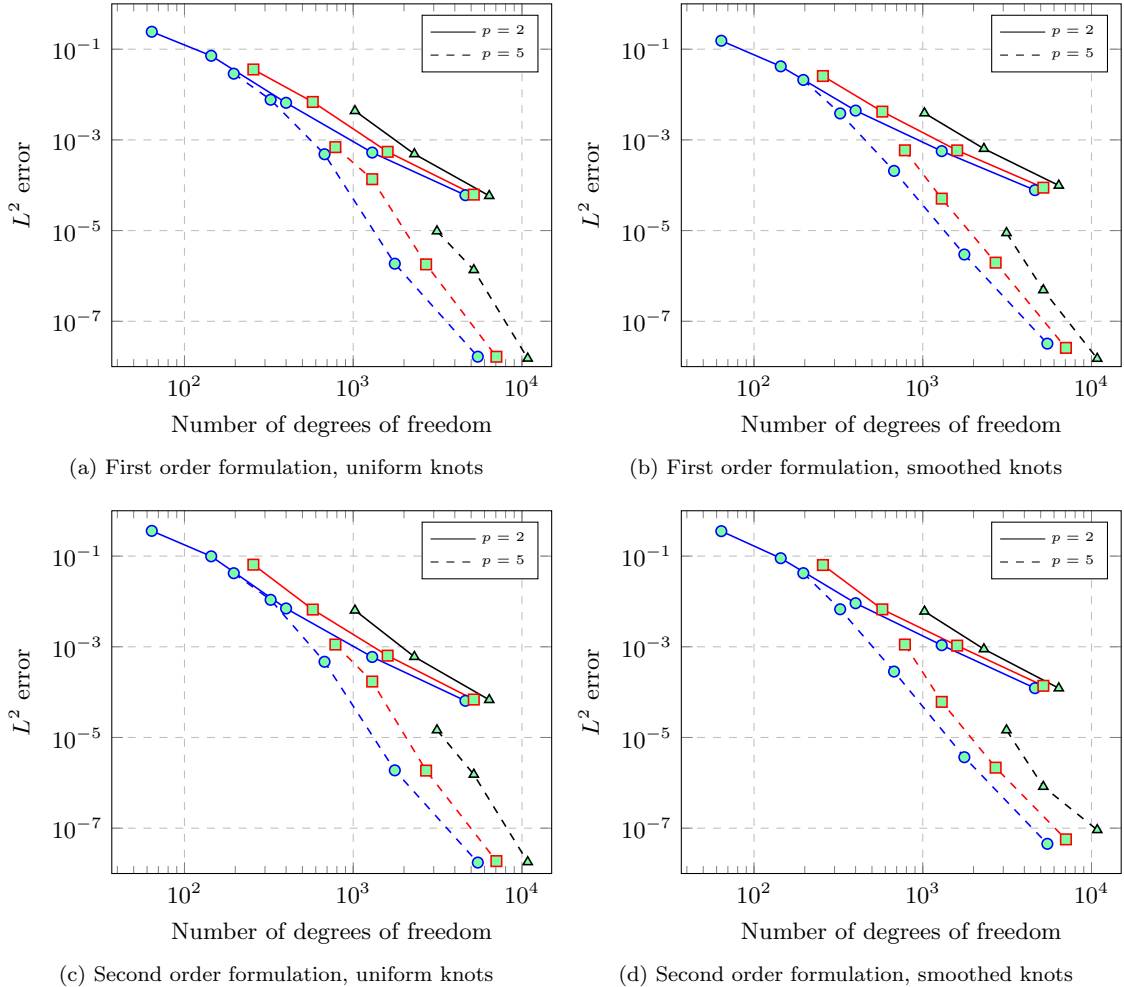

\reviewerTwo{We first compare $L^2$ pressure errors when using the curvilinear mass matrix inverse and the weight adjusted mass matrix inverse.  In both cases, uniform splines of degree $p=4$ are used.  Figure~\ref{fig:l2vswadgerror} shows $L^2$ errors for both first and second order formulations for both a lightly warped curvilinear mesh ($\alpha = 1/8$) and a heavily warped curvilinear mesh ($\alpha = .28$).   We observe that the $L^2$ errors are visually indistinguishable, and that the $L^2$ difference between the two solutions is several orders of magnitude smaller than the $L^2$ discretization error.  This mirrors the qualitative behavior observed in Figures~\ref{fig:l2vswadg} and \ref{fig:l2vswadg2}; however, the $L^2$ difference between the solutions obtained using the curvilinear and weight-adjusted mass matrix inverses is significantly smaller than the difference between $L^2$ and weight-adjusted projections. }  $L^2$ errors are shown in Figure~\ref{fig:convergence2D}, where, for clarity, we show results only for $p=2$ and $p=5$.

For first order DG formulations, we observe that $L^2$ errors for uniform and smoothed knot vectors are closer \reviewerOne{than in} the one-dimensional case, which appears to be due to \reviewerOne{the presence of the non-affine (curvilinear) geometric mapping}.  \reviewerOne{This suggests that the use of smoothed knots instead of uniform knots can improve solver efficiency, since it is possible to take a larger maximum stable timestep without reducing accuracy.  A more rigorous study of this claim is still required.}  Under the second order DG formulation, $L^2$ errors for smoothed knot vectors behave less consistently.  As with the first order formulation, the use of smoothed knot vectors reduces errors in the pre-asymptotic range for the second order formulation with $p > 2$.  However, in contrast with the first order formulation, $L^2$ errors for smoothed knots become larger than errors for uniform knot vectors under continued patch or mesh refinement.\footnote{Smoothed knot vectors tend to offer less of an advantage over uniform knot vectors for second order formulations.  This may be due to the fact that smoothed knot vectors are constructed to optimize the $L^2$ sup-inf in which first order DG formulations are more naturally convergent  \cite{warburton2013low, chan2016weight1}, while second order DG formulations converge in a stronger energy norm.  However, the precise reasons for this difference are not immediately clear, and will be investigated in future work.}

Finally, examine the effect of mesh warping on multi-patch DG solutions using both uniform and smoothed knot vectors.  We use the mapping shown in Figure~\ref{fig:mapped}, and compute $L^2$ errors as the warping factor $\alpha$ increases from $0$ to $1/4$.  We fix $K = 8$ and compute errors for $p = 3, 4, 5$ on a single patch domain.  For both first and second order formulations of the acoustic wave equation, smooth knots deliver slightly lower $L^2$ errors than uniform knots.  As with experiments in Figure~\ref{fig:wavelengthconverge}, increasing $K$ increases the value of $\alpha$ (the ``crossover point'') at which smoothed knots become more accurate than uniform knots.  

\begin{figure}
\centering
\subfloat[First order formulation]{
\begin{tikzpicture}
\begin{semilogyaxis}[
width=.45\textwidth,
    xlabel={Warping factor $\alpha$},   
    ylabel={$L^2$ error},
    ymin=1e-3, ymax=.1,    
    legend pos=south east, legend cell align=left, legend style={font=\tiny},	
       xmajorgrids=true,  ymajorgrids=true, grid style=dashed,
] 
\addplot[color=blue,mark=*, semithick, mark options={solid,fill=markercolor}]
coordinates{(0,0.00462693)(0.0125,0.007073)(0.025,0.0116375)(0.0375,0.0166377)(0.05,0.021786)(0.0625,0.0270765)(0.075,0.0326103)(0.0875,0.0385547)(0.1,0.0451154)(0.1125,0.0525108)(0.125,0.0609548)};
\addplot[color=red,mark=*, dashed, semithick, mark options={solid,fill=markercolor}]
coordinates{(0,0.00458676)(0.0125,0.0055373)(0.025,0.00775311)(0.0375,0.0105429)(0.05,0.0136791)(0.0625,0.0171546)(0.075,0.0210896)(0.0875,0.0257146)(0.1,0.0313525)(0.1125,0.0383814)(0.125,0.0471834)};

\legend{Uniform knots, Smoothed knots}
\end{semilogyaxis}
\end{tikzpicture}
}
\subfloat[Second order formulation]{
\begin{tikzpicture}
\begin{semilogyaxis}[
width=.45\textwidth,
    xlabel={Warping factor $\alpha$},   
    ylabel={$L^2$ error},
    ymin=1e-3, ymax=.1,    
    legend pos=south east, legend cell align=left, legend style={font=\tiny},	
       xmajorgrids=true,  ymajorgrids=true, grid style=dashed,
] 
\addplot[color=blue,mark=*, semithick, mark options={solid,fill=markercolor}]
coordinates{(0,0.00504237)(0.0125,0.00834838)(0.025,0.0142055)(0.0375,0.0205381)(0.05,0.0269995)(0.0625,0.0335261)(0.075,0.0400686)(0.0875,0.0466709)(0.1,0.0536183)(0.1125,0.0616203)(0.125,0.0708213)};
\addplot[color=red,mark=*, dashed, semithick, mark options={solid,fill=markercolor}]
coordinates{(0,0.0113237)(0.0125,0.0118539)(0.025,0.013336)(0.0375,0.0157542)(0.05,0.0189697)(0.0625,0.0229092)(0.075,0.027875)(0.0875,0.0339693)(0.1,0.0406825)(0.1125,0.0473484)(0.125,0.0549665)};

\legend{Uniform knots, Smoothed knots}
\end{semilogyaxis}
\end{tikzpicture}
}
\caption{Comparison of $L^2$ errors for uniform and smooth knots as a function of the warping parameter $\alpha$ using $p = 3, K= 8$ splines on a single patch.}
\label{fig:warpingconverge}
\end{figure}
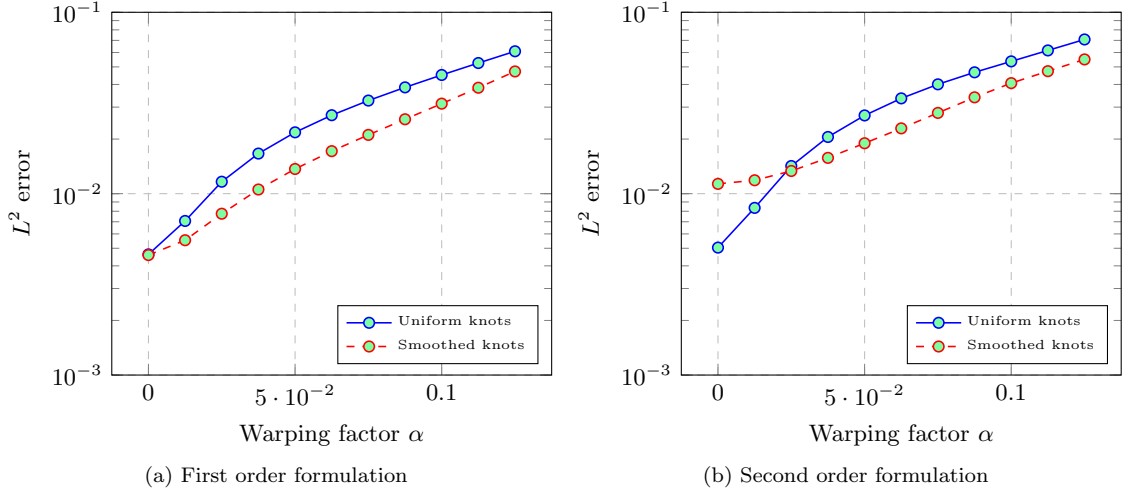

%\subsection{Spectral properties of $\bm{A}_h$ for the advection equation}
\subsection{Spectral properties of spline discretizations of the advection and first-order acoustic wave equation}

In the next few sections, we examine the spectral properties of the DG discretization matrix $\bm{A}_h$ for the \reviewerOne{one-dimensional} advection equation, such as the spectral radius (which can be used to estimate a maximum stable timestep) and discrete dispersion relations.  %Sections~\ref{sec:specadv} and \ref{sec:disp} focus on the advection equation, while Section~\ref{sec:eigerr} focuses on the second order discretization of the Laplacian.  

\subsubsection{Spectral radius of $\bm{A}_h$}
\label{sec:specadv}
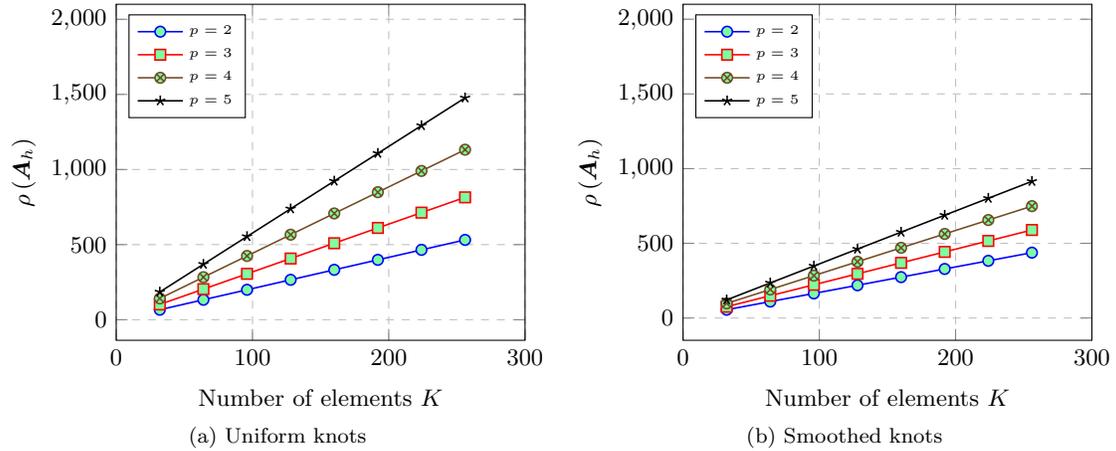
\begin{figure}
\centering
\subfloat[Uniform knots]{
\begin{tikzpicture}
\begin{axis}[
width=.425\textwidth,
    xlabel={Number of elements $K$},   
    ylabel={$\rho\LRp{\bm{A}_h}$},       
    xmin=0, xmax=300,
    ymax=0,ymax=2100,    
    legend pos=north west, legend cell align=left, legend style={font=\tiny},	
       xmajorgrids=true,  ymajorgrids=true, grid style=dashed,
] 
\addplot+[semithick, mark options={solid,fill=markercolor}]
coordinates{(32,66.4194)(64,132.839)(96,199.258)(128,265.678)(160,332.097)(192,398.517)(224,464.936)(256,531.355)};
\addplot+[semithick, mark options={solid,fill=markercolor}]
coordinates{(32,101.838)(64,203.677)(96,305.515)(128,407.354)(160,509.192)(192,611.031)(224,712.869)(256,814.707)};
\addplot+[semithick, mark options={solid,fill=markercolor}]
coordinates{(32,141.497)(64,282.993)(96,424.49)(128,565.987)(160,707.484)(192,848.98)(224,990.477)(256,1131.97)};
\addplot+[semithick, mark options={solid,fill=markercolor}]
coordinates{(32,184.614)(64,369.228)(96,553.842)(128,738.456)(160,923.07)(192,1107.68)(224,1292.3)(256,1476.91)};

\legend{$p = 2$, $p=3$, $p=4$, $p = 5$}
\end{axis}
\end{tikzpicture}
}
\subfloat[Smoothed knots]{
\begin{tikzpicture}
\begin{axis}[
width=.425\textwidth,
    xlabel={Number of elements $K$},   
    ylabel={$\rho\LRp{\bm{A}_h}$},   
    xmin=0, xmax=300,
    ymax=0,ymax=2100,    
    legend pos=north west, legend cell align=left, legend style={font=\tiny},	
       xmajorgrids=true,  ymajorgrids=true, grid style=dashed,
] 

\addplot+[semithick, mark options={solid,fill=markercolor}]
coordinates{(32,55.5039)(64,110.088)(96,164.669)(128,219.248)(160,273.827)(192,328.405)(224,382.984)(256,437.563)};
\addplot+[semithick, mark options={solid,fill=markercolor}]
coordinates{(32,75.8822)(64,149.237)(96,222.58)(128,295.92)(160,369.26)(192,442.599)(224,515.938)(256,589.276)};
\addplot+[semithick, mark options={solid,fill=markercolor}]
coordinates{(32,97.8724)(64,190.915)(96,283.926)(128,376.931)(160,469.932)(192,562.932)(224,655.931)(256,748.929)};
\addplot+[semithick, mark options={solid,fill=markercolor}]
coordinates{(32,121.298)(64,234.744)(96,348.126)(128,461.491)(160,574.85)(192,688.206)(224,801.56)(256,914.912)};

\legend{$p = 2$, $p=3$, $p=4$, $p = 5$}
\end{axis}
\end{tikzpicture}
}
\caption{Growth of $\rho\LRp{\bm{A}_h}$ for the advection equation using $\tau = .5$ and spline spaces of degree $p = 2, \ldots, 5$. }
\label{fig:rhoA}
\end{figure}

We begin by verifying that the growth of the spectral radius of the DG discretization matrix $\rho\LRp{\bm{A}_h}$ is $O(p/h)$.  Figure~\ref{fig:rhoA} shows the growth of the spectral radius for the advection equation (with weakly enforced periodic boundary conditions) using $\tau = 1/2$ and both uniform and smoothed knots.  We observe that the growth of the spectral radius matches very closely with $O(p/h)$, and constants are estimated in Table~\ref{table:rhoA}.  The use of smoothed knots reduces this growth by a factor slightly less than two.  The growth of $\rho\LRp{\bm{A}_h}$ for the acoustic wave equation is virtually identical, and we do not show it.  Similarly, the growth of $\rho\LRp{\bm{A}_h}$ for the non-dissipative case when $\tau = 0$ is very similar, and is also not shown.  

\begin{table}[!h]
\centering
\begin{tabular}{|c||c|c|c|c|}
\hline
& $p = 2$ & $p = 3$ & $p=4$ & $p=5$\\
\hhline{|=|=|=|=|=|}
Uniform knots & 2.0756 &   3.1824 &   4.4218 &   5.7692\\
\hline
Smoothed knots & 1.7092  &  2.3019 &   2.9255 &   3.5739\\
\hline
\end{tabular}
\caption{Estimated rate of growth of $\rho\LRp{\bm{A}_h}$ for the advection equation with respect to $1/h$.}
\label{table:rhoA}
\end{table}

\begin{figure}
\centering
\subfloat[Uniform knots]{
\begin{tikzpicture}
\begin{axis}[
width=.425\textwidth,
    xlabel={Number of elements $K$},   
        ylabel={$\rho\LRp{\bm{A}_h}$},       
    xmin=0, xmax=300,
    ymax=0,ymax=700,    
    legend pos=north west, legend cell align=left, legend style={font=\tiny},	
       xmajorgrids=true,  ymajorgrids=true, grid style=dashed,
] 
\addplot+[semithick, mark options={solid,fill=markercolor}]
coordinates{(32,32.9239)(64,66.339)(96,99.6489)(128,132.926)(160,166.191)(192,199.449)(224,232.703)(256,265.956)};
\addplot+[semithick, mark options={solid,fill=markercolor}]
coordinates{(32,36.0829)(64,72.8362)(96,109.357)(128,145.911)(160,182.425)(192,218.937)(224,255.45)(256,291.953)};
\addplot+[semithick, mark options={solid,fill=markercolor}]
coordinates{(32,51.525)(64,77.0896)(96,115.797)(128,154.475)(160,193.14)(192,231.798)(224,270.454)(256,309.109)};
\addplot+[semithick, mark options={solid,fill=markercolor}]
coordinates{(32,69.1513)(64,91.7232)(96,120.356)(128,160.533)(160,200.771)(192,240.95)(224,281.151)(256,321.333)};
\addplot+[semithick, mark options={solid,fill=markercolor}]
coordinates{(32,92.592)(64,121.947)(96,154.335)(128,198.112)(160,246.867)(192,296.187)(224,345.548)(256,394.912)};
\addplot+[semithick, mark options={solid,fill=markercolor}]
coordinates{(32,112.952)(64,150.924)(96,193.087)(128,254.178)(160,317.703)(192,381.243)(224,444.784)(256,508.325)};
\addplot+[semithick, mark options={solid,fill=markercolor}]
coordinates{(32,144.171)(64,189.195)(96,241.188)(128,311.773)(160,388.78)(192,466.469)(224,544.209)(256,621.953)};

%\addplot[color=black,dashed,semithick, mark options={solid,fill=markercolor}]
%coordinates{(32,32)(64,64)(96,96)(128,128)(160,160)(192,192)(224,224)(256,256)}; % O(2/h) line
\legend{$p = 2$, $p=3$, $p=4$, $p = 5$, $p = 6$, $p = 7$, $p=8$}
%\legend{$\frac{1}{h}$, $\frac{1.5}{h}$}
\end{axis}
\end{tikzpicture}
}
\subfloat[Smoothed knots]{
\begin{tikzpicture}
\begin{axis}[
width=.425\textwidth,
    xlabel={Number of elements $K$},   
        ylabel={$\rho\LRp{\bm{A}_h}$},       
    xmin=0, xmax=300,
    ymax=0,ymax=700,    
    legend pos=north west, legend cell align=left, legend style={font=\tiny},	
       xmajorgrids=true,  ymajorgrids=true, grid style=dashed,
] 

%\addplot[color=black,semithick, mark options={solid,fill=markercolor}]
%coordinates{(32,32)(64,64)(96,96)(128,128)(160,160)(192,192)(224,224)(256,256)};% [yshift=-3pt] node[below, pos=.95,color=black] {$y = 1/h$};% O(1/h) line
%\addplot[color=red,semithick, mark options={solid,fill=markercolor}]
%coordinates{(32,48)(64,96)(96,144)(128,192)(160,240)(192,288)(224,336)(256,384)}; %[yshift=3pt] node[above, pos=.95,color=black] {$y = 1.5/h$};
%\pgfplotsset{cycle list set=0}
\addplot+[semithick, mark options={solid,fill=markercolor}]
coordinates{(32,33.8023)(64,67.3347)(96,100.665)(128,133.949)(160,167.218)(192,200.479)(224,233.737)(256,266.992)}; %p = 2
\addplot+[semithick, mark options={solid,fill=markercolor}]
coordinates{(32,37.5425)(64,74.8178)(96,111.545)(128,148.129)(160,184.671)(192,221.196)(224,257.712)(256,294.223)}; % p = 3
\addplot+[semithick, mark options={solid,fill=markercolor}]
coordinates{(32,40.4264)(64,80.3236)(96,119.264)(128,158.017)(160,196.716)(192,235.393)(224,274.059)(256,312.72)}; % p = 4
\addplot+[semithick, mark options={solid,fill=markercolor}]
coordinates{(32,50.3593)(64,84.3661)(96,125.088)(128,165.445)(160,205.705)(192,245.927)(224,286.133)(256,326.328)}; % p = 5
\addplot+[semithick, mark options={solid,fill=markercolor}]
coordinates{(32,64.6167)(64,86.4858)(96,129.798)(128,173.13)(160,216.469)(192,259.811)(224,303.155)(256,346.5)}; % p = 6
\addplot+[semithick, mark options={solid,fill=markercolor}]
coordinates{(32,76.7086)(64,97.5475)(96,133.384)(128,177.856)(160,222.354)(192,266.858)(224,311.365)(256,355.876)}; %p = 7
\addplot+[semithick, mark options={solid,fill=markercolor}]
coordinates{(32,94.4519)(64,116.577)(96,136.311)(128,181.68)(160,227.101)(192,272.535)(224,317.975)(256,363.42)};  % p =8
%\addplot[color=black,dashed,semithick, mark options={solid,fill=markercolor}]
%coordinates{(32,32)(64,64)(96,96)(128,128)(160,160)(192,192)(224,224)(256,256)}; % O(2/h) line
\legend{$p = 2$, $p=3$, $p=4$, $p = 5$, $p = 6$, $p = 7$, $p=8$}
%\legend{${1}/{h}$, ${1.5}/{h}$}
\end{axis}
\end{tikzpicture}
}
\caption{Growth of $\rho\LRp{\bm{A}_h}$ for the advection equation using an exact upwind flux ($\tau = 1$) and spline spaces of degree $p = 2, \ldots, 8$. }
\label{fig:rhoAupwind}
\end{figure}
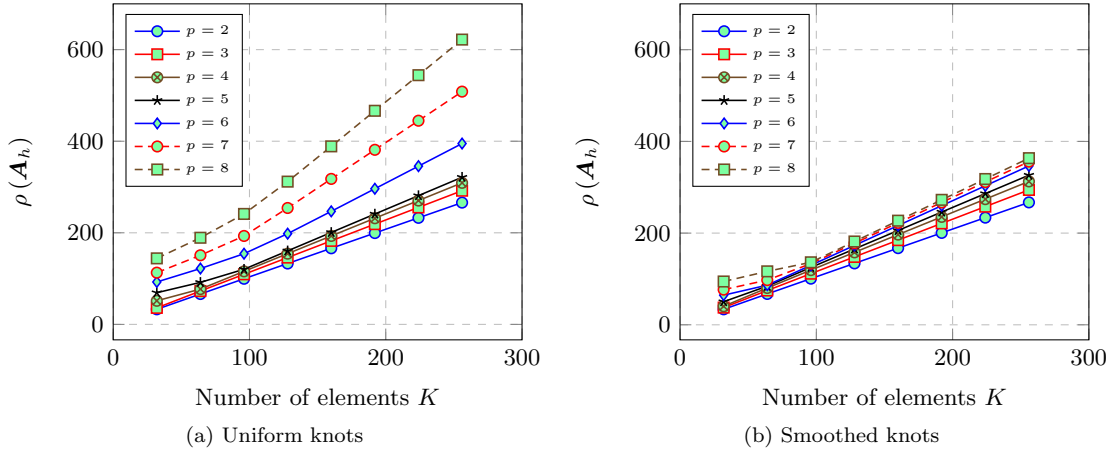

The case when $\tau = 1$ is special.  When taking a uniform wavespeed and $\tau = 1$, the penalty flux is identical to the upwind flux derived from Riemann problems \cite{hesthaven2007nodal}.  For first order formulations using an exact upwind flux, the dependence of $\rho\LRp{\bm{A}_h}$ on $h$ changes significantly for $h$ small and $p \leq 5$.  Figure~\ref{fig:rhoAupwind} shows the growth of the spectral radius with the number of elements on a single patch for uniform knot vectors and $p = 2, \ldots, 7$.  The rate of growth of $\rho\LRp{\bm{A}_h}$ with $1/h$ does not change significantly for $p \leq 5$, increasing at the expected rate of $O(p/h)$ only as $p > 5$.   %The same behavior is observed when uniformly refining multiple spline patches.  
If smoothed knots are used in lieu of uniform knot vectors when $\tau = 1$, we observe that $\rho\LRp{\bm{A}_h}$ is significantly smaller in the pre-asymptotic range of $1/h$ and slightly larger in the asymptotic range of $1/h \rightarrow \infty$, though the rate of growth of the spectral radius is unchanged for $1/h$ sufficiently large.  When using smoothed knots and exact upwind fluxes, we also observe that the rate of growth $\rho\LRp{\bm{A}_h}$ with $1/h$ is slower than the estimated $O(p/h)$ rate up to degree 8.

This phenomena appears to be restricted to the specific case of \reviewerOne{$\tau=1$, where the penalty flux coincides with the upwind flux}, and the growth of $\rho\LRp{\bm{A}_h}$ returns to $O(p/h)$ when the penalty parameter deviates by more than $O(1\%)$ from $\tau = 1$.  \reviewerOne{However, numerical experiments also indicate that this phenomena is present for the acoustic wave equation, and persists when using multiple patches.  
%even discontinuous advection $a(x)$, provided that the upwind numerical flux is used
%\[
%\avg{a \phi} - \frac{1}{2}\LRp{|a|^+\phi^+ - |a|^-\phi^-}
%\]
Future work will study whether this phenomena is present for more complex settings (such as curvilinear meshes and discontinuous wavespeeds).}

\begin{figure}
\centering
\subfloat[Uniform knots, $\tau = 0,.5,1$]{\includegraphics[height=.3\textheight]{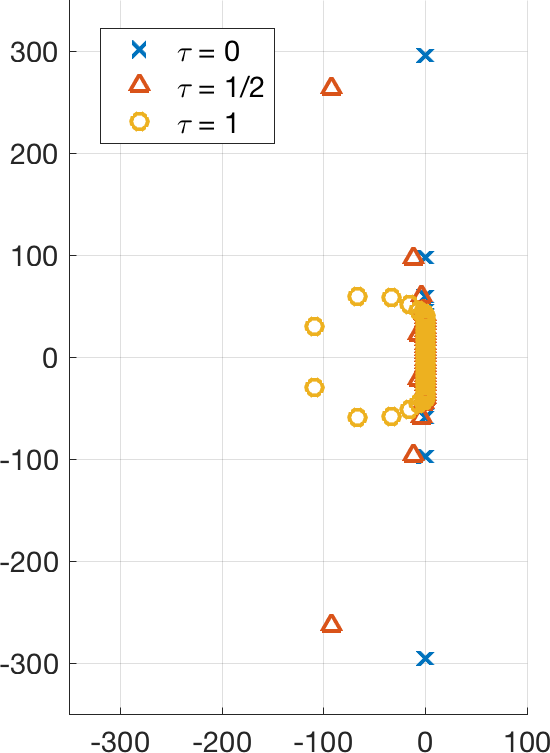}}
\hspace{2em}
\subfloat[Uniform vs smoothed knots, $\tau = 1$]{\includegraphics[height=.3\textheight]{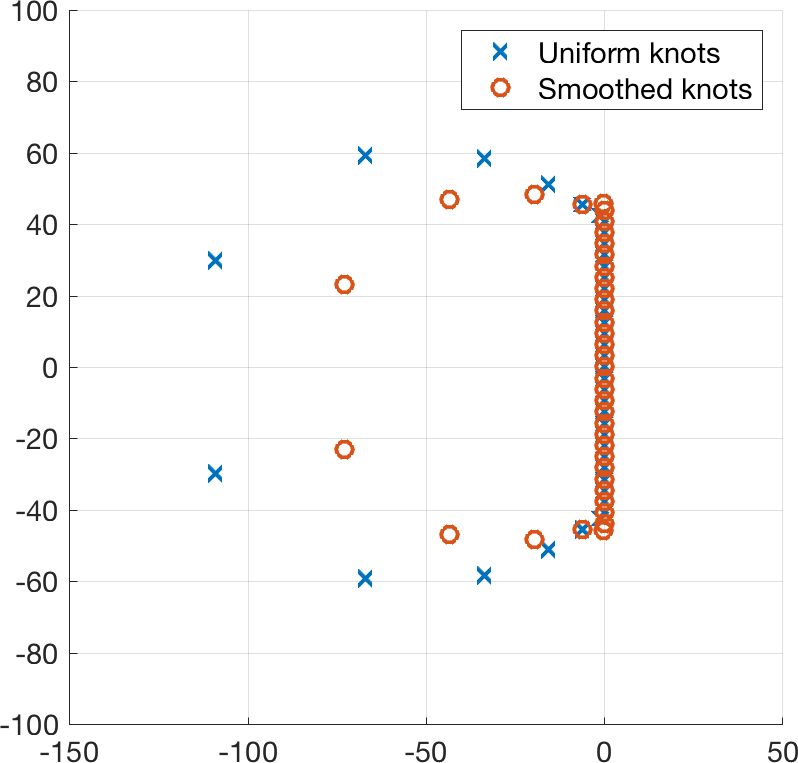}}
\caption{Spectra for the advection equation using a spline space with $p=7$, $K = 32$. }
\label{fig:spectra}
\end{figure}

\reviewerOne{A closer examination of the spectra of $\bm{A}_h$ offers some insight into the growth of $\rho(\bm{A}_h)$ and its dependence on $\tau$ and knot smoothing.  We note that the distribution of the spectra can also impact the maximum stable timestep, which depends on the choice of time-stepping method and its corresponding region of stability.  
Figure~\ref{fig:spectra} shows various spectra for the advection equation for a $p=7$, $K=32$ spline space.  It can be observed that $\rho(\bm{A}_h)$ is much larger for $\tau = 0, .5$ than for $\tau = 1$ due to the presence of two outlying eigenvalues with large imaginary part.  For $\tau = 1$, these extremal eigenvalues fall into a tight semi-circular distribution in the left half plane, resulting in a significant reduction in the spectral radius.  When knot smoothing is applied, the eigenvalues contract into a semi-circular distribution of slightly smaller radius.}

\subsubsection{Discrete dispersion relations }
\label{sec:disp}
\begin{figure}
\centering
\subfloat[Uniform knots]{
\includegraphics[width=.45\textwidth]{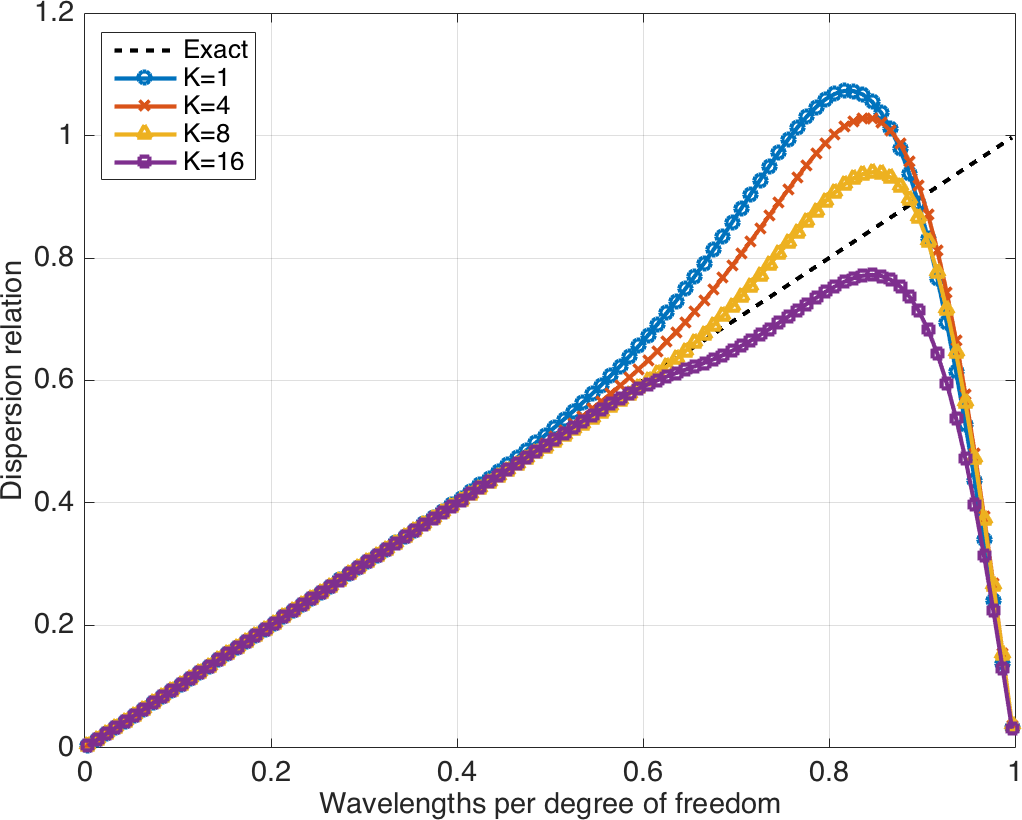}
\label{subfig:disp1}
}
\subfloat[Smoothed knots]{
\includegraphics[width=.45\textwidth]{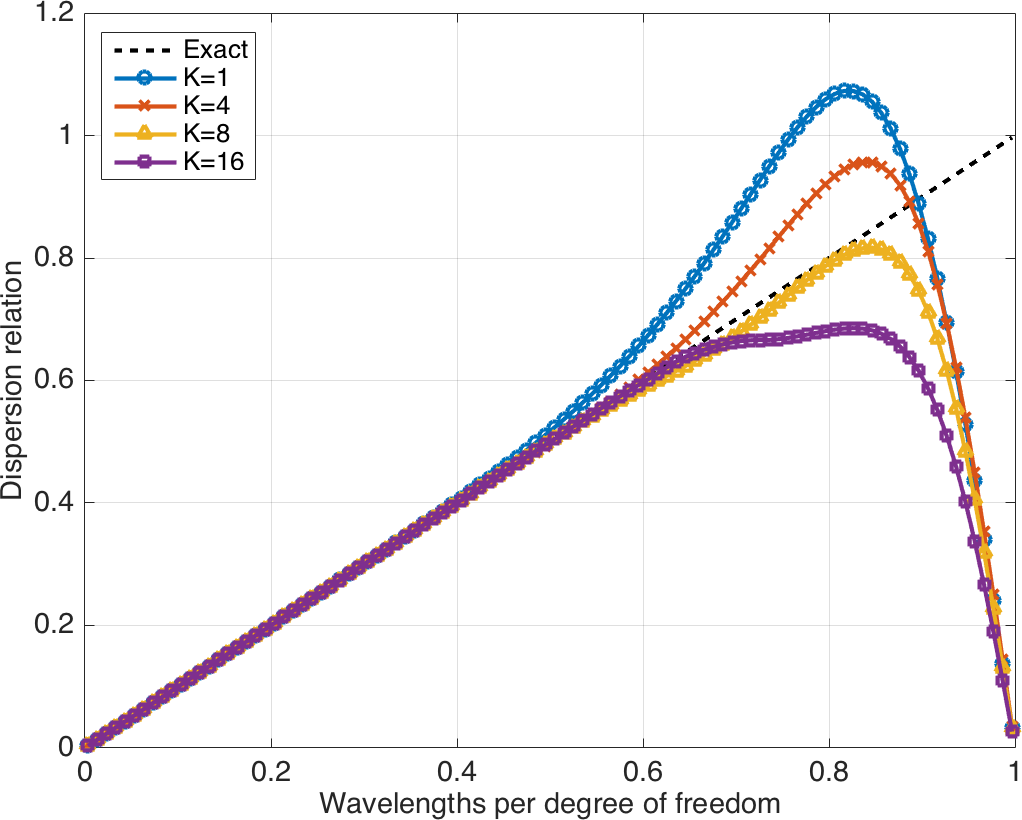}
\label{subfig:disp2}
}\\
\subfloat[Dispersion errors (uniform knots)]{
\includegraphics[width=.45\textwidth]{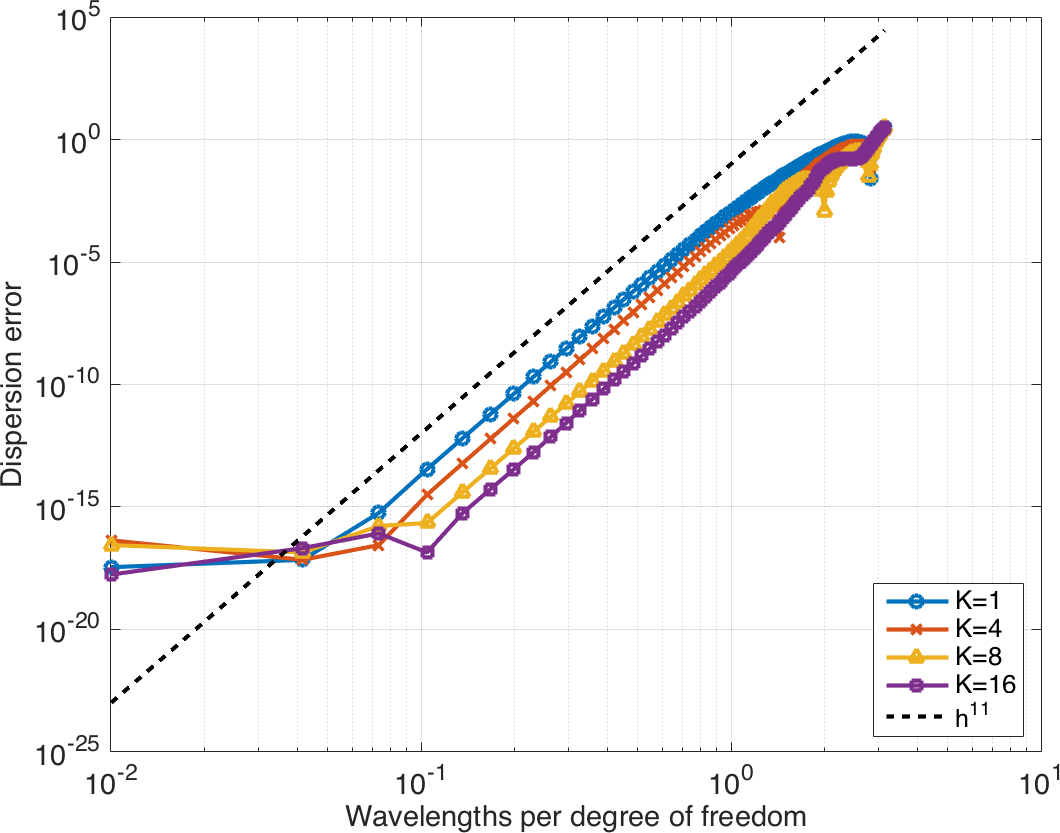}
}
\subfloat[Dispersion errors (smoothed knots)]{
\includegraphics[width=.45\textwidth]{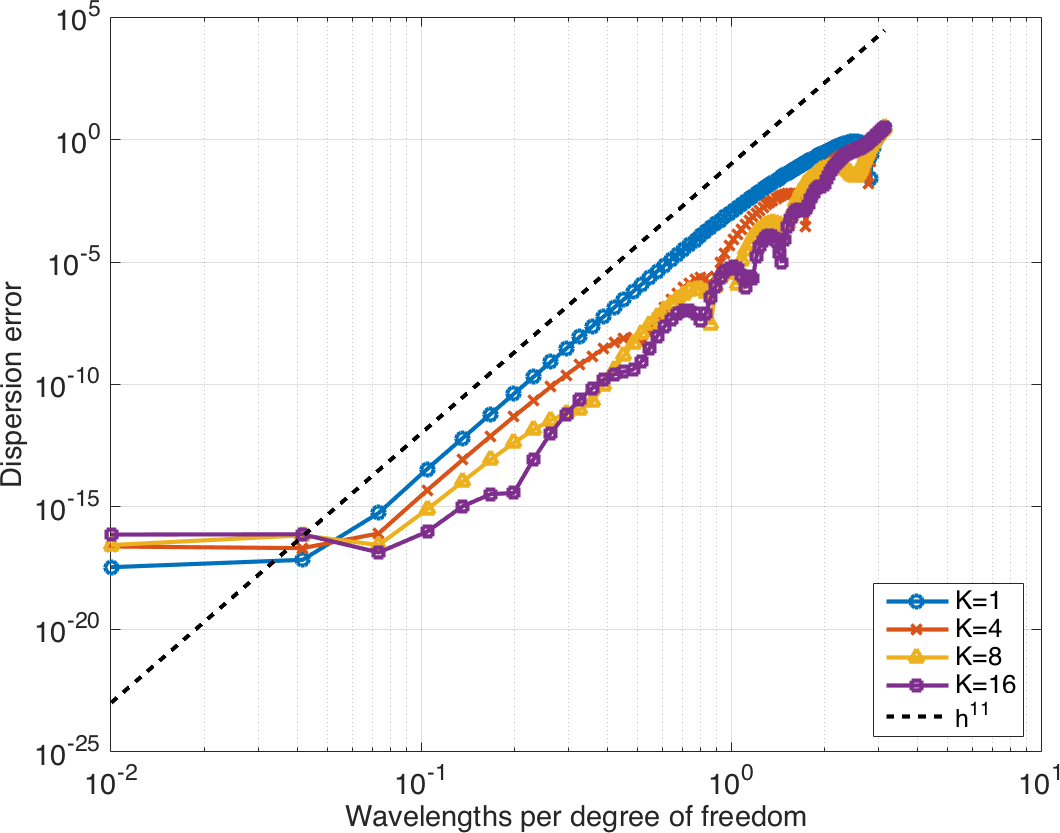}
}
\caption{Numerical dispersion relations and dispersion error $\LRb{{\rm Re}(\omega) - {\rm Re}(\omega_h)}$ for spline spaces with $p = 4, K = 1, 4, 8, 16$ and $\tau = 1$.  The exact dispersion relation is given as a black dotted line in Figure~\ref{subfig:disp1}, \reviewerOne{ and the theoretical rate of convergence of the dispersion error (for $p=4$, this is $O(h^{11})$) is plotted as a black dotted line in Figure~\ref{subfig:disp2}}.  }
\label{fig:eigerr1}
\end{figure}

We compute numerical dissipation and dispersion relations for the periodic advection equation using an upwind flux ($\tau = 1$).\footnote{\reviewerOne{For $\tau = 0$ (which corresponds to a central flux), dispersion errors are similar to the $\tau =1$ case for a small number of wavelengths per degree of freedom.  For a larger number of wavelengths per degree of freedom, dispersion errors are larger for $\tau = 0$ than for $\tau = 1$.  The qualitative behavior of the discrete dispersion relation is similar to that of polynomial DG discretizations \cite{ainsworth2004dispersive}.}}  Assuming a discrete solution of the form $u_h(x,t) = e^{i (kx - \omega_h t)}$, we seek a relation between the discrete frequency $k$ and wavenumber $\omega$.  Inserting this ansatz in the semi-discrete variational formulation for advection over a single patch yields a generalized eigenvalue problem, which can be solved for $\omega_h$ \cite{hesthaven2007nodal}.  While the exact relation between frequency and wavenumber is $\omega = k$, the discrete wavenumber $\omega_h$ will differ due to discretization errors.  Figure~\ref{fig:eigerr1} compares ${\rm Re}(\omega_h)$ (which corresponds to numerically introduced dispersion) with the exact relation $\omega = k$.  

We note that increasing $p$ improves the rate of convergence of the dispersion error with the number of wavelengths per degree of freedom, while increasing $K$ mainly shifts the dispersion error curve \reviewerOne{by a constant number of wavelengths per degree of freedom}.  The dispersion error converges at a rate of $(2p+3)$ for $p$ even ($(2p+2)$ for $p$ odd) with respect to the number of wavelengths per degree of freedom \cite{ainsworth2004dispersive}.  As expected, since spline spaces contain polynomials of degree $p$, we observe this rate for both uniform and smoothed knots.  Both uniform and smoothed knot vectors result in comparable dispersion errors, though smoothed knot vectors result in slightly decreased dispersion errors at a higher number of wavelengths per degree of freedom.  We note that these observations are in line with the approximation results of Figure~\ref{fig:ppw}.  

%\subsection{Eigenvalue and eigenvector errors for the second order Laplacian}
\subsection{Spectral properties of spline discretizations of the second-order acoustic wave equation}
\label{sec:eigerr}

\begin{figure}
\centering
\subfloat[Uniform knots]{
\begin{tikzpicture}
\begin{axis}[
width=.425\textwidth,
    xlabel={Number of elements $K$},   
    ylabel={$\sqrt{\rho\LRp{\bm{A}_h}}$},       
    xmin=0, xmax=300,
    ymax=0,ymax=2100,    
    legend pos=north west, legend cell align=left, legend style={font=\tiny},	
       xmajorgrids=true,  ymajorgrids=true, grid style=dashed,
] 
\addplot+[semithick, mark options={solid,fill=markercolor}]
coordinates{(32,83.0595)(64,171.538)(96,259.969)(128,348.389)(160,436.804)(192,525.216)(224,613.628)(256,702.038)};

\addplot+[semithick, mark options={solid,fill=markercolor}]
coordinates{(32,129.151)(64,267.997)(96,406.744)(128,545.468)(160,684.182)(192,822.892)(224,961.599)(256,1100.3)};
\addplot+[semithick, mark options={solid,fill=markercolor}]
coordinates{(32,179.72)(64,374.676)(96,569.456)(128,764.195)(160,958.918)(192,1153.63)(224,1348.34)(256,1543.05)};
\addplot+[semithick, mark options={solid,fill=markercolor}]
coordinates{(32,233.735)(64,489.525)(96,745.037)(128,1000.48)(160,1255.9)(192,1511.31)(224,1766.7)(256,2022.1)};

\legend{$p = 2$, $p=3$, $p=4$, $p = 5$}
\end{axis}
\end{tikzpicture}
}
\subfloat[Smoothed knots]{
\begin{tikzpicture}
\begin{axis}[
width=.425\textwidth,
    xlabel={Number of elements $K$},   
    ylabel={$\sqrt{\rho\LRp{\bm{A}_h}}$},   
    xmin=0, xmax=300,
    ymax=0,ymax=2100,    
    legend pos=north west, legend cell align=left, legend style={font=\tiny},	
       xmajorgrids=true,  ymajorgrids=true, grid style=dashed,
] 

\addplot+[semithick, mark options={solid,fill=markercolor}]
coordinates{(32,66.9105)(64,138.075)(96,209.18)(128,280.27)(160,351.355)(192,422.436)(224,493.516)(256,564.595)};
\addplot+[semithick, mark options={solid,fill=markercolor}]
coordinates{(32,91.7975)(64,190.07)(96,288.195)(128,386.284)(160,484.359)(192,582.426)(224,680.489)(256,778.55)};
\addplot+[semithick, mark options={solid,fill=markercolor}]
coordinates{(32,117.432)(64,243.924)(96,370.121)(128,496.246)(160,622.341)(192,748.422)(224,874.495)(256,1000.56)};
\addplot+[semithick, mark options={solid,fill=markercolor}]
coordinates{(32,143.557)(64,299.146)(96,454.218)(128,609.163)(160,764.056)(192,918.925)(224,1073.78)(256,1228.62)};

\legend{$p = 2$, $p=3$, $p=4$, $p = 5$}
\end{axis}
\end{tikzpicture}
}
\caption{Growth of $\sqrt{\rho\LRp{\bm{A}_h}}$ for the second order wave equation using spline spaces of degree $p = 2, \ldots, 5$. }
\label{fig:rhoA2}
\end{figure}
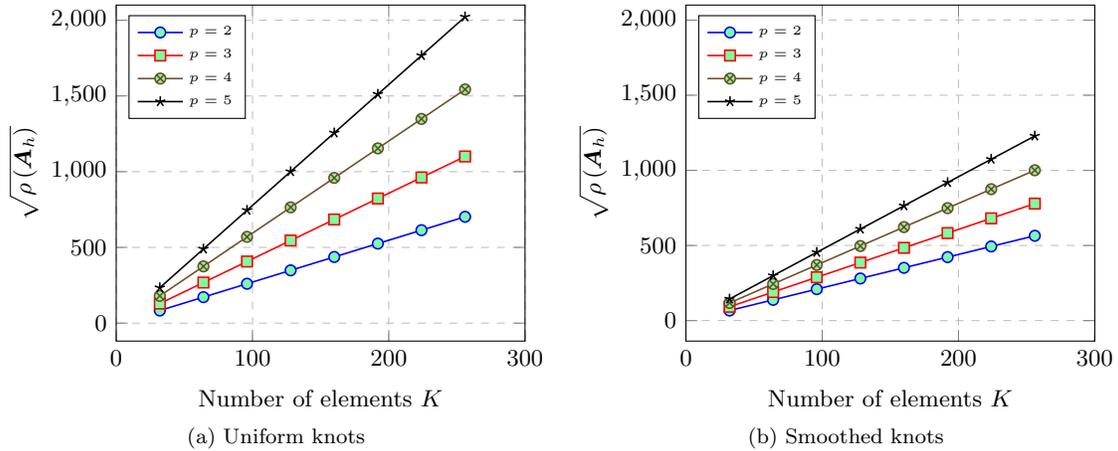

\reviewerTwo{We now turn our attention to the second order formulation of the acoustic wave equation.  
Recall that the stable timestep restriction for the second order formulation is estimated by the square root of the spectral radius $\sqrt{\rho\LRp{\bm{A}_h}}$.  Figure~\ref{fig:rhoA2} shows the value of $\sqrt{\rho\LRp{\bm{A}_h}}$ for various $p$ and $K$, and it can be seen that the values of $\sqrt{\rho\LRp{\bm{A}_h}}$ for the second order formulation are slightly larger than the values of $\rho\LRp{\bm{A}_h}$ observed for the first order formulation in Figure~\ref{fig:rhoA}.  Table~\ref{table:rhoA2} also plots the estimated slope of the growth of $\sqrt{\rho\LRp{\bm{A}_h}}$ with respect to $1/h$.  The slopes are roughly $O(p)$, though the values of the estimated slopes are also slightly higher than those shown in Table~\ref{table:rhoA} for the first order formulation.
}

\begin{table}[!h]
\centering
\begin{tabular}{|c||c|c|c|c|}
\hline
& $p = 2$ & $p = 3$ & $p=4$ & $p=5$\\
\hhline{|=|=|=|=|=|}
Uniform knots & 2.7632  &  4.3353  &  6.0859  &  7.9831\\
\hline
Smoothed knots & 2.2217     &3.0655     & 3.9419    &4.8429\\
\hline
\end{tabular}
\caption{Estimated rate of growth of $\sqrt{\rho\LRp{\bm{A}_h}}$ for the second order wave equation with respect to $1/h$.}
\label{table:rhoA2}
\end{table}

\reviewerTwo{We next examine eigenvalue and eigenvector errors for spline discretizations of the following one-dimensional generalized eigenproblem:
\[
-\frac{d^2p_k}{dx^2} = \lambda_k p_k
\]
over the unit interval $(0,1)$ subject to the boundary conditions $p_k(0) = p_k(1) = 0$.  The above eigenproblem admits an infinite number of eigenvalue and eigenvector solutions $\left(\lambda_k,p_k(x) \right) = \left(k^2\pi^2, \sin(k\pi x) \right)$ for $k \in \mathbb{N}$.  A finite number of discrete eigenvalues and eigenvectors $\left(\lambda_{h,k}, u_{h,k}(x)\right)$ with $0 < \lambda_{h,1} \leq \lambda_{h,2} \leq \ldots \leq \lambda_{h,N_{\rm dofs}}$ are attained from a standard Galerkin discretization on a single patch.  Homogeneous Dirichlet boundary conditions are imposed in a strong fashion as to allow for a more direct comparison to earlier works in the literature \cite{hughes2008duality, hughes2014finite}.  We have also examined the impact of obtaining boundary conditions in a weak manner and found that this did not significantly affect solution quality.  In prior work, it has been established that the error associated with a standard Galerkin discretization of the second order wave equation is directly tied to the eigenvalue and eigenvector errors associated with the above generalized eigenproblem \cite{hughes2014finite}, inspiring us to examine the errors associated with each of the eigenvalues and eigenvectors for a given discretization.}

%We now examine eigenvalue and eigenvector errors for the second order discretization of the Laplacian.  We compute discrete eigenvalues and eigenvectors $\lambda_{h,k}, u_{h,k}(x)$ from a standard Galerkin discretization on a single patch, where homogeneous Dirichlet boundary conditions are imposed in a strong fashion.  We note that the choice of strong imposition is to allow for a more direct comparison to earlier works in the literature \cite{hughes2008duality, hughes2014finite}, and that boundary conditions can also be imposed in a weak manner without significantly affecting these results.  

\begin{figure}
\centering
\subfloat[Eigenvalue errors]{
\includegraphics[width=.45\textwidth]{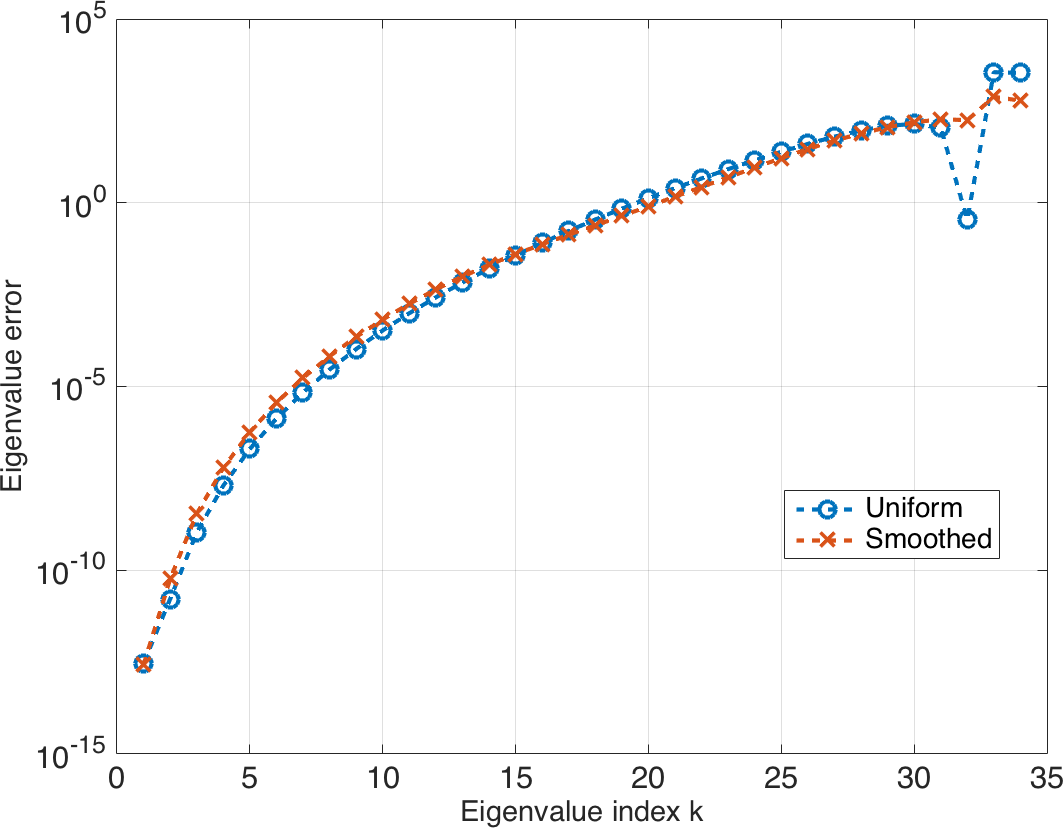}
}
\subfloat[Eigenvector errors]{
\includegraphics[width=.45\textwidth]{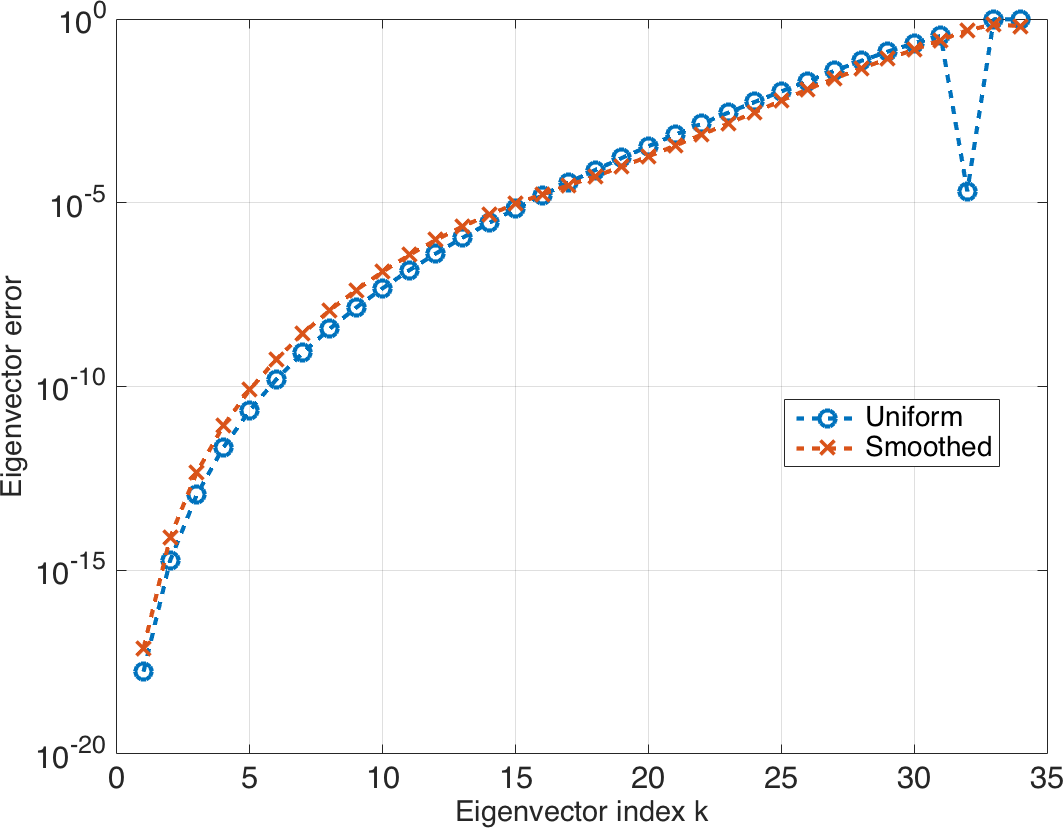}
}
\caption{Eigenvalue errors $\LRb{\lambda_k-\lambda_{h,k}}$ and eigenvector errors $\nor{p_k-p_{h,k}}_{L^2}^2$ for an $p = 4, K = 32$ spline space using uniform and smoothed knots.  }
\label{fig:eigerr2}
\end{figure}

Figure~\ref{fig:eigerr2} shows errors in eigenvalues $\LRb{\lambda_k-\lambda_{h,k}}$ { and eigenfunctions $\LRb{u_k-u_{h,k}}$} using splines with both uniform and smoothed knots { for the specific case of $p = 4$ and $K = 32$}.  We observe that smoothed knots are slightly less accurate for low $k$ and slightly more accurate for high $k$.  { Note furthermore that the last two eigenvalues are poorly approximated using uniform knots.  This is due to the fact that the last two discrete eigenvalues correspond to outlier frequencies \cite{hughes2008duality}.  The last two eigenvalues are slightly better approximated and smaller in magnitude for smoothed knots, and are the reason why the maximum stable timestep is larger with smoothed knots than with uniform knots.  
}

\subsection{A three-dimensional problem on a multi-patch geometry}

\begin{figure}
\centering
\subfloat[Multi-patch model]{\includegraphics[width=.35\textwidth]{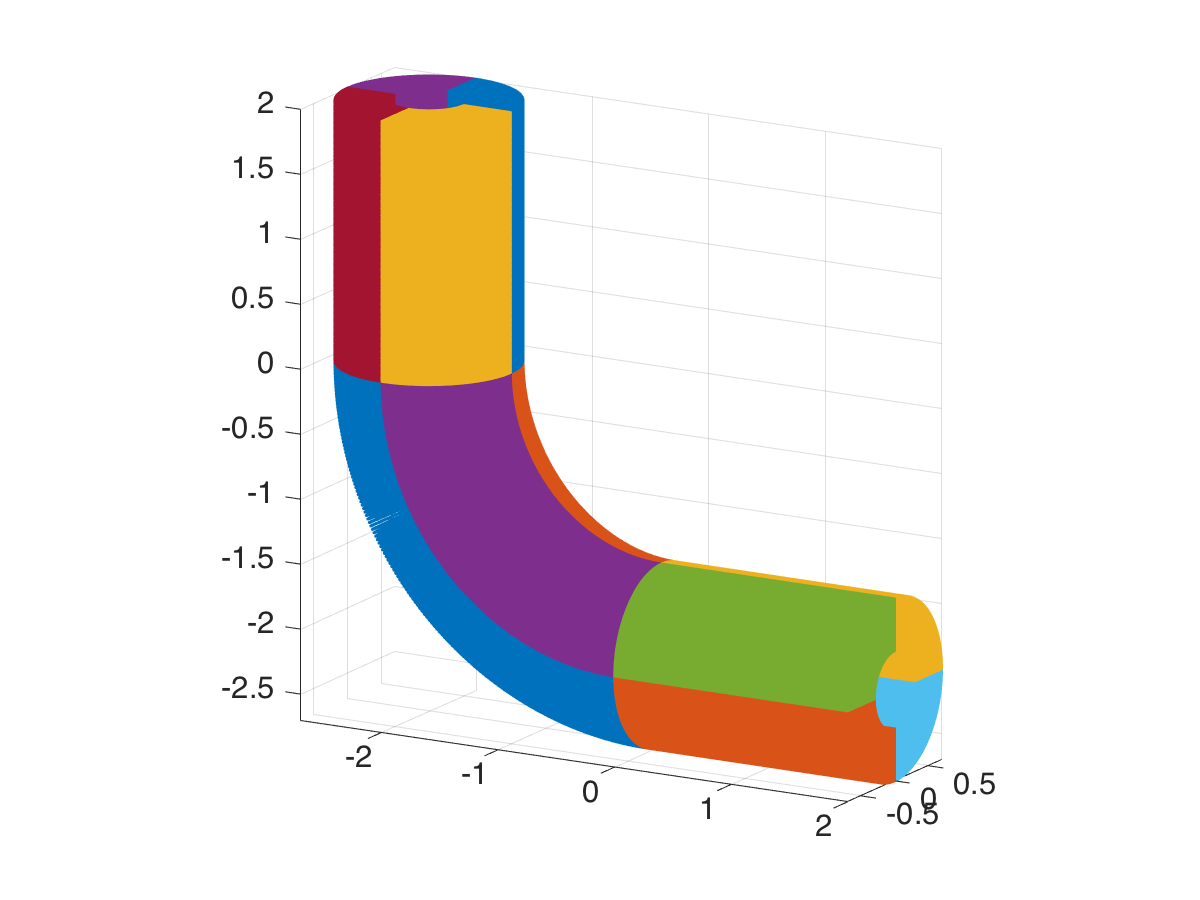}}
\subfloat[Pressure]{\includegraphics[width=.33\textwidth]{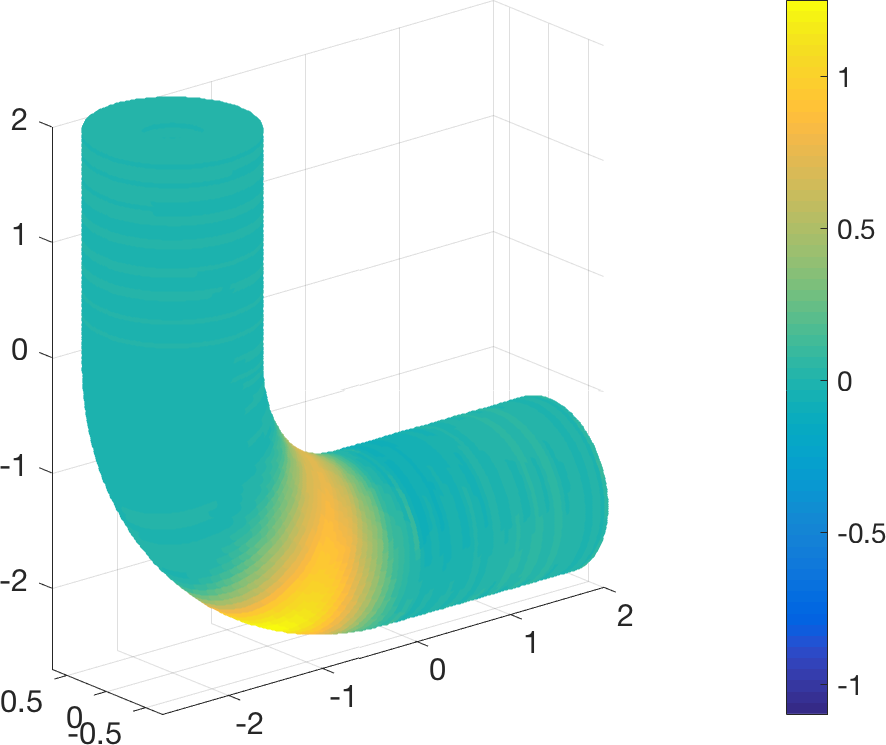}}
\subfloat[Pressure (cut view)]{\includegraphics[width=.33\textwidth]{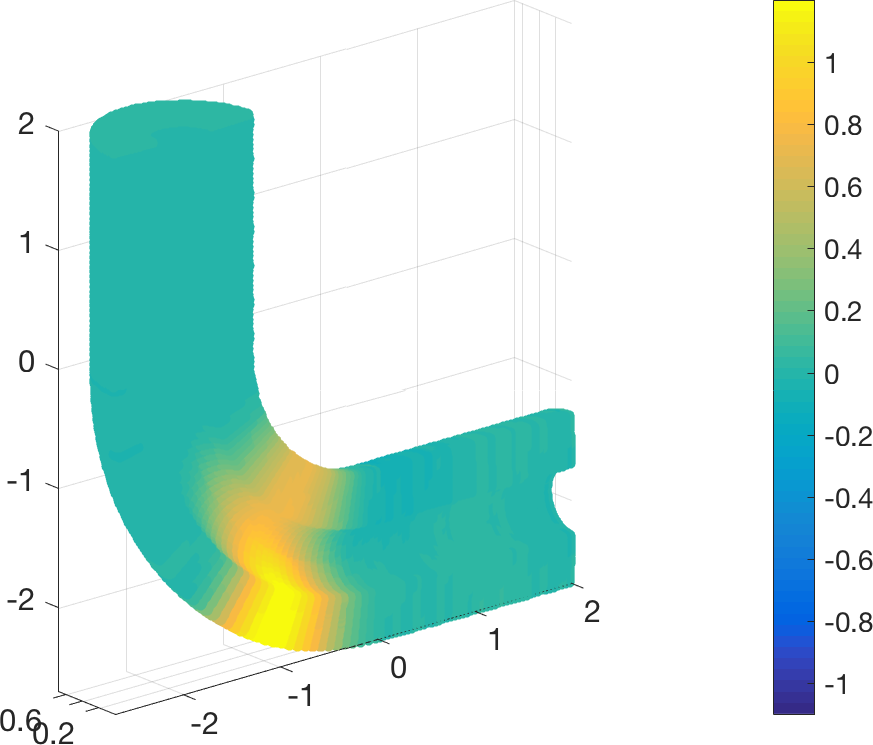}}
\caption{A multi-patch model of a 3D pipe elbow, and pressure solution at time $T = 4$. }
\label{fig:pipe3D}
\end{figure}

Finally, we show solutions to a model wave propagation problem using a non-trivial geometry.  Figure~\ref{fig:pipe3D} shows a multi-patch discretization of a three-dimensional pipe elbow, constructed using 12 different patches, as well as a first order DG simulation of acoustic wave propagation through the pipe.  We assume $c = 1$ and a zero initial condition, and impose a pulse velocity boundary condition at $x=2$
\[
u(\bm{x},0) = \begin{cases}
1-\cos(2\pi t /t_0), &t < t_0\\
0, &\text{ otherwise}
\end{cases}
\]
where $u$ is the velocity in the $x$-direction, and we take $t_0 = 2$.  We impose zero Neumann boundary conditions at all other boundary faces.  
%\note{Talk about pipe geometry and setup.  Show plots.}

We discretize each patch using a degree $p=6$ spline basis (using a smoothed knot vector) and $K = 16$ elements in each coordinate direction.  We note that, due to the stretched nature of some of the patches, additional efficiency might be gained by using an anisotropic discretization with a varying number of elements in each coordinate direction, which we will explore in a future manuscript.  

\section{Conclusions}
\label{sec:conc}

This work presents a strategy for applying NURBS-based finite element discretizations to transient hyperbolic problems using explicit time-stepping.  We utilize a multi-patch DG discretization, and apply a weight-adjusted approximation to the mass matrix inverse over each patch.  Additionally, the approximation is efficient to invert, involving only one-dimensional operations due to the tensor product structure in approximation of the mass matrix inverse.  The resulting methods are energy stable and high order accurate under assumptions on the regularity of the geometric mapping, and numerical experiments show that the approximate weight-adjusted mass matrix inverse delivers $L^2$ errors which are virtually identical to those using the full mass matrix inverse.   

We also investigate the timestep restriction associated with NURBS-based discretizations, and show numerically that (for patches where $h$ is sufficiently small with respect to $p$) the maximum stable timestep decreases as $O(h/p)$ instead of the $O(h/p^2)$ associated with $C^0$ and DG finite element discretizations.  Finally, we investigate the use of smoothed knot vectors (which approximate $n$-width optimal knot vectors) in spline discretizations.  Numerical experiments show that spline spaces under smoothed knot vectors are slightly more accurate than splines with uniform knot vectors for approximations involving oscillatory functions or warped geometric mappings.  We remark that, while the use of smoothed knots does not result in drastic decreases in $L^2$ error for high frequencies or curved mappings, their improved accuracy and larger stable time-step make their use attractive for time-domain solvers based on explicit time integration.  

We note that several areas related to this work remain to be explored.  For example, it is not immediately clear when splines can achieve a computational advantage over very high order polynomials.  Numerical experiments suggest that in certain situations (for example, on curved domains or approximating highly oscillatory functions) splines yield lower $L^2$ errors than polynomials for a similar number of degrees of freedom.  Additionally, for the same resolution, splines can take larger time steps than an equivalent finite element $p$-method, making it possible to achieve the same level of error with a reduced number of time steps, assuming a sufficiently accurate time stepping scheme.  However, these advantages must be balanced with the fact that, in order to integrate spline spaces \reviewerTwo{to} sufficient accuracy, the number of quadrature points per degree of freedom is larger for splines than for polynomials.  This additional cost can be significantly reduced through the use of optimal and near-optimal spline quadrature rules \cite{hughes2010efficient, auricchio2012simple} or by specialized techniques for the assembly and application of spline finite element matrices \cite{calabro2016fast}.  Finally, we will explore ways to fully remove any inverse scaling of the maximum stable timestep by $p$.  It may be advantageous to combine the nonlinear mapping used in \cite{hughes2008duality} with optimal or smoothed knot vectors.  Additionally, preliminary results in this work suggest that, when using smoothed knot vectors, a first order formulation combined with an appropriately chosen dissipative flux result in discretization matrices whose spectral radius is nearly independent of $p$ for $K$ sufficiently large.   

The conformity requirements of multi-patch DG methods also remain to be explored.  In this work, all solution components are approximated using the same spline space; however, it is well-known that such a choice of approximation space can introduce problems such as locking or spurious modes to finite element discretizations of problems such as Maxwell's equations \cite{buffa2010isogeometric, ratnani2012arbitrary}.  DG methods sidestep these issues by using numerical fluxes to penalize or dissipate away non-conforming components of the solution \cite{Warburton20063205, chan2016short}.  However, for multi-patch DG methods, it may also be necessary to utilize locally conforming approximation spaces within each patch, which we will consider in future work.

\section{Acknowledgments}

The authors thank Joseph Benzaken help in constructing the three-dimensional pipe elbow geometry.  Jesse Chan is supported by the National Science Foundation under awards DMS-1719818 and DMS-1712639.  John A. Evans was partially supported by the Air Force Office of Scientific Research under Grant No. FA9550-14-1-0113.  

\bibliographystyle{unsrt}
\bibliography{dgpenalty}

\end{document}